\documentclass{article}
\usepackage{graphicx} 
\usepackage{amssymb}
\usepackage{enumitem}
\usepackage{hyperref}
\usepackage{amsmath}
\usepackage{parskip}
\usepackage{tikz-cd}

\usepackage{amsthm}
\usepackage[titletoc]{appendix}
\usepackage{amsfonts}
\usepackage{mathrsfs}
\usepackage{afterpage}
\usepackage{float}
\usepackage{caption}
\usepackage{multirow}
\usepackage{graphicx}
\usepackage{color}
\usepackage[margin=1.0in]{geometry}
\newcommand*{\SHom}{\mathscr{H}\kern -.5pt om}
\newcommand*{\SExt}{\mathscr{E}\kern -.5pt xt}

\title{On operadic open-closed maps in characteristic $p$}
\author{Zihong Chen}
\date{}

\begin{document}
\theoremstyle{definition}
\newtheorem{mydef}{Definition}[section]
\numberwithin{mydef}{section}
\newtheorem{rmk}[mydef]{Remark}
\newtheorem{conj}[mydef]{Conjecture}
\theoremstyle{plain}
\newtheorem{cor}[mydef]{Corollary}
\newtheorem{lemma}[mydef]{Lemma}
\newtheorem{thm}[mydef]{Theorem}
\newtheorem{prop}[mydef]{Proposition}
\maketitle

\begin{abstract}
Consider a closed monotone symplectic manifold $(M,\omega)$. \cite{Gan2} constructed a cyclic open-closed map, which goes from the cyclic homology of the Fukaya category of $M$ to the $S^1$-equivariant quantum cohomology of $M$. In this paper, we show that with mod $p$ coefficients, Ganatra's cyclic open-closed map is compatible with a certain $\mathbb{Z}/p$-equivariant open-closed map under the natural $\mathbb{Z}/p$-Gysin type comparison map for Hochschild homology. This is a key technical result that appears in the author's subsequent work \cite{Che} concerning Fukaya categorical approaches to (equivariant) Gromov-Witten theory with mod $p$ coefficients, where we used it to prove the unramified exponential type conjecture, under a certain generation assumption, via reduction mod $p$ methods. Along with the proof of the technical result, this paper gives a new homotopy theoretic framework for studying open-closed maps in symplectic topology. The main insights of this paper are: 1) a $\mathbb{Z}/p$-Gysin comparison result for ($\mathcal{A}_{\infty}$-) cyclic objects, 2) new construction of the open-closed map using the language of cyclic objects and \cite{AGV}'s operadic Floer theory, and 3) comparisons of the new constructions with its classical counterparts.
\end{abstract}

\renewcommand{\theequation}{1.\arabic{equation}}
\setcounter{equation}{0}

\section{Introduction}
1.1. \textbf{Motivation and background}. Let $(M,\omega)$ be a closed monotone symplectic manifold. There are two types of algebraic invariants that one can build out of the symplectic topology of $M$. One is called the \emph{open-string theory}, which studies the Floer homology of Lagrangian submanifolds of $M$, packaged into an $\mathcal{A}_{\infty}$-category called the (monotone) Fukaya category \cite{Oh},\cite{She}. The other is called the \emph{closed-string theory}, which studies operations on the quantum cohomology of $M$ \cite{MS}, i.e. the singular cohomology of $M$ equipped with a deformed cup product that encodes its genus $0$ Gromov-Witten invariants.  \par\indent
The open and closed-string theory of $M$ are closely related, and the following instance is of specific interest to us. Let $CC^{S^1}$ denote the negative cyclic homology chain complex of an $\mathcal{A}_{\infty}$-category, cf. (2.15). \cite{Gan2} defined the (negative) \emph{cyclic open-closed map}, which is a chain map
\begin{equation}
OC^{S^1}: CC^{S^1}(\mathrm{Fuk}(M)_{\lambda})\rightarrow QH(M)^{S^1}=QH(M)[[t]]
\end{equation}
of degree $\frac{1}{2}\dim_{\mathbb{R}}M$, where $QH$ denotes a chain model computing the quantum cohomology of $M$. Here, $t$ is the $S^1$-equivariant formal variable of degree $2$, corresponding to a degree $2$ generator of $H^*(BS^1)$. The cyclic open-closed map plays an important role in the study of noncommutative Hodge structures in symplectic topology, as it relates the variation of semi-infinite Hodge structures
on $S^1$-equivariant quantum cohomology of $M$ with that on the negative cyclic homology of $\mathrm{Fuk}(M)_{\lambda}$, cf. \cite{GPS}. As a special instance, one has the following conjecture of \cite{GPS}, which was proved by \cite[Theorem 1.7]{Hug} and \cite[Theorem 6.3.5]{PS} separately in different settings.
\begin{thm}
At the level of cohomology, $OC^{S^1}$ intertwines the Getzler-Gauss-Manin connection
on $CC^{S^1}(\mathrm{Fuk}(M)_{\lambda})$ with the quantum connection on $QH(M)^{S^1}$.
\end{thm}
Theorem 1.1 can be used as a powerful tool to compute the cyclic open-closed map in characteristic $0$, since in this situation an algebraic differential equation imposes strong constraints on its solutions, cf. \cite[section 6]{Hug}, \cite{Che}. \par\indent
This paper is devoted to the study of the cyclic open-closed map and its $\mathbb{Z}/p$-equivariant analogue in the context of mod $p$ coefficients.
The latter, which we call the \emph{$\mathbb{Z}/p$-equivariant open-closed map}, is used by the author in an upcoming work \cite{Che} to give a Fukaya categorical interpretation of the Quantum Steenrod operations, cf. Theorem 1.3 in loc.cit. These operations, originally due to Fukaya \cite{Fuk}, come from quantization of the classical Steenrod operations via mod $p$ counts of certain genus $0$ (equivariant) Gromov-Witten invariants. Its structural properties were systematically studied in
in \cite{Wil}, and it has since seen various applications to Hamiltonian dynamic \cite{Shel}, arithmetic mirror symmetry \cite{Sei2} and the quantum differential equation in characteristic $p$ \cite{SW}, \cite{Lee}.\par\indent
One important consequence of the Fukaya categorical approach in \cite{Che} is that
Quantum Steenrod operations, which are intrinsically defined over a field of characteristic $p$, preserve certain structures imposed by the quantum differential equation in characteristic $0$, cf. \cite[Corollary 1.10]{Che}. An interesting application of this result is a proof of the unramified exponential type conjecture of \cite{KKP} and \cite{GGI}, under the assumption that $X$ satisfies Abouzaid's generation criterion over $\overline{\mathbb{Q}}$, using a reduction mod $p$ argument, cf. \cite[Theorem 1.2]{Che}.   \par\indent
The main objective of the current paper is to complete the key technical result that will be used in \cite{Che} to enable this characteristic $0$/$p$ interaction. Namely, the $\mathbb{Z}/p$-equivariant open-closed map defined in \cite{Che} and the cyclic
open-closed map of \cite{Gan2} are compatible under a $\mathbb{Z}/p$-Gysin comparison map for Hochschild homology. We now describe the setup of the main theorem. \par\indent
For the rest of the paper, fix a field $k$ of characteristic $p$, where $p$ is an \emph{odd} prime; the case $p=2$ requires a separate treatment, which we omit in this paper. We now recall some background in algebra.
Let $A$ be an $\mathcal{A}_{\infty}$-category and $N$ be an $A-A$-bimodule, cf. \cite[Definition 2.12]{Gan1} for a definition.
\begin{mydef}
The \emph{Hochschild chain complex} (or \emph{cyclic bar complex}) of $A$ with coefficients
in $N$ is
\begin{equation}
CC_*(A,N):=\bigoplus_{X_0,X_1,\cdots,X_k} N(X_k,X_0)\otimes A(X_0,X_1)\otimes\cdots\otimes A(X_{k-1},X_k),
\end{equation}
with grading given by $\deg(y\otimes x_1\otimes\cdots\otimes x_k)=|y|+\sum_{i=1}^k\|x_i\|$ and differential given by
\begin{equation*}
b(y\otimes x_1\otimes\cdots\otimes x_k)=\sum (-1)^{\sharp_j^i}\mu_{N}^{j|1|i}(x_{k-j+1},\cdots,x_k,y,x_1,\cdots,x_i)\otimes x_{i+1}\otimes\cdots\otimes x_{k-j}
\end{equation*}
\begin{equation}
+\sum (-1)^{\maltese_{-k}^{-(s+j+1)}}\mathbf{m}\otimes x_1\otimes \cdots\otimes x_s\otimes \mu_{A}^j(x_{s+1},\cdots,x_{s+j})\otimes x_{s+j+1}\otimes\cdots\otimes x_k,
\end{equation}
where the signs $\sharp$ and $\maltese$ are given as in \cite[Definition 2.28]{Gan1}. When $N=A_{\Delta}$ is the diagonal bimodule (cf. \cite[Definition 2.21]{Gan1}), we denote $CC_*(A):=CC_*(A,A_{\Delta})$
\end{mydef}
The following definition will play a key role in the rest of paper.
\begin{mydef}
The \emph{$p$-fold Hochschild complex} of an $\mathcal{A}_{\infty}$-category $A$ is defined as
\begin{equation}
{}_pCC_*(A):=CC_*(A,A_{\Delta}\otimes_A\cdots\otimes_AA_{\Delta}),
\end{equation}
where on the right hand side we take the $p$-fold bimodule tensor product (cf. \cite[Definition 2.19]{Gan1}). We denote its differential as $b^p$.
\end{mydef}
Explicitly, the underlying graded vector space of ${}_pCC_*(A)$ is given by
\begin{equation}
\bigoplus_{X^i_1,\cdots,X^i_{k_i}, 1\leq i\leq p} A(X^p_{k_p},X^1_1)\otimes A(X^1_1,X^1_2)[1]\otimes\cdots\otimes A(X^1_{k_1-1},X^1_{k_1})[1]\otimes \cdots\otimes A(X^{p-1}_{k_{p-1}},X^p_1)\otimes A(X^{p}_1,X^p_2)[1]\otimes A(X^{p}_{k_{p}-1},X^p_{k_p})[1]
\end{equation}
In (1.5), the morphism spaces without degree shift should be considered as coming from the $p$ diagonal bimodules entries. In particular, ${}_pCC_*(A)$ is spanned by elements of the form
\begin{equation}
\mathbf{x}^1\otimes x_1^1\otimes \cdots\otimes x^1_{k_1}\otimes\mathbf{x}^2\otimes x_1^2\otimes \cdots\otimes x^2_{k_2}\otimes \cdots\otimes \mathbf{x}^p\otimes x_1^p\otimes \cdots\otimes x^p_{k_p} \in{}_pCC_*(A),
\end{equation}
where $\mathbf{x}^1,\cdots,\mathbf{x}^p$ are the bimodule entries.
When $A$ is cohomologically unital, ${}_pCC_*(A)$ is quasi-isomorphic to $CC_*(A)$ as there is a quasi-isomorphism of $A-A$-bimodules $A_{\Delta}\otimes_AA_{\Delta}\simeq A_{\Delta}$.
Moreover, there is a chain level $\mathbb{Z}/p$-action on ${}_pCC_*(A)$ given by cyclically permuting the $p$-fold bimodule tensor product. Explicitly, on an element of the form (1.6)
the standard generator $\tau\in\mathbb{Z}/p$ acts by
$$\tau: \mathbf{x}^1\otimes x_1^1\otimes \cdots\otimes x^1_{k_1}\otimes\mathbf{x}^2\otimes x_1^2\otimes \cdots\otimes x^2_{k_2}\otimes \cdots\otimes \mathbf{x}^p\otimes x_1^p\otimes \cdots\otimes x^p_{k_p}\mapsto $$
\begin{equation}
(-1)^{\dag}\mathbf{x}^p\otimes x_1^p\otimes \cdots\otimes x^p_{k_p}\otimes \mathbf{x}^1\otimes x_1^1\otimes \cdots\otimes x^1_{k_1}\otimes\cdots\otimes \mathbf{x}^{p-1}\otimes x_1^{p-1}\otimes \cdots\otimes x^{p-1}_{k_{p-1}},
\end{equation}
where
\begin{equation}
\dag=\big(|\mathbf{x}^p|+\sum_{i=1}^{k_p} \|x^p_i\|\big)\cdot\big(\sum_{j=1}^{p-1} |\mathbf{x}^j|+\sum_{j=1}^{p-1}\sum_{i=1}^{k_j}\|x^j_{i}\|\big)
\end{equation}
is the Koszul sign.
\begin{mydef}
The
\emph{negative $\mathbb{Z}/p$-equivariant Hochschild complex of $A$} is
\begin{equation}
CC^{\mathbb{Z}/p}_*(A):={}_pCC_*(A)[[t,\theta]], |t|=2,\;|\theta|=1,\;\theta^2=0,
\end{equation}
equipped with the $t$-linear differential
\begin{equation}
\begin{cases}
x\mapsto b^px+(-1)^{|x|}(\tau x-x),\\
x\theta\mapsto b^px\,\theta+(-1)^{|x|}(x+\tau x+\cdots+\tau^{p-1}x)t
\end{cases}.
\end{equation}
\end{mydef}
There is an action of $k[[t,\theta]]$ on $CC^{\mathbb{Z}/p}_*(A)$ defined on the chain level where $t$ acts by
\begin{equation}
\begin{cases}
xt^k\mapsto (-1)^{|x|}xt^{k+1},\\
xt^k\theta\mapsto (-1)^{|x|}xt^{k+1}\theta
\end{cases}
\end{equation}
and $\theta$ acts by
\begin{equation}
\begin{cases}
xt^k\mapsto (-1)^{|x|}xt^{k}\theta,\\
xt^k\theta\mapsto (-1)^{|x|}(\tau-1)^{p-2}xt^{k+1}\theta
\end{cases}.
\end{equation}
This action descends to cohomology and makes $HH^{\mathbb{Z}/p}_*(A)$ a $k[[t,\theta]]$-module. \par\indent
Analogous to the cyclic open-closed map, in section 2 we define the \emph{$\mathbb{Z}/p$-equivariant open-closed map}, which is a chain map
\begin{equation}
OC^{\mathbb{Z}/p}: CC^{\mathbb{Z}/p}(\mathrm{Fuk}(M)_{\lambda})\rightarrow QH(M)[[t,\theta]], |t|=2, |\theta|=1, \theta^2=0,
\end{equation}
of degree $\frac{1}{2}\dim_{\mathbb{R}}M$. \par\indent
1.2. \textbf{The $\mathbb{Z}/p$-Gysin Comparison map for open-closed maps}.
The main results of this paper are the following.
\begin{prop}
For a cohomologically unital $\mathcal{A}_{\infty}$-category $A$, there exists a quasi-isomorphism
\begin{equation}
\Phi_p: CC^{S^1}(A)\oplus CC^{S^1}(A)\theta\simeq {}_pCC^{\mathbb{Z}/p}(A),
\end{equation}
where
$\theta$ is a formal variable of degree $1$.
\end{prop}
\begin{thm}
The following diagram is homotopy commutative:
\begin{equation}
\begin{tikzcd}[row sep=1.2cm, column sep=0.8cm]
{}_pCC^{\mathbb{Z}/p}(\mathrm{Fuk}(M)_{\lambda})\arrow[rrr,"{OC^{\mathbb{Z}/p}}"]& & &QH(M)[[t,\theta]]  \\
CC^{S^1}(\mathrm{Fuk}(M)_{\lambda})\oplus CC^{S^1}(\mathrm{Fuk}(M)_{\lambda})\theta \arrow[u,"{\Phi_p}"]\arrow[urrr,"{OC^{S^1}\oplus OC^{S^1}\theta}"]& & &
\end{tikzcd}
\end{equation}
\end{thm}
\begin{rmk}
The gist of Proposition 1.5 and Theorem 1.6 is analogous to the following classical fact in topology: the homology of the $\mathbb{Z}/p\subset S^1$ homotopy fixed points of an $S^1$-space $X$ is isomorphic to two copies of the homology of its $S^1$ homotopy fixed points, one with shifted down degree by $1$. The proof is a simple application of the Gysin long exact sequence associated to the $S^1$-bundle $X\times_{\mathbb{Z}/p} ES^1\rightarrow X\times_{S^1}ES^1$. Proposition 1.5 can be seen as a chain complex version of that statement, applied to the `$C_*(S^1)$-complex' $CC(A)$. We hence call it the $\mathbb{Z}/p$-Gysin comparison. In spirit, Theorem 1.6 comes from the naturality of the $\mathbb{Z}/p$-Gysin comparison under an $S^1$-map. \par\indent
In particular, to prove Theorem 1.6, we would like to apply the $\mathbb{Z}/p$-Gysin comparison map to the map $OC: HH_*(\mathrm{Fuk}(M)_{\lambda})\rightarrow QH^*(M)$. Some of the difficulties in realizing this vision to an actual proof are:
\begin{itemize}
    \item In the symplectic literature, there is a preferred chain model $CC(A)$ for the Hochschild homology of an $\mathcal{A}_{\infty}$-algebra $A$, known as the \emph{cyclic bar complex}.
On the chain level, $C_*(S^1)$ does not act on this complex. Instead, $CC(A)$ admits an action of the graded algebra $k[\epsilon]:=k[\epsilon]/\epsilon^2, |\epsilon|=-1$, where $\epsilon$ acts as
the Connes' differential, cf. section 2. There is an $\mathcal{A}_{\infty}$-quasi-equivalence $k[\epsilon]\simeq C_*(S^1)$, where $C_*(S^1)$ is equipped with the Pontryagin product coming from the multiplicative structure of $S^1$, and thus we think of the action as coming from $S^1$.
Moreover, the open-closed map $OC$ is a map of $\mathcal{A}_{\infty}$ $k[\epsilon]/\epsilon^2$-modules, where $k[\epsilon]$ acts trivially on $QH(M)$. For us, the issue of using $k[\epsilon]$ as a small model for $C_*(S^1)$ is that there is no chain level inclusion $k[\mathbb{Z}/p]\subset k[\epsilon]/\epsilon^2$
in order to extract the `induced $\mathbb{Z}/p$-action' from the `$S^1$-action' on the cyclic bar complex $CC(A)$.
    \item Existing solutions in the literature to the above issue would be to use a larger chain model for $HH(A)$ that is equipped with a chain level $\mathbb{Z}/p$-action. For instance, one can consider the chain complex $CC(A)\otimes^{\mathbb{L}}_{k[\epsilon]} C_*(S^1)$ or $CC(A)\otimes^{\mathbb{L}}_{k[\epsilon]} k[\tau,\sigma]$; see (5.20) for the definition of $k[\tau,\sigma]$, which is a small dg model for $C_*(S^1)$ that sees the $p$-th root of unities. On both complexes, $\mathbb{Z}/p$ acts diagonally: trivial on the first tensor factor, and acts via the natural inclusions $k[\mathbb{Z}/p]\subset C_*(S^1), k[\mathbb{Z}/p]\subset k[\tau,\sigma]$ on the second tensor factor. Analogues of these complexes in the context of symplectic cohomology were used in \cite{Sen} to prove a compatibility result of the $\mathbb{Z}/p$ and $S^1$ action on symplectic cohomology, cf. \cite[Proposition 3.4.2]{Sen}.
However, our application to Quantum Steenrod operations requires a particular chain model for computing Hochschild homology, namely the $p$-fold cyclic bar complex ${}_pCC(A)$.
Thus, one would for instance need to prove an equivalence of chain complexes $CC(A)\otimes^{\mathbb{L}}_{k[\epsilon]}C_*(S^1)\simeq
{}_pCC(A)$ or $CC(A)\otimes^{\mathbb{L}}_{k[\epsilon]} k[\tau,\sigma]\simeq {}_pCC(A)$ that intertwines the $\mathbb{Z}/p$-actions (and moreover show that this quasi-isomorphism intertwines the open-closed maps), which as far as the author is concerned, is complicated to construct using explicit homological algebra (and prove using explicit TQFT arguments).
\end{itemize}
In view of these issues, we turn to the formulation of Hochschild homology as a cyclic object. This approach solves the above issues by being both
\begin{itemize}
    \item computable: using the standard cyclic bar complex (resp. $p$-fold cyclic bar complex) associated to a cyclic module (resp. finite $p$-cyclic module), see section 3, which generalizes
the chain complexes $CC(A)$ and ${}_pCC(A)$ for an $\mathcal{A}_{\infty}$-category $A$.
    \item sufficiently functorial: given a cyclic module, which models a chain complex with a circle action, there is a way to extract the underlying finite cyclic subgroup action
using ideas of edgewise subdivision.
\end{itemize}
\end{rmk}
We now briefly summarize the logic of the proof. \par\indent
First, we restate Proposition 1.5 and Theorem 1.6 in a more homotopy coherent framework by using a simplicial model for Hochschild homology. This perspective was initiated by Connes' introduction of the cyclic category $\Lambda$ in \cite{Con}, and his
observation that Hochschild homology of an associative algebra can be modeled as a cyclic module (i.e. functors out of $\Lambda$), and has since been widely used in e.g. \cite{Lod},\cite{Hoy} and \cite{NS}.
Building on this perspective, the main novelties in this paper can be summarized as follows:
\begin{enumerate}[label=\arabic*)]
    \item The construction of an $\mathcal{A}_{\infty}$-version of the cyclic category, $\Lambda\rtimes \mathcal{A}_{\infty}^{dg}$, as well as its finite cyclic versions. Generalizing the classical construction of Connes, we constructed
an $\mathcal{A}_{\infty}$-cyclic module $A^{\sharp}$ modelling the Hochshild chain complex of an $\mathcal{A}_{\infty}$-category $A$.
    \item The proof of a $\mathbb{Z}/p$-Gysin comparison theorem for cyclic modules. This gives the quasi-isomorphism $\Phi_p$ of Proposition 1.5. \item The construction of an \emph{operadic open-closed maps} for a closed monotone symplectic manifold. This can be viewed as a simplicial lift of the classical open-closed map
$OC: HH_*(\mathrm{Fuk}(M)_{\lambda})\rightarrow QH^*(M)$. A simple but key observation in this construction is that chains on the moduli space of disks with $1$ interior marked point and several boundary marked points, $\overline{\mathcal{R}}^1_{\bullet}$, have a
natural ($\mathcal{A}_{\infty}$-) cocyclic structure coming from its geometry. From the operadic open-closed map one can extract the operadic cyclic and $\mathbb{Z}/p$-equivariant open-closed maps by taking
suitable homotopy colimits. Using this framework, we prove a version of the compatibility result of Theorem 1.6 involving operadic cyclic and $\mathbb{Z}/p$-equivariant open-closed maps.
     \item Finally we show that the operadic cyclic and $\mathbb{Z}/p$-equivariant open-closed map agree with their classical counterparts constructed in \cite{Gan2} and \cite{Che}, respectively. This includes the preliminary algebraic comparison result between the operadic cyclic chain complex (resp. operadic finite cyclic chain complex), which are constructed
as certain abstract homotopy (co)limits, with the explicit chain complexes $CC^{S^1}(A)$ (resp. $CC^{\mathbb{Z}/p}(A)$). This generalizes a theorem of \cite{Hoy}.
\end{enumerate}
As a consequence of 1) to 4), we prove Proposition 1.5 and Theorem 1.6. 1), 2) and part of 4), which are results purely in algebra, may have independent interest.\par\indent
The organization of this paper is as follows. In section 2 we recall the definition of the monotone Fukaya category and the cyclic and $\mathbb{Z}/p$-equivariant open-closed maps. In section 3 we define an $\mathcal{A}_{\infty}$-version of Connes cyclic category, $\Lambda\rtimes \mathcal{A}_{\infty}^{dg}$, and its variants (notably the finite cyclic version ${}_p\Lambda\rtimes \mathcal{A}_{\infty}^{dg}$).
We prove that certain homotopy (co)limits of $\mathcal{A}_{\infty}$-(finite)cyclic modules can be computed by explicit chain complexes involving the ($p$-fold) cyclic bar construction, generalizing \cite[Theorem 2.3]{Hoy}. In section 4
we start by reviewing the setup of \cite{AGV}'s operadic Floer theory. We then use this to define the operadic open-closed maps. In section 5, we prove a $\mathbb{Z}/p$-Gysin comparison theorem
for cyclic modules. Combining this with the result of section 4, we show that the compatibility result of Theorem 1.6 holds for the operadic cyclic and $\mathbb{Z}/p$-equivariant open-closed map.  In section 6, we prove that the operadic cyclic open-closed maps (resp. operadic $\mathbb{Z}/p$-equivariant open-closed maps)
agree with the usual cyclic open-closed maps defined by \cite{Gan2} (resp. $\mathbb{Z}/p$-equivariant open-closed maps), and conclude the proof of Theorem 1.6.
\section*{Acknowledgements}
First and foremost, I would like to thank my advisor Paul Seidel for his continuous support and guidance throughout my graduate study. I would also like to thank Mohammed Abouzaid,
Nathaniel Bottman, Marc Hoyois and Hiro Tanaka for helpful discussions at various points. This work was partially supported by the Simons Foundation through a Simons Investigator grant (256290).

\renewcommand{\theequation}{2.\arabic{equation}}
\setcounter{equation}{0}

\section{Fukaya categories and open-closed maps}
2.1. \textbf{The monotone Fukaya category}. In this subsection, we review the definition of the monotone Fukaya category associated to a closed monotone symplectic manifold $(M,\omega)$. For a detailed treatment on the subject, see \cite[section 2.3]{She}. In fact, to each $\lambda\in k$, one can associate a $k$-linear $\mathbb{Z}/2$-graded $\mathcal{A}_{\infty}$-category $\mathrm{Fuk}(M)_{\lambda}$. \par\indent
Let $L$ be an oriented spin monotone Lagrangian submanifold $L\subset M$ equipped with a $k^*$-local system. Recall that monotonicity means that $\mu(L)=[\omega]$ considered as classes in $H^2(M,L)$, where $\mu$ denotes the Maslov class. Orientability implies that the minimal Maslov number is $\geq 2$. By an abuse of notation,
we denote this datum simply by its underlying Lagrangian $L$. Let $\mathcal{J}$ denote the space of compatible almost complex structures and $\mathcal{H}:=C^{\infty}(M,\mathbb{R})$ the space of Hamiltonians. For each $L$, we fix $J_L\in \mathcal{J}$. For each pair $(L_0,L_1)$, we fix $J_t\in C^{\infty}([0,1],\mathcal{J})$ and $H_t\in C^{\infty}([0,1],\mathcal{H})$ such that $J(i)=J_{L_i}$.
If the $k^*$ local systems on both Lagrangians are trivial, the morphism space $CF^*(L_0,L_1)$ is the $k$-vector space generated by time-$1$ Hamiltonian chords of $H_t$ from $L_0$ to $L_1$; in general, it is the direct sum of hom spaces between the fibers of the local systems at the startpoint and endpoint of the chord. \par\indent
Fix Lagrangians $L_0,L_1$. By standard transversality arguments, for generic almost complex structure $J_{L_0}, J_{L_1}$ and one parameter family $(H_t,J_t), t\in [0,1]$ such that $J_0=J_{L_0}, J_1=J_{L_1}$:
\begin{enumerate}[label=R\arabic*)]
    \item The moduli space $\mathcal{M}_1(L_0)$ of Maslov index $2$ $J$-holomorphic disks with one boundary marked point and boundary on $L_0$ is regular.
    \item The moduli space $\mathcal{M}_1(J_t)$ of pairs $(t,u)$, where $t\in [0,1]$ and $u$ is a Chern number $1$ $J_t$-holomorphic sphere with one marked point, is regular.
    \item For any time $1$ Hamiltonian chord $\gamma: [0,1]\rightarrow X$ starting on $L_0$ and ending on $L_1$, the map
\begin{equation}
(\gamma\circ t, \mathrm{ev}): \mathcal{M}_1(J_t)\rightarrow X\times X
\end{equation}
avoids the diagonal. In other words, all $J_t$ holomorphic spheres avoid $\gamma$.
    \item For Hamiltonian chords $x,y$ from $L_0$ to $L_1$, the moduli space $\mathcal{M}(x,y,J_t,H_t)$ of strips satisfying Floer's equation
\begin{equation}
(du-X_{H_t}\otimes dt)^{0,1}_{J_t}=0
\end{equation}
is regular.
\end{enumerate}
For each $L$,  since $L$ has minimal Maslov number $\geq 2$, the only possible nodal configuration in the Gromov compactification of $\mathcal{M}_1(L)$ is a $J_L$ holomorphic sphere of Chern number $1$ attached to a constant disk on $L$, and the moduli space of those has codimension $2$. In particular, $\mathcal{M}_1(L)$ has a well-defined pseudocycle fundamental class.
If $L$ is equipped with the trivial local system, we define $w(L)\in \mathbb{Z}$ by
\begin{equation}
\mathrm{ev}_*[\mathcal{M}_1(L)]=w(L)[L]\in H_n(L,\mathbb{Z}).
\end{equation}
More generally, $\mathrm{ev}_*$ is weighted by the monodromy of the local system around the boundary of the disc, and $w(L)$ defines an element of $k$. \par\indent
The Floer differential $\mu^1:CF^*(L_0,L_1)\rightarrow CF^*(L_0,L_1)$ is defined as follows. Let $x_-,x_+\in CF^*(L_0,L_1)$,
then the coefficient of $x_-$ in $\mu^1(x_+)$ is the signed count (again, weighted by monodromy) of isolated elements of $\mathcal{M}(x_-,x_+)/\mathbb{R}$ (when $\mathcal{M}(x_-,x_+)$ has dimension $1$),
where $\mathcal{M}(x_-,x_+)$ is the moduli space of $u:\mathbb{R}\times [0,1]\rightarrow M$ such that
\begin{equation}
\begin{cases}
\partial_su+J_t(\partial_tu-X_{H_t})=0\\
u(s,0)\in L_0,\;\; u(s,1)\in L_1\\
\lim_{s\rightarrow \pm\infty}u(s,\cdot)=x_{\pm}.
\end{cases}
\end{equation}
By Gromov compactness and monotonicity, when $\mathcal{M}(x_-,x_+)$ is one dimensional (i.e. the Maslov index of $u$ is $1$), the space $\mathcal{M}(x_-,x_+)/\mathbb{R}$ is compact. When $\mathcal{M}(x_-,x_+)$ is $2$-dimensional, its Gromov compactification consists of broken strips $u_1,u_2$, each with Maslov index $1$, as well as a Maslov index $2$ disk bubbling of a Maslov index $0$ (hence constant in $s$) strip. For generic $J_t$, sphere bubbling cannot occur by regularity assumption R3). Therefore, we have
\begin{equation}
\mu^1(\mu^1(x))=(w(L_0)-w(L_1))x.
\end{equation}
Hence, if $w(L_0)=w(L_1)$, then $(\mu^1)^2=0$.  \par\indent
We define the objects of $\mathrm{Fuk}(M)_{\lambda}$ to be oriented spin monotone Lagrangian submanifolds $L\subset M$ equipped with a $k^*$-local system, such that $w(L)=\lambda$. The morphism chain complexes
are defined as $(CF^*(L_0,L_1),\mu^1)$. \par\indent
We now describe the higher $\mathcal{A}_{\infty}$ operations in $\mathrm{Fuk}(M)_{\lambda}$. Let $S$ be a surface with boundary and interior marked points (the boundary marked points are thought of as punctures).
Given a Lagrangian labeling $\mathbf{L}$ of the boundary components of $S$, a \emph{labeled Floer datum} for $S$ consists of the following data:
\begin{itemize}
    \item for each boundary marked point $\zeta$, a strip-like end $\epsilon_{\zeta}: \mathbb{R}^{\pm}\times [0,1]\rightarrow S$ at $\zeta$ (the strip-like end being positive or negative depending on whether $\zeta$ is an input or an output);
    \item a choice of $K\in \Omega^1(S,\mathcal{H})$ and $J\in C^{\infty}(S,\mathcal{J})$ such that $K(\xi)|_{L_C}=0$ for all $\xi\in TC$, where $C$ is a boundary component and $L_C$ is the corresponding Lagrangian label. Moreover, $K,J$ are compatible with strip-like ends in the sense that
\begin{equation}
\epsilon_{\zeta}^*K=H_{\zeta}(t)dt,\quad J(\epsilon_{\zeta}(s,t))=J_{\zeta}(t),
\end{equation}
where $H_{\zeta}, J_{\zeta}$ are the chosen Hamiltonian and almost complex structure for the pair of Lagrangians meeting at $\zeta$. We also require $J=J_L$ when restricted to a boundary component labeled $L$. The pair $(K,J)$ is called a \emph{perturbation datum}.
\end{itemize}
The higher $\mathcal{A}_{\infty}$-operations of $\mathrm{Fuk}(M)_{\lambda}$ are governed by the Delign-Mumford moduli space of disks with boundary marked points. Let $\mathcal{R}^{d+1}$ be the moduli space of disks with one boundary output and $d$ boundary inputs. It admits a compactification to a manifold with corners $\overline{\mathcal{R}}^{d+1}$ given by
\begin{equation}
\overline{\mathcal{R}}^{d+1}=\coprod_T \mathcal{R}^T,
\end{equation}
where $T$ ranges over all planar stable $d$-leafed trees and $\mathcal{R}^T:=\prod_{v\in \mathrm{Ve}(T)}\mathcal{R}^{|v|}$.\par\indent
We make a \emph{consistent choice of labeled Floer data} for $\mathcal{R}^{d+1}, d\geq 2$, meaning it is compatible with the product of Floer data of lower dimensional $\mathcal{R}^{d'}$'s near a boundary stratum, see \cite[section (9g),(9i)]{Sei1}.\par\indent
The higher operations $\mu^d, d\geq 2$ are then defined by counting isolated elements of the parametrized moduli space $\mathcal{M}(y_1,\cdots,y_d;y_-)$, which is the space of $(r,u)$, $r\in \mathcal{R}^{d+1}, u:\mathcal{S}_r\rightarrow M$ satisfying
\begin{equation}
(du-Y_K)^{0,1}_J=0,
\end{equation}
with appropriate Lagrangian boundary and asymptotic conditions, where $Y_K$ is the one-form on $S$ with value the Hamiltonian vector field associated to $K$. For a generic choice of Floer data, this moduli space is regular (\cite[section (9k)]{Sei1}).\par\indent 
2.2. \textbf{Ganatra's cyclic open-closed map}. Let $A=\mathrm{Fuk}(M)_{\lambda}$. Pick a Morse-Smale function $f$ on $M$ and consider its Morse complex $CM^*(f)$, equipped with the small quantum cup product defined using $3$-pointed Gromov Witten
invariants (with incidence constraints on the stable and unstable manifolds of critical points of $f$), see \cite[section 3]{SW}. We use this as a chain model for quantum cohomology.
There is a chain map, cf. for instance \cite[section 5]{Gan1}, called the \emph{open-closed map},
\begin{equation}
OC: CC_*(\mathrm{Fuk}(M)_{\lambda})\rightarrow CM^*(f)
\end{equation}
of degree $\frac{1}{2}\dim_{\mathbb{R}}M$. It is defined by counting rigid solutions to the Floer equation (2.8), as the domain varies over \begin{equation}
\mathcal{R}^1_{d+1},
\end{equation}
the parameter space of
disks with $1$ interior marked point and $d+1$ boundary marked points. The boundary marked points of the disk are asymptotic to an element of $CC_*(\mathrm{Fuk}(M)_{\lambda})$,
and the interior marked points is incidence to the stable manifold of a critical point of $f$.\par\indent
We also use the \emph{non-unital Hochschild complex}, cf. \cite[section 3.1]{Gan2}. As a graded vector space, it is defined to be
\begin{equation}
CC^{nu}_*(A):=CC_*(A)\oplus CC_*(A)[1].
\end{equation}
To define the differential, one considers the \emph{bar differential}
\begin{equation}
b'(x_0\otimes\cdots\otimes x_d):=\sum_{-1\leq s\leq k-1,1\leq j\leq k-s} (-1)^{\maltese_{-d}^{-(s+j+1)}} x_0\otimes\cdots\otimes x_s\otimes \mu(x_{s+1},\cdots,x_{s+j})\otimes x_{s+j+1}\otimes\cdots\otimes x_d
\end{equation}
and
\begin{equation}
d_{\wedge\vee}(x_0\otimes\cdots\otimes x_d):=(-1)^{\maltese_0^{d-1}+\|x_d\|\cdot\maltese_0^{d-1}+1} x_d\otimes x_0\otimes\cdots\otimes x_{d-1}+(-1)^{\maltese^{d}_1}x_0\otimes\cdots\otimes x_d.
\end{equation}
Then, the differential $b^{nu}$ is defined as
\[b^{nu}:=
\begin{pmatrix}
b&d_{\wedge\vee}\\
0& b'
\end{pmatrix},
\]
where $b$ is the Hochschild differential (1.3). When $A$ is homologically unital, the natural inclusion $CC_*(A)\hookrightarrow CC^{nu}_*(A)$ is a quasi-isomorphism, cf. \cite[Proposition 2.2]{Gan1}. \par\indent
There is a chain level unital $k[\epsilon]/\epsilon^2$-action on $CC^{nu}_*(A)$, where $\epsilon$ (of degree $-1$) acts as the \emph{non-unital Connes operator}, given by
\begin{equation}
B^{nu}(x_0\otimes \cdots\otimes x_k,y_0\otimes\cdots\otimes y_l):=\sum_i (-1)^{\maltese_{i+1}^k\maltese_0^i+\|x_0\|+\maltese_0^k+1}(0,x_{i+1}\otimes \cdots\otimes x_k\otimes x_0\otimes \cdots\otimes \cdots\otimes x_i).
\end{equation}
The \emph{(non-unital) negative cyclic chain complex} of $A$ is defined as
\begin{equation}
CC^{S^1}_*(A):=CC^{nu}_*(A)[[t]],\;\;|t|=2,
\end{equation}
equipped with the differential $b^{nu}+B^{nu}t$. \par\indent
The open-closed map can be enhanced to an $S^1$-equivariant version $OC^{S^1}$, called the \emph{negative cyclic open-closed map}. This is a chain map
\begin{equation}
OC^{S^1}:CC^{S^1}_*(\mathrm{Fuk}(M)_{\lambda})\rightarrow CM^*(f)[[t]]
\end{equation}
of degree $\frac{1}{2}\dim_{\mathbb{R}}M$ of the form
\begin{equation}
OC^{S^1}=\sum_{k\geq 0} (\check{OC}^k\oplus \hat{OC}^k)t^k,
\end{equation}
where
\begin{equation}
\check{OC}^k:CC_*(\mathrm{Fuk}(M)_{\lambda})\rightarrow CM^{*+\frac{1}{2}\dim_{\mathbb{R}}M-2k}(f)
\end{equation}
and
\begin{equation}
\hat{OC}^k:CC_*(\mathrm{Fuk}(M)_{\lambda})[1]\rightarrow CM^{*+\frac{1}{2}\dim_{\mathbb{R}}M-2k}(f)
\end{equation}
are operations
defined using parametrized moduli problems associated to the parameter spaces of domains ${}_k\check{\mathcal{R}}^1_{d+1}$ and ${}_k\hat{\mathcal{R}}^1_d$, respectively, cf. \cite[section 5.5]{Gan2}. We briefly recall the description of these spaces. \par\indent
First, we define the space
\begin{equation}
\mathcal{R}^{d+1,f_i}
\end{equation}
of disks with forgotten marked point to be a copy of $\mathcal{R}^{d+1}$, except that the $i$-th marked point is labelled as `forgotten'. Topologically, the compactification $\overline{\mathcal{R}}^{d+1,f_i}$ is also just $\overline{\mathcal{R}}^{d+1}$, but one keeps track of the forgotten marked point.
A universal choice of Floer data for $\overline{\mathcal{R}}^{d+1,f_i}, d\geq 2$ is a universal choice of Floer data for the underlying $\overline{\mathcal{R}}^{d+1}$, that further satisfy the following conditions:
\begin{itemize}
    \item the Floer data on the unique element of $\overline{\mathcal{R}}^{3,f_i}$ is translation invariant after forgetting the $i$-th point,
    \item for $d>2$, the Floer data on $\overline{\mathcal{R}}^{d+1,f_i}$ is pulled back along the forgetful map $\overline{\mathcal{R}}^{d+1,f_i}\rightarrow \overline{\mathcal{R}}^{d}$.
\end{itemize}
Fixing a universal choice of Floer data for $\overline{\mathcal{R}}^{d+1}, d\geq 2$, the two conditions above uniquely determines a universal choice of Floer data for $\overline{\mathcal{R}}^{d+1,f_i}$. We now consider four types of parameter spaces of disks involved in the definition of the cyclic open-closed map.\par\indent
The first space ${}_k\check{\mathcal{R}}^1_{d+1}$ is the moduli space of discs with $d+1$ positive boundary marked points $z_0,\cdots, z_d$ labeled in counterclockwise order, 1 interior negative puncture $z_{out}$, and $k$ additional interior marked points $p_1,\cdots, p_k$. Choosing a representative of an element of this moduli space which fixes $z_0$ at 1 and $z_{out}$ at $0$ on the unit disc, the $p_i$ should be strictly radially ordered; that is,
\begin{equation}
0<|p_1|<\cdots<|p_k|<\frac{1}{2}.
\end{equation}
The second space ${}_{k-1}\check{\mathcal{R}}_{d+1}^{S^1}$ topologically is just ${}_{k-1}\check{\mathcal{R}}^1_{d+1}\times S^1$, but we view it as the sublocus of ${}_k\check{\mathcal{R}}^1_{d+1}$ where the interior point $p_k$ is constrained by $|p_k|=\frac{1}{2}$.\par\indent
The codimension 1 boundary components of the compactification ${}_k\overline{\check{\mathcal{R}}^1}_{d+1}$ are of the following three types
\begin{equation}
\mathcal{R}^{s+1}\times {}_k\check{\mathcal{R}}^1_{d-s+1}
\end{equation}
\begin{equation}
{}_{k-1}\check{\mathcal{R}}_{d+1}^{S^1}
\end{equation}
\begin{equation}
{}_k^{i,i+1}\check{\mathcal{R}}^1_{d+1}
\end{equation}
(2.22) corresponds to a disk of type $\mathcal{R}^s$ bubbling off the boundary of the main component. (2.23) is the locus where the norm of $p_k$ goes to $\frac{1}{2}$. Finally, in (2.24), ${}_k^{i,i+1}check{\mathcal{R}}^1_{d+1}$ denotes the sublocus of ${}_k\check{\mathcal{R}}^1_{d+1}$ such that $|p_i|=|p_{i+1}|$. \par\indent
The compactification ${}_{k-1}\overline{\check{\mathcal{R}}}_{d+1}^{S^1}$ is by definition ${}_{k-1}\overline{\check{\mathcal{R}}}_{d+1}^{1}\times S^1$.\par\indent
The third space ${}_k\hat{\mathcal{R}}^1_d$ is the moduli space of discs with $d+1$ positive boundary marked points $z_f,z_0,\cdots, z_{d-1}$ labeled in counterclockwise order, 1 interior negative puncture $z_{out}$, and $k$ additional interior marked points $p_1,\cdots, p_k$. Choosing a representative of an element of this moduli space which fixes $z_f$ at 1 and $z_{out}$ at $0$ on the unit disc, the $p_i$ should be strictly radially ordered:
\begin{equation}
0<|p_1|<\cdots<|p_k|<\frac{1}{2}.
\end{equation}
Topologically, this is the same as ${}_k\check{\mathcal{R}}^1_{d+1}$, but we would like to view $z_f$ as a special boundary marked point. \par\indent
The fourth space ${}_{k-1}\hat{\mathcal{R}}^{S^1}_d$ is the sublocus of ${}_k\hat{\mathcal{R}}^1_d$ where $|p_k|=\frac{1}{2}$, and is topologically just ${}_{k-1}\hat{\mathcal{R}}^1_d\times S^1$.
\begin{figure}[H]
 \centering
 \includegraphics[width=1\textwidth]{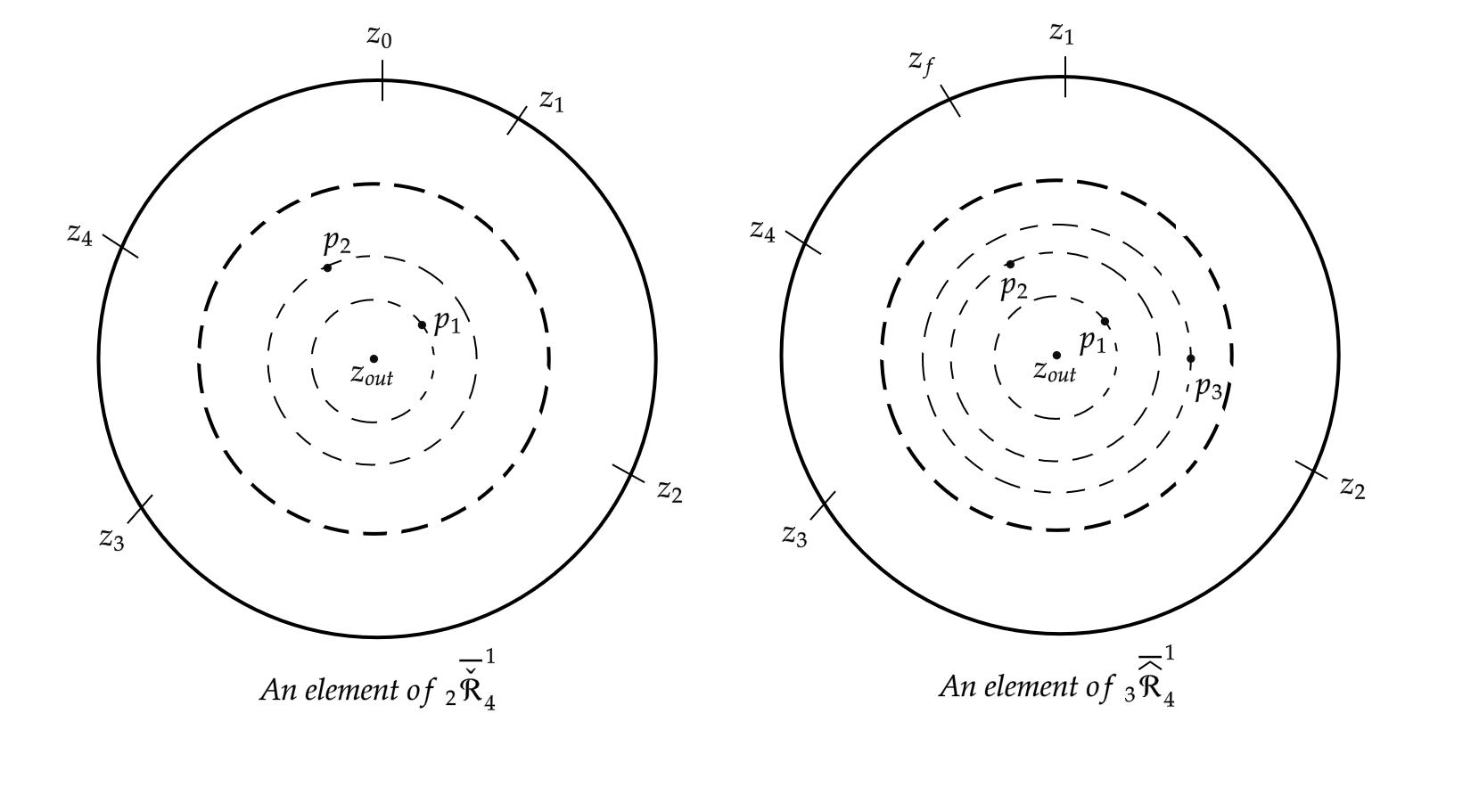}
 \caption{}
\end{figure}
The codimension 1 boundary components of
${}_k\hat{\mathcal{R}}^1_{d}$ are of the following four types
\begin{equation}
\mathcal{R}^{s+1}\times {}_k\hat{\mathcal{R}}^1_{d-s+1}
\end{equation}
\begin{equation}
\mathcal{R}^{m+1,f_i}\times_{d-m} {}_k\check{\mathcal{R}}^1_{d-m+1}
\end{equation}
\begin{equation}
{}_{k-1}\hat{\mathcal{R}}_{d}^{S^1}
\end{equation}
\begin{equation}
{}_k^{i,i+1}\hat{\mathcal{R}}^1_{d}
\end{equation}
For a nodal configuration where a disk bubbles off the boundary, if $z_f$ lies on the main component, this configuration gives rise to an element of (2.26); if $z_f$ lies on the disk bubble, it gives rise to an element of (2.27). (2.28) corresponds to the sublocus of ${}_k\hat{\mathcal{R}}^1_d$ where $|p_k|=\frac{1}{2}$. In (2.29), ${}_k^{i,i+1}\hat{\mathcal{R}}^1_{d}$ denotes the sublocus of
${}_k^{i,i+1}\hat{\mathcal{R}}^1_{d}$ where $|p_i|=|p_{i+1}|$.\par\indent
The compactification ${}_{k-1}\overline{\hat{\mathcal{R}}}_{d}^{S^1}$ is by definition ${}_{k-1}\overline{\hat{\mathcal{R}}}_{d}^{1}\times S^1$.\par\indent
Now, we fix a universal choice of Floer data for the spaces $\overline{\mathcal{R}}^{d+1}, \overline{\mathcal{R}}^{d+1,f_i}, d\geq 2, 0\leq i\leq d$. A \emph{Floer data for the cyclic open-closed map} is an inductive choice of Floer data on each surface in the parameter spaces ${}_k\check{\mathcal{R}}^1_{d+1}, {}_k\check{\mathcal{R}}^{S^1}_{d+1} {}_k\hat{\mathcal{R}}^1_{d}, {}_k\hat{\mathcal{R}}^{S^1}_{d}$, such that the Floer datum on a nodal surface agrees with the product of prior chosen Floer datum of the nodal components, and that they satisfy the following conditions (cf. \cite[(5.148)-(5.155)]{Gan2}):
\begin{enumerate}[label=\arabic*)]
    \item For $S_0\in {}_k\overline{\check{\mathcal{R}}^1_{d+1}}$,
       \begin{enumerate}[label=1\alph*)]
           \item If $S_0$ belongs to a boundary stratum of type (2.24), the Floer datum on it is pulled back via the map $\pi_i: {}_k^{i,i+1}\check{\mathcal{R}}^1_{d+1}\rightarrow {}_{k-1}\check{\mathcal{R}}^1_{d+1}$ that forgets $p_{i+1}$.
        \end{enumerate}
    \item For $S_1\in {}_k\overline{\check{\mathcal{R}}^{S^1}_{d+1}}$,
     \begin{enumerate}[label=2\alph*)]
         \item If $S_1$ belongs to the codimension $1$ sublocus $ {}_k\check{\mathcal{R}}^{S^1_i}_{d+1}$ where $p_{k+1}$ points in the direction of $z_i$, the Floer datum on $S_1$ is pulled back along the map $\tau_i:{}_k\check{\mathcal{R}}^{S^1_i}_{d+1}\rightarrow {}_k\check{\mathcal{R}}_{d+1}$ which forgets $p_{k+1}$ and cyclically permutes the boundary labels so that $z_i$ becomes $z_0$.
         \item The Floer datum is equivariant under cyclically permuting the boundary labels.
         \item If $S_1$ belongs to a boundary stratum of the form ${}_k^{i,i+1}\check{\mathcal{R}}^{S^1}_{d+1}={}_k^{i,i+1}\check{\mathcal{R}}^1_{d+1}\times S^1$, the Floer datum on it is pulled back along the map $\pi^{S^1}_i:{}_k^{i,i+1}\check{\mathcal{R}}^{S^1}_{d+1}\rightarrow {}_{k-1}\check{\mathcal{R}}^{S^1}_{d+1}$ that forgets $p_{i+1}$.
     \end{enumerate}
     \item For $S_2\in {}_k\overline{\hat{\mathcal{R}}}^1_{d}$,
      \begin{enumerate}[label=2\alph*)]
          \item  Let ${}_k\check{\mathcal{R}}^{S^1_{d,0}}_{d+1}$ denote the sublocus of ${}_k\check{\mathcal{R}}^{S^1}_{d+1}$ where $p_{k+1}$ points in between $z_d$ and $z_0$. There is a diffeomorphism $\pi_f:{}_k\hat{\mathcal{R}}^1_{d}\rightarrow \pi_f:{}_k\check{\mathcal{R}}^{S^1_{d-1,0}}_{d}$ that adds a unique interior point $p_{k+1}$ of norm $\frac{1}{2}$ that points in the direction $z_f$, and then forgets $z_f$. The Floer datum on a surface $S_2$ in ${}_k\hat{\mathcal{R}}^1_{d}$ is required to be pulled back along $\pi_f$.
          \item If $S_2$ belongs to a boundary stratum of type (2.28), the Floer datum on it is pulled back along the map $\pi_{\mathrm{boundary}}: {}_{k-1}\hat{\mathcal{R}}^{S^1}_{d}\rightarrow {}_{k-1}\check{\mathcal{R}}^{S^1}_{d}$ that forgets $z_f$.
          \item If $S_2$ belongs to a boundary stratum of type (2.29), the Floer datum on it is pulled back along the map $\hat{\pi}_i: {}_k^{i,i+1}\hat{\mathcal{R}}^1_{d}\rightarrow {}_{k-1}\hat{\mathcal{R}}^1_{d}$ that forgets $p_{i+1}$.
        \end{enumerate}
\end{enumerate}
By \cite[Proposition 10]{Gan2}, a Floer data for the cyclic open-closed map exists. Given such a choice, one can define $\check{OC}^k$ and $\hat{OC}^k$ as the operations obtained from counting parametrized moduli problems associated to the spaces ${}_k\check{\mathcal{R}}^1_{d+1}$ and ${}_k\hat{\mathcal{R}}^1_{d+1}$, respectively; cf. \cite[section 5.5]{Gan2} for the detailed definition.  \par\indent
2.3. \textbf{The $p$-fold open-closed map}. We briefly review the definition of the \emph{$p$-fold open-closed map} but defer the proofs of its properties to \cite{Che}. The $p$-fold open-closed map
\begin{equation}
{}_pOC: {}_pCC_*(\mathrm{Fuk}(M)_{\lambda})\rightarrow CM^*(f),
\end{equation}
which is a chain map of degree $\frac{1}{2}\dim_{\mathbb{R}}M$; cf. Definition 1.3 for ${}_pCC_*$. Consider
\begin{equation}
\mathcal{R}^1_{k_1,\cdots,k_p},
\end{equation}
the moduli space of disks with one interior output marked point $y_{out}$ and $k_1+\cdots+k_p+p$ boundary input marked points $z^1,z^1_1,\cdots,z^1_{k_1},z^2,z^2_1,\cdots,z^2_{k_2},\cdots,z^p,z^p_1,\cdots,z^p_{k_p}$ in counterclockwise order such that up to automorphism of the disk, $y_{out},z^1,z^2,\cdots,z^p$ lie at $0,\zeta,\zeta^2,\cdots,\zeta^p$, where $\zeta=e^{2\pi i/p}$. $z^1,\cdots,z^p$ are called the \emph{distinguished inputs}. When we define ${}_pOC$ the boundary marked points of $\mathcal{R}^1_{k_1,\cdots,k_p}$ will be asymptotic to a $p$-fold Hochschild chain, with the $p$-distinguished marked points incident to the $p$ distinguished bimodule entries, hence the name.\par\indent
The codimension $1$ boundary components of $\overline{\mathcal{R}}^1_{k_1,\cdots,k_p}$ are of type
\begin{equation}
\mathcal{R}^{k'_i+1}\times \mathcal{R}^1_{k_1,\cdots,k_i-k_i'+1,\cdots,k_p},\quad 1\leq i\leq p,
\end{equation}
\begin{equation}
\mathcal{R}^{k_{i-1}'+k_{i}'+2}\times \mathcal{R}^1_{k_1,\cdots,k_{i-1}-k_{i-1}',k_{i}-k_{i}',\cdots,k_p},\quad 1\leq i\leq p.
\end{equation}
(2.32) corresponds to a disk bubble at the boundary where all the marked points on the bubble are non-distinguisehd; (2.33) corresponds to a disk bubble where one distinguised marked point (and some non-distinguished ones) bubbles off. \par\indent
Shorthand $\mathcal{F}$ for $\mathrm{Fuk}(M)_{\lambda}$. Given an element
\begin{equation}
\mathbf{x}=\mathbf{x}^1\otimes x^1_1\otimes\cdots\otimes x^1_{k_1}\otimes\cdots\otimes \mathbf{x}^p\otimes x^p_1\otimes \cdots\otimes x^p_{k_p}\in (\mathcal{F}_{\Delta}\otimes \mathcal{F}[1]^{\otimes k_1}\otimes\cdots\otimes \mathcal{F}_{\Delta}\otimes \mathcal{F}[1]^{\otimes k_p})^{\mathrm{diag}}\subset {}_pCC_*(\mathcal{F}),
\end{equation}
the coefficient of $y_{out}\in\mathrm{crit}(f)$ in
\begin{equation}
{}_pOC(\mathbf{x})\in CM^*(f)
\end{equation}
is given by counting rigid solutions to a parametrized moduli problem associated to $\mathcal{R}^1_{k_1,\cdots,k_p}$, with appropriate Lagrangian boundary conditions, the interior marked point constrained at $W^u(y_{out})$, and asymptotic conditions specified by $\mathbf{x}^1, x^1_1,\cdots, x^1_{k_1},\cdots, \mathbf{x}^p, x^p_1,\cdots x^p_{k_p}$.
\begin{prop}{\cite{Che}}
${}_pOC: {}_pCC_*(\mathrm{Fuk}(M)_{\lambda})\rightarrow CM^*(f)$ defines a chain map of degree $\frac{1}{2}\dim_{\mathbb{R}}M$.
\end{prop}
2.4. \textbf{The $\mathbb{Z}/p$-equivariant open-closed map}.
Similar to how $OC$ is enhanced to its $S^1$-equivariant version $OC^{S^1}$, ${}_pOC$ can be enhanced to the \emph{$\mathbb{Z}/p$-equivariant open-closed map}, which is a chain map
\begin{equation}
OC^{\mathbb{Z}/p}: CC_*^{\mathbb{Z}/p}(\mathrm{Fuk}(M)_{\lambda})\rightarrow CM^*(f)[[t,\theta]]
\end{equation}
of degree $\frac{1}{2}\dim_{\mathbb{R}}M$, where $CC_*^{\mathbb{Z}/p}$ is the negative $\mathbb{Z}/p$-equivariant Hochschild complex, cf. Definition 1.4. This chain map is the form
\begin{equation}
OC^{\mathbb{Z}/p}=\sum_{k\geq 0}(OC^{\mathbb{Z}/p}_{2k}+OC^{\mathbb{Z}/p}_{2k+1}\theta)t^k,
\end{equation}
for certain operations
\begin{equation}
OC^{\mathbb{Z}/p}_{2k}: {}_pCC_*(\mathrm{Fuk}(M)_{\lambda})\rightarrow CM^{*-\frac{1}{2}\dim_{\mathbb{R}}M-2k}(f)
\end{equation}
and
\begin{equation}
OC^{\mathbb{Z}/p}_{2k+1}: {}_pCC_*(\mathrm{Fuk}(M)_{\lambda})\rightarrow CM^{*-\frac{1}{2}\dim_{\mathbb{R}}M-2k-1}(f)
\end{equation}
of degree $0$. These operations are
defined via certain equivariant parametrized moduli spaces which we now recall in more detail.\par\indent
There is a $\mathbb{Z}/p$-action on
\begin{equation}
\coprod_{k_1,\cdots,k_p} \mathcal{R}^1_{k_1,\cdots,k_p}
\end{equation}
such that for $r\in\mathcal{R}^1_{k_1,\cdots,k_p}$, the standard generator $\tau\in\mathbb{Z}/p$ acts by counterclockwisely rotating the standard representative of $r$
by $e^{2\pi i/p}$ (hence $\tau(r)\in \mathcal{R}^1_{k_p,k_1,\cdots,k_{p-1}}$).
This action uniquely extends to the compactification $\coprod_{k_1,\cdots,k_p} \overline{\mathcal{R}}^1_{k_1,\cdots,k_p}$. Moreover, this $\mathbb{Z}/p$-action extends to the fiber bundle of Floer data over $\coprod_{k_1,\cdots,k_p} \overline{\mathcal{R}}^1_{k_1,\cdots,k_p}$, by pulling back the Floer
datum on $\Sigma_r$ to $\Sigma_{\tau(r)}$ along the $\mathbb{Z}/p$-rotation, denoted $\tau_r$. \par\indent
Consider the infinite sphere, regarded as a topological model for $E\mathbb{Z}/p$:
\begin{equation}
S^{\infty}:=\{w=(w_0,w_1,\cdots)\in\mathbb{C}^{\infty}: w_k=0\;\textrm{for}\;k\gg0, \|w\|^2=1\}.
\end{equation}
The free $\mathbb{Z}/p$-action on $S^{\infty}$ is given by
\begin{equation}
\tau(w_0,w_1,\cdots)=(e^{2\pi i/p} w_0,e^{2\pi i/p}w_1,\cdots).
\end{equation}
Consider the cells
\begin{equation}
\Delta_{2k}=\{w\in S^{\infty}: w_k\geq 0, w_{k+1}=w_{k+2}=\cdots=0\},
\end{equation}
\begin{equation}
\Delta_{2k+1}=\{w\in S^{\infty}: e^{i\theta}w_k\geq 0\;\textrm{for some}\;\theta\in[0,2\pi/p], w_{k+1}=w_{k+2}=\cdots=0\}.
\end{equation}
With suitably chosen orientation, one has
\begin{equation}
\partial\Delta_{2k}=\Delta_{2k-1}+\tau\Delta_{2k-1}+\cdots+\tau^{p-1}\Delta_{2k-1},
\end{equation}
\begin{equation}
\partial\Delta_{2k+1}=\tau\Delta_{2k}-\Delta_{2k}.
\end{equation}
It is proved in \cite{Che} that one can choose Floer data for the spaces $\overline{\mathcal{R}}^1_{k_1,\cdots,k_p}$ of disk, parametrized by $w\in S^{\infty}$, that satisfy the usual consistency conditions along the boundary strata
as well as $\mathbb{Z}/p$-equivariance. The latter condition means that, for a pair $(w,r)\in S^{\infty}\times \overline{\mathcal{R}}^1_{k_1,\cdots,k_p}$ and $\tau\in\mathbb{Z}/p$ the standard generator,
\begin{equation}
\epsilon_{\tau(w),r}=\epsilon_{w,\tau(r)}\circ \tau_r,\;\; J_{\tau(w),r}=J_{w,\tau(r)}\circ \tau_r,\;\;K_{\tau(w),r}=\tau_r^*K_{w,\tau(r)},
\end{equation}
where $(\epsilon, K,J)$ is a Floer datum on $\Sigma_r$, cf. section 2.1. We call such an $S^{\infty}$-dependent choice of Floer data \emph{a choice of Floer data for the $\mathbb{Z}/p$-equivariant open-closed map}. \par\indent
Then, $OC^{\mathbb{Z}/p}_{2k}$ is defined as follows. For $y_{out}\in \mathrm{crit}(f)$ and $\mathbf{x}=x^1\otimes x^1_1\otimes\cdots\otimes x^1_{k_1}\otimes\cdots\otimes x^p\otimes x^p_1\otimes \cdots\otimes x^p_{k_p}\in {}_pCC(\mathrm{Fuk}(M)_{\lambda})$, the $y_{out}$ coefficient of $OC^{\mathbb{Z}/p}_{2k}(\mathbf{x})$
is given by the count of rigid elements, weighted by monodromy along the disk boundary, of the parametrized moduli space
\begin{equation}
\mathcal{M}(\Delta_{2k}\times \mathcal{R}^1_{k_1,\cdots,k_p},y_{out},\mathbf{x})
\end{equation}
of tuples $(w,r,u)$, where $w\in \Delta_{2k}$, $r\in \mathcal{R}^1_{k_1,\cdots,k_p}$ and $u: \Sigma_r\rightarrow M$ satisfying Floer's equation (2.8); the boundary marked points
are asymptotic to $\mathbf{x}^1, x^1_1,\cdots, x^1_{k_1},\cdots, \mathbf{x}^p, x^p_1,\cdots x^p_{k_p}$; the boundary components lies on the corresponding Lagrangians; and the interior marked point is incident to the unstable manifold of $y_{out}$. $OC^{\mathbb{Z}/p}_{2k+1}$ is defined analogously using the parametrized moduli spaces
\begin{equation}
\mathcal{M}(\Delta_{2k+1}\times \mathcal{R}^1_{k_1,\cdots,k_p},y_{out},\mathbf{x}).
\end{equation}
\begin{prop}\cite{Che}
$OC^{\mathbb{Z}/p}$ is a chain map.
\end{prop}

\renewcommand{\theequation}{3.\arabic{equation}}
\setcounter{equation}{0}

\section{Cyclic and finite cyclic objects}
Connes' cyclic category $\Lambda$ is a mixture of the opposite simplex category $\Delta^{op}$ and the collection
of finite cyclic groups. The notion of \emph{cyclic objects}, i.e. functors $\Lambda\rightarrow \mathcal{C}$, gives a combinatorial way to package a circle action. A classical example
is the Hochschild homology of an associative algebra, cf. \cite{Con},\cite{Lod}, which we generalize to the case of an $A_{\infty}$-algebra. \par\indent
In 3.1, we start by reviewing the classical definition of $\Lambda$ due to Connes. Next, we introduce a new small combinatorial category ${}_p\Lambda$ together with a functor $j: {}_p\Lambda\rightarrow \Lambda$.\footnote{Note our ${}_p\Lambda$ is different from the $p$-cyclic category $\Lambda_p$ considered in e.g. \cite{Ka}, \cite{NS}.}
Analogous to the fact that functors out of $\Lambda$ model `objects with $S^1$-action', functors out of ${}_p\Lambda$ give a model for `objects with $\mathbb{Z}/p$-action'. Moreover, restriction along $j$ is the
`restriction of a circle action to the underlying $\mathbb{Z}/p\subset S^1$'. We also introduce variants of those categories $\overrightarrow{\Lambda},\overrightarrow{{}_p\Lambda}$, which are obtained from
$\Lambda,{}_p\Lambda$ respectively by removing the `degeneracy maps', and are useful for studying algebras that are not strictly unital. \par\indent
In 3.2, we introduce a dg version of the cyclic category without degeneracy maps $\overrightarrow{\Lambda}\rtimes \mathcal{A}_{\infty}^{dg}$ (and similarly $\overrightarrow{{}_p\Lambda}\rtimes \mathcal{A}_{\infty}^{dg}$).
We show that to an $\mathcal{A}_{\infty}$-algebra $A$ one can associate its \emph{Hochschild functor}, which is a dg functor $A^{\sharp}: \overrightarrow{\Lambda}\rtimes \mathcal{A}^{dg}_{\infty} \rightarrow\mathrm{Mod}_k$, that encodes the classical Hochschild invariants of $A$. \par\indent
In 3.3 we prove the main results of section 3, Proposition 3.11 (resp. Proposition 3.17), which shows that an explicit chain complex called the cyclic complex (resp.finite $p$-cyclic complex) of $A$
computes the negative cyclic homology (resp. negative finite $p$-cyclic homology) of $A$ defined in terms of an abstract $S^1$ (resp. $\mathbb{Z}/p$) homotopy fixed point. Proposition 3.11 is an adaptation of \cite[Theorem 2.3]{Hoy}. \par\indent
Finally, in 3.4 we briefly discuss cocyclic and finite cocyclic objects.\par\indent
3.1. \textbf{Connes' cyclic category} $\Lambda$. We give a geometric description of $\Lambda$ following \cite{Con}, \cite[E.6.1.2]{Lod}:
\begin{itemize}
    \item Objects of $\Lambda$ are labeled by natural numbers $[n], n\geq 0$. One thinks of $[n]$ as a configuration of $n+1$ marked points $z_0<\cdots<z_n$ on a circle.
    \item Morphisms from $[n]$ to $[m]$ are homotopy classes of degree $1$ nondecreasing maps from $S^1$ to $S^1$ that send marked points to marked points.
\end{itemize}
Combinatorially, a morphism from $[n]$ to $[m]$ is uniquely determined by
\begin{itemize}
    \item a cyclic reordering of $[n]$ given by $[n]_{\sigma}=\{\sigma_0<\sigma_1<\cdots<\sigma_n\}$, where $\sigma\in \frac{\mathbb{Z}}{n+1}\subset S_{n+1}$ and
    \item  an (ordered) partition of $[n]_{\sigma}$ into $m+1$ subsets, i.e. (possibly empty) order subsets $f_0,\cdots,f_m\subset [n]_{\sigma}$ such that $[n]_{\sigma}=f_0\star f_1\star\cdots\star f_m$, where
$\star$ denotes the join of partially ordered sets.
\end{itemize}
To make the identification, given a homotopy class of $f:S^1\rightarrow S^1$, let $f_i$ be the set of (indices of) marked points that are sent to $z_i$ by $f$. By abuse of notation, we also write the ordered set $f_i$ as $f^{-1}(i)$.\par\indent
Schematically, we represent a morphism $f: [n]\mapsto [m]$ by drawing a circle with $m+1$ ordered marked points (thought of as the object $[m]$), and write the elements of $f^{-1}(i)$ in order next to the $i$-th marked point. There is an inclusion $i:\Delta^{op}\subset\Lambda$ as the subcategory whose morphisms preserve the zeroth marked point. For instance, the face maps $d_i: [n]\rightarrow [n-1], 0\leq i\leq n$ and degeneracy maps $s_i:[n-1]\rightarrow [n], 0\leq i\leq n-1$ are presented in Figure 2.
\begin{figure}[H]
 \centering
 \includegraphics[width=0.6\textwidth]{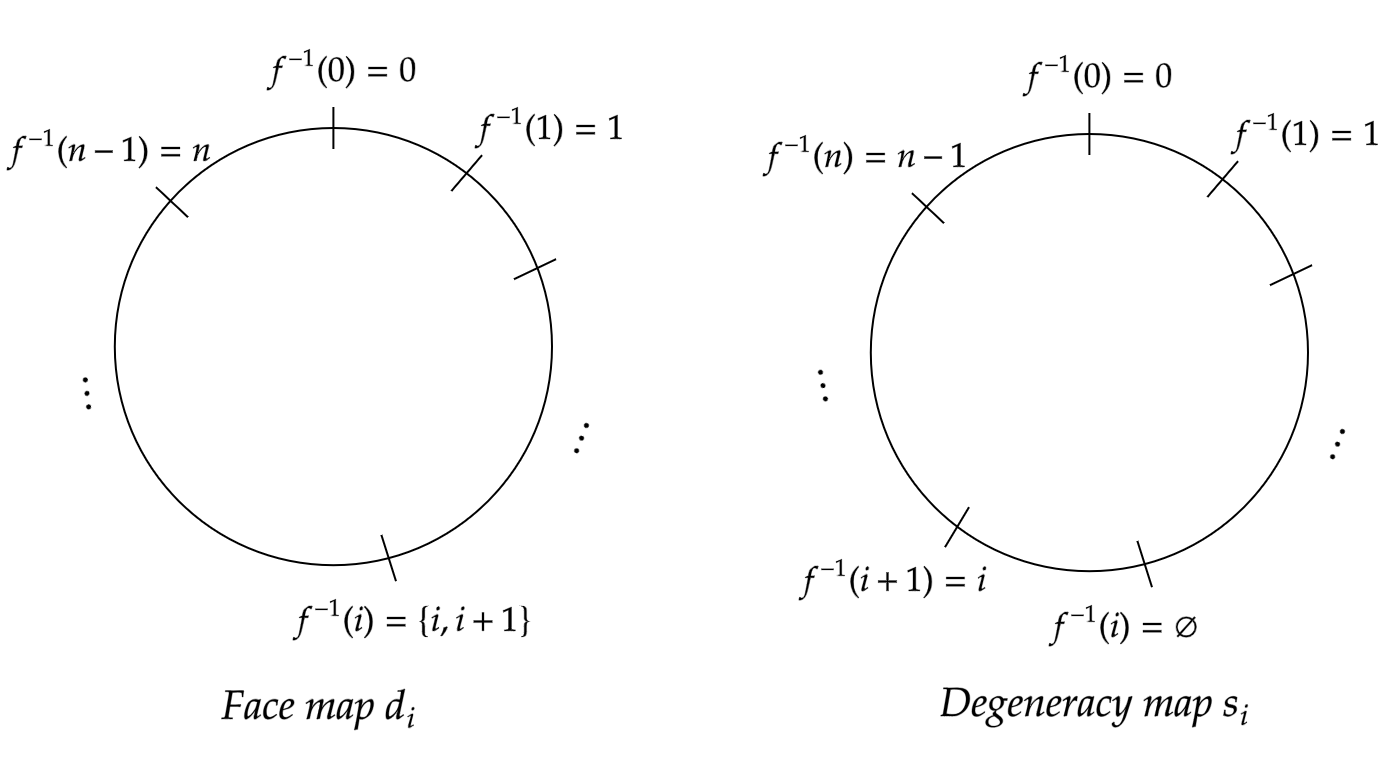}
 \caption{}
\end{figure}
Let $\tau\in \frac{\mathbb{Z}}{n+1}$ be the standard generator. Note there is a natural identification $\mathrm{Aut}_{\Lambda}([n])=\frac{\mathbb{Z}}{n+1}$, and any morphism $f\in \Lambda([n],[m])$ can be uniquely expressed as $f=\tau^k\circ g$ for some $k$ and $g\in \Delta^{op}([n],[m])$. \par\indent
\textbf{The finite $p$-cyclic category ${}_p\Lambda$}. Let $p$ be a prime number. The \emph{finite $p$-cyclic category} ${}_p\Lambda$ has
\begin{itemize}
    \item Objects $p$-tuples of natural numbers $[k_1,\cdots,k_p]$, thought of as a configuration of $k_1+\cdots+k_p+p$ counterclockwise marked points on the circle, with the $0, k_1+1, k_1+k_2+2,\cdots,k_1+\cdots+k_{p-1}+p-1$-th points marked as distinguished.
    \item Morphisms from $[k_1,\cdots,k_p]$ to $[k_1',\cdots,k_p']$ are homotopy classes of degree $1$ nondecreasing maps $f: S^1\rightarrow S^1$ that sends marked points to marked points, and is further required to send the distinguished marked points bijectively to distinguished marked points. In particular, by the nondecreasing condition, $f$ must act as a $\mathbb{Z}/{p}$ cyclic permutation on the distinguished points.
\end{itemize}
There is an obvious functor $j:{}_p\Lambda\rightarrow \Lambda$ which on objects sends $[k_1,\cdots,k_p]$ to $[k_1+\cdots+k_p+p-1]$, where the $p$ distinguished marked points becomes ordinary marked points. Moreover, $i_p:(\Delta^{op})^p\subset {}_p\Lambda$ sits as the subcategory whose morphisms fix each distinguished marked point. These categories fit into a diagram (with the left vertical arrow being the ordinal sum, cf Appendix C)
\begin{equation}
\begin{tikzcd}[row sep=1.2cm, column sep=0.8cm]
\Delta^{op}\arrow[r,"i"]& \Lambda\\
(\Delta^{op})^p\arrow[r,"{i_p}"]\arrow[u,"o"]& {}_p\Lambda\arrow[u,"j"]
\end{tikzcd}.
\end{equation}
\textbf{The variants $\overrightarrow{\Lambda}$ and $\overrightarrow{{}_p\Lambda}$}. We define $\overrightarrow{\Lambda}\subset\Lambda$ to be the subcategory with the same set of objects as $\Lambda$, but morphisms the homotopy classes of $f:S^1\rightarrow S^1$ that are surjective on marked points (i.e. no `degeneracy maps'). Similarly, we define $\overrightarrow{\Delta^{op}}, \overrightarrow{{}_p\Lambda}$ and $\overrightarrow{(\Delta^{op})^p}$ by removing the morphisms that are not surjective
on marked points.\par\indent
3.2. Fix a symmetric monoidal category $(\mathcal{V},\otimes)$ with finite coproducts, and a non-symmetric operad $\{Q_j\}_{j\geq 1}$ valued in $\mathcal{C}$. The cases of interest
to us are $\mathcal{V}=\mathrm{Top}$ or $\mathrm{Ch}_k$. Here $\mathrm{Ch}_k$ denotes the category of chain complexes over $k$ equipped with the standard symmetric monoidal tensor product. To make the distinction, we will denote $\mathrm{Mod}_k$ for the dg catgeory of chain complexes over $k$.
\begin{mydef}
Define $\overrightarrow{\Lambda}\rtimes Q$ to be the following $\mathcal{V}$-enriched category.
\begin{itemize}
    \item Objects are given by $[n], n\geq 0$.
    \item The morphism object in $\mathcal{V}$ from $[n]$ to $[m]$ is defined to be
\begin{equation}
\coprod_{f\in \overrightarrow{\Lambda}([n],[m])} \bigotimes_{i=0}^m Q_{|f^{-1}(i)|}.
\end{equation}
\end{itemize}
Note that by definition of $\overline{\Lambda}$, one has $|f^{-1}(i)|>0$. Composition combines the composition of maps in $\overrightarrow{\Lambda}$ with the operadic structure of $Q$. More explicitly, for $f\in \overrightarrow{\Lambda}([n],[m]), g\in \overrightarrow{\Lambda}([m],[r])$, we define the $(g,f)$-component of the composition $\mathrm{Hom}_{\overrightarrow{\Lambda}\rtimes Q}([n],[m])\otimes \mathrm{Hom}_{\overrightarrow{\Lambda}\rtimes Q}([m],[r])\rightarrow \mathrm{Hom}_{\overrightarrow{\Lambda}\rtimes Q}([n],[r])$ to be
the map
\begin{equation}
\bigotimes_{j=0}^r Q_{|g^{-1}(j)|}\otimes\bigotimes_{i=0}^m Q_{|f^{-1}(i)|}\rightarrow \bigotimes_{j=0}^r Q_{|(g\circ f)^{-1}(j)|}=\bigotimes_{j=0}^r Q_{|f^{-1}(i^j_1)|+\cdots+|f^{-1}(i^j_{|g_j|})|}
\end{equation}
induced by the operadic compositions ($0\leq j\leq r$),
\begin{equation}
\circ: Q_{|g^{-1}(j)|}\otimes \bigotimes_{k=1}^{|g^{-1}(j)|}Q_{|f^{-1}(i^j_k)|}\rightarrow Q_{|f^{-1}(i^j_1)|+\cdots+|f^{-1}(i^j_{|g_j|})|},
\end{equation}
where $g^{-1}(j)=\{i^j_1<i^j_2<\cdots<i^j_{|g_j|}\}\subset [m]$. The ordering of the big monoidal product in (3.4) should agree with the ordering of $g^{-1}(j)$. \par\indent
Similar, one can define $\overrightarrow{{}_p\Lambda}\rtimes Q, \overrightarrow{\Delta^{op}}\rtimes Q$, etc.
\end{mydef}
\emph{Example 1: The Associahedron}. Let $\overline{\mathcal{R}}^{d+1}, d\geq 2$ be the compactified space of disks with $d$ boundary inputs and $1$ boundary output from section 2.1. We set $(\mathcal{A}_{\infty})_d:=\overline{\mathcal{R}}^{d+1}, d\geq 2$ and $(\mathcal{A}_{\infty})_1:=pt$. Then, the collection of spaces
$\mathcal{A}_{\infty}=\{(\mathcal{A}_{\infty})_d\}_{d\geq 1}$ form a topological operad known as the \emph{Stasheff associahedron}. The operadic structures involving $(\mathcal{A}_{\infty})_d, d\geq 2$ are given by concatenation of disks, and $(\mathcal{A}_{\infty})_1=\mathrm{pt}$ acts as the operadic unit. \par\indent
By (2.7), there is a natural cellular structure on
$\overline{\mathcal{R}}^{d+1}$ where the cells are given by $\mathcal{R}^T$ for some planar stable $d$-leafed tree $T$. With respect to these cellular structures (ranging over $d$), the operadic composition maps of $\mathcal{A}_{\infty}$ is cellular, as concatenation of (stable) disks respects the concatenation of
the underlying tree types. Let $\mathcal{A}^{dg}_{\infty}:=\{C_{-*}^{cell}((A_{\infty})_d)\}_{d\geq 1}$ be the dg operad of cellular chains on the associahedron. We call this dg operad the \emph{dg associahedron}.\par\indent
Let $\mu^d\in C^{cell}_{d-2}((\mathcal{A}_{\infty})_d)=C^{cell}_{d-2}(\overline{\mathcal{R}}^{d+1})$ be the top dimensional cell (fixing a preferred orientation) corresponding to the unique $d$-leafed tree with one internal vertex. An algebra over the dg operad
$\mathcal{A}^{dg}_{\infty}$ is easily seen to recover the usual notion of a non-unital $A_{\infty}$-algebra, with the higher structure maps given by the operations corresponding to $\mu_d$'s (the notations are intentionally made to agree). \par\indent 
Following Definition 3.1, we can form the dg category $\overrightarrow{\Lambda}\rtimes \mathcal{A}^{dg}_{\infty}$. Let $A$ be an $k$-linear $\mathcal{A}_{\infty}$-algebra.
\begin{mydef}
The \emph{Hochschild functor} of $A$ is the dg functor
\begin{equation}
A^{\sharp}: \overrightarrow{\Lambda}\rtimes \mathcal{A}^{dg}_{\infty}\rightarrow \mathrm{Mod}_k
\end{equation}
given by
\begin{itemize}
    \item On objects, $A^{\sharp}([n]):=(A,\mu^1)^{\otimes n+1}$ as chain complex.
    \item On morphisms, $A^{\sharp}$ sends $(f\in \overrightarrow{\Lambda}([n],[m]),\otimes_{i=0}^m \phi_i\in \bigotimes_{i=0}^m (\mathcal{A}^{dg}_{\infty})_{|f^{-1}(i)|})$ to the map $A^{\otimes n+1}\rightarrow A^{\otimes m+1}$ given by
\begin{equation}
a_0\otimes a_1\otimes \cdots\otimes a_n\mapsto \phi_0(\otimes_{j\in f^{-1}(0)}a_j)\otimes \phi_1(\otimes_{j\in f^{-1}(1)}a_j)\otimes\cdots\otimes \phi_m(\otimes_{j\in f^{-1}(m)}a_j).
\end{equation}
\end{itemize}
Here, the tensor product $\otimes_{j\in f^{-1}(i)}a_j$ is taken with respect to the canonical ordering on $f^{-1}(i)$.
\end{mydef}
Now we describe a way to reproduce the standard Hochschild and cyclic homology chain complex from $A^{\sharp}$, in a way that generalizes to all dg functors $\overrightarrow{\Lambda}\rtimes \mathcal{A}^{dg}_{\infty}\rightarrow \mathrm{Mod}_k$. First we fix some notations.
\begin{itemize}
    \item We denote the unique element of $(\mathcal{A}_{\infty})_1=pt$, which becomes a generator of $(\mathcal{A}_{\infty}^{dg})_1=C^{cell}_0((\mathcal{A}_{\infty})_1)$, as $\mathrm{id}$.
    \item Recall that the standard generator $\tau\in\frac{\mathbb{Z}}{n+1}$ gives an automorphism of $[n]$ in $\overrightarrow{\Lambda}$. By abuse of notation, we also use $\tau$ to denote the morphism
\begin{equation}
(\tau,\mathrm{id}\otimes \mathrm{id}\otimes\cdots\otimes\mathrm{id})\in \mathrm{Hom}_{\overrightarrow{\Lambda}\rtimes\mathcal{A}_{\infty}^{dg}}([n],[n]).
\end{equation}
\end{itemize}
\begin{mydef}
A dg functor $Q: \overrightarrow{\Lambda}\rtimes \mathcal{A}^{dg}_{\infty}\rightarrow \mathrm{Mod}_k$ is called a \emph{non-unital} $\mathcal{A}_{\infty}$-\emph{cyclic $k$-module}.
\end{mydef}
To distinguish notation from shifts, we denote $Q([n])$ as $Q_n$ and denote $d_{Q_n}$ its differential.
\begin{mydef}
Let $Q$ be a non-unital $\mathcal{A}_{\infty}$-cyclic $k$-module. As a graded $k$-vector space, the \emph{cyclic bar complex} of $Q$ is
\begin{equation}
CC(Q):=\bigoplus_{n\geq 0} Q_n[n],
\end{equation}
where $[n]$ denotes shifting the degree by $n$ (i.e. $X[n]_k=X_{n+k}$). We define the following operations on $CC(Q)$:
\begin{itemize}
    \item The \emph{bar differential} $b'_Q$ is the degree $1$ differential on $CC(Q)$ defined by: for $x\in Q_n[n]\subset CC(Q)$,
\begin{equation}
b'_Q(x):=d_{Q_n}    (x)+\sum_{m<n}\sum_{\substack{f\in \overrightarrow{\Lambda}([n],[m]): \min f^{-1}(0)=0\\\textrm{and}\;\exists!\,j\in[m]\;\textrm{s.t.}\;|f^{-1}(j)|>1}}  (-1)^{\dag} Q(\mathrm{id}\otimes\mathrm{id}\otimes\cdots\otimes \mu_{|f^{-1}(j)|}\otimes \cdots\otimes\mathrm{id})(x).
\end{equation}
The sign is given by $\dag=r+|f^{-1}(j)|t$, where $r, t$ are the number of $\mathrm{id}$'s before and after the term $\mu_{|f^{-1}(j)|}$, respectively, in (3.9).
     \item The \emph{cyclic bar differential} $b_Q$ is the degree $1$ differential on $CC(Q)$
\begin{equation}
b_Q:=b'_Q+w,
\end{equation}
where $w$ is the \emph{wrapped around term} defined by: for $x\in Q_n[n]\subset CC(Q)$,
\begin{equation}
w(x):=\sum_{m<n}\sum_{\substack{f\in\overrightarrow{\Lambda}([n],[m]): 0\in f^{-1}(0)\backslash\{\min f^{-1}(0)\}\\\textrm{and}\;|f^{-1}(j)|=1\;\textrm{for all}\;j\neq 0}} (-1)^{\ddag} Q(\mu_{|f^{-1}(0)|}\otimes \mathrm{id}\otimes\cdots\otimes \mathrm{id})(x).
\end{equation}
The sign is $\ddag$ is defined as follows. Let $f^{-1}(0)=\{i_{-r}<\cdots<i_{-1}<0<i_1<\cdots<i_t\}$, then $\ddag=n+r+t+tn$.
\end{itemize}
$(CC(Q),b_Q)$ is called the \emph{cyclic bar complex} of $Q$. When $Q=A^{\sharp}$ for an $A_{\infty}$-algebra $A$, $b'_Q$ and $b_Q$ are exactly the bar differential and cyclic bar (Hochschild) differential on $CC(A)$; cf. Definition 1.2.
\end{mydef}
\begin{mydef}
The \emph{negative cyclic complex} $CC^{S^1}(Q)$ is the product totalization of
\begin{equation}
0\rightarrow (CC(Q),b'_Q)\xrightarrow{\tau-1} (CC(Q),b_Q)\xrightarrow{N} (CC(Q),b'_Q)\xrightarrow{\tau-1} \cdots,
\end{equation}
where it is understood that on $Q_d[d]\subset CC(Q)$, $\tau$ acts as the standard generator of $\frac{\mathbb{Z}}{d+1}$, and $N$ acts as $1+\tau+\cdots+\tau^{d}$. Alternatively, we can write the negative cyclic complex as
\begin{equation}
CC(Q)[[u,e^+]],
\end{equation}
where $|u|=2, |e^+|=-1, (e^+)^2=0$, and the $u$-linear differential is given by
\begin{equation}
\begin{cases}
x\mapsto b_Q(x)+N(x)ue^+\\
xe^+\mapsto b'_Q(x)e^++(\tau-1)x.
\end{cases}
\end{equation}
When $Q=A^{\sharp}$ for a (non-unital) $\mathcal{A}_{\infty}$-algebra $A$, $CC^{S^1}(A^{\sharp})$ is isomorphic to $CC^{S^1}(A)$ from (2.15).
\end{mydef}

\textbf{Construction}. We briefly discuss an enhancement of the previous constructions when the $\mathcal{A}_{\infty}$-algebra at hand is \emph{strictly unital}. To incorporate the unit operadically, we can define an augmented nonsymmetric operad (i.e. a nonsymmetric operad with $0$-ary operations), \emph{the unital dg associahedron}, $\mathcal{A}_{\infty}^{dg,u}=\{(\mathcal{A}_{\infty}^{dg,u})_i\}_{i\geq 0}$ by $(\mathcal{A}_{\infty}^{dg,u})_i=(\mathcal{A}_{\infty}^{dg})_i$ if $i>0$ and $(\mathcal{A}_{\infty}^{dg,u})_0=k\mathbf{1}$,
where $\mathbf{1}$ is an $0$-ary operation acting as `inserting a strict unit'. More precisely, the operadic structure maps of $\mathcal{A}_{\infty}^{dg,u}$ are given by the operadic structure maps of $\mathcal{A}_{\infty}^{dg}$ if the element $\mathbf{1}$ is not involved. When $\mathbf{1}$ is involved, we define the operadic
structure map as follows. If $T$ is a rooted planar tree representing a cell in
$C_*^{cell}(\overline{\mathcal{R}}^{d+1})$, we define $T\circ_i \mathbf{1}=0$ unless the $i$-th leaf of $T$ is connected to a trivalent vertex, in which case the rule is
\begin{figure}[H]
 \centering
 \includegraphics[width=1\textwidth]{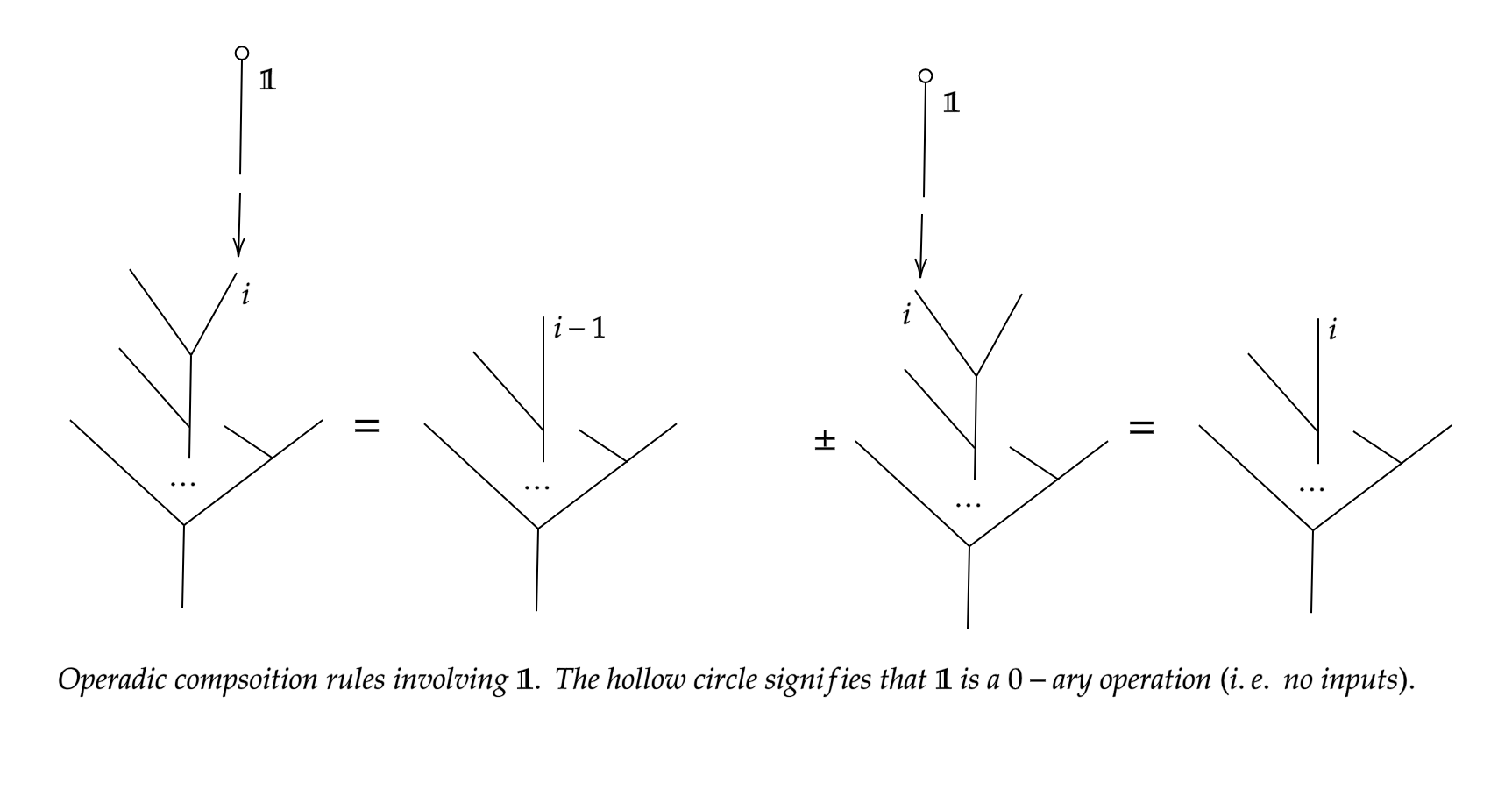}
 \caption{}
\end{figure}
It is easy to see that an algebra over the dg operad $\mathcal{A}_{\infty}^{dg,u}$ recovers the notion of \emph{a strictly unital $\mathcal{A}_{\infty}$-algebra}, cf. \cite[section (2a)]{Sei1}.\par\indent
Since $\mathcal{A}_{\infty}^{dg,u}$ contains $0$-ary operations, we can define a dg category
\begin{equation}
\Lambda\rtimes \mathcal{A}_{\infty}^{dg,u}
\end{equation}
analogous to Definition 3.1, but in contrast, allowing maps $f:S^1\rightarrow S^1\in \Lambda([m],[n])$ that are not surjective onto marked points (and we insert $0$-ary operations at those skipped points).\par\indent 
In light of the above discussion, we make the following definition.
\begin{mydef}
A dg functor
\begin{equation}
\Lambda\rtimes \mathcal{A}_{\infty}^{dg,u}\rightarrow \mathrm{Mod}_k
\end{equation}
is called a \emph{strictly unital $\mathcal{A}_{\infty}$-cyclic $k$-module}.
\end{mydef}
One of key properties of a strictly unital $\mathcal{A}_{\infty}$-cyclic $k$-module is the following.
\begin{lemma}
Let $Q$ be a strictly unital $\mathcal{A}_{\infty}$-cyclic $k$-module, then its bar complex $(CC(Q),b')$ is acyclic.
\end{lemma}
\noindent\emph{Proof}. The extra degeneracy map provides a contracting homotopy for the bar complex, cf. the proof of \cite[1.1.12]{Lod}.\qed\par\indent
This property is formalized by the following notion due to Wodzicki.
\begin{mydef}
$Q\in \mathrm{Fun}^{dg}(\overrightarrow{\Lambda}\rtimes \mathcal{A}^{dg}_{\infty},\mathrm{Mod}_k)$ is called \emph{H-unital} if its bar complex $(CC(Q),b'_Q)$ is acyclic.
\end{mydef}

In practice, strict units rarely arise from geometric constructions. For instance, it is well known that the Fukaya category only has a unit at the cohomological level (which can be enhanced to stronger notion of a homotopy unit), cf. \cite[section 10]{Gan1}.
Luckily, having a cohomological unit suffices because of the following lemma.
\begin{lemma}
If $A$ is a cohomologically unital $\mathcal{A}_{\infty}$-algebra, then $A^{\sharp}$ is H-unital.
\end{lemma}
\noindent\emph{Proof}. See the proof of \cite[Prop 2.2]{Gan1}.\qed\par\indent
3.3. \textbf{Negative cyclic homology as homotopy fixed point}. We adapt the results of \cite{Hoy} to the $k$-linear $\mathcal{A}_{\infty}$ setting, which allow one to reinterpret the chain complexes $CC(Q), CC^{S^1}(Q)$ as computing certain abstract homotopy (co)limits.  \par\indent
We fix some notations for functors among the combinatorial categories introduced in section 3.1. For a category $\mathcal{C}$, let $\mathcal{C}_k$ be the free $k$-linear category associated to $\mathcal{C}$, viewed as a dg category whose morphism spaces are concentrated in degree $0$.
\begin{itemize}
    \item ($i$): Recall from section 3.1 that there is an inclusion of categories $i:\Delta^{op}\subset \Lambda$. This indu ces an inclusion of categories $\overrightarrow{i}: \overrightarrow{\Delta}^{op}\subset  \overrightarrow{\Lambda}$, as well as
inclusions of dg categories $\overrightarrow{i_k}:(\overrightarrow{\Delta}^{op})_k\subset \overrightarrow{\Lambda}_k$ and $\overrightarrow{i_{\mathcal{A}_{\infty}^{dg}}}:\overrightarrow{\Delta}^{op}\rtimes\mathcal{A}_{\infty}^{dg}\subset \overrightarrow{\Lambda}\rtimes \mathcal{A}_{\infty}^{dg}$, etc.
      \item ($i_p$): Analogously, there is an inclusion of categories $i_p:(\Delta^{op})^p\subset {}_p\Lambda$, which induces an inclusion of categories $\overrightarrow{i_p}: (\overrightarrow{\Delta}^{op})^p\subset  \overrightarrow{{}_p\Lambda}$, as well as
inclusions of dg categories $(\overrightarrow{i_p})_k:(\overrightarrow{\Delta}^{op})^p_k\subset (\overrightarrow{{}_p\Lambda})_k$ and $(\overrightarrow{i_p})_{\mathcal{A}_{\infty}^{dg}}:(\overrightarrow{\Delta}^{op})^p\rtimes\mathcal{A}_{\infty}^{dg}\subset \overrightarrow{{}_p\Lambda}\rtimes \mathcal{A}_{\infty}^{dg}$, etc.
      \item ($P$): There is a dg functor $P_{\overrightarrow{\Lambda}}: \overrightarrow{\Lambda}\rtimes \mathcal{A}_{\infty}^{dg}\rightarrow \overrightarrow{\Lambda}_k$ which is the identity on objects, and on morphism spaces is induced by the augmentation maps $\varepsilon: C_{-*}^{cell}(\overline{\mathcal{R}}^{d+1})\rightarrow k$.
Similarly, we have dg functors $P_{\overrightarrow{\Delta}^{op}}:\overrightarrow{\Delta}^{op}\rtimes \mathcal{A}_{\infty}^{dg}\rightarrow (\overrightarrow{\Delta}^{op})_k$,
$P_{\overrightarrow{{}_p\Lambda}}:\overrightarrow{{}_p\Lambda}\rtimes \mathcal{A}_{\infty}^{dg}\rightarrow (\overrightarrow{{}_p\Lambda})_k$, etc. Since $\overline{\mathcal{R}}^{d+1}$ is contractible, each such $P$ is an equivalence of dg categories.
       \item ($\iota$): Let $\iota_{\Lambda}:\overrightarrow{\Lambda}\subset \Lambda$ denote the inclusion of the subcategory $\overrightarrow{\Lambda}$. Similarly, we have
$\iota_{\Delta^{op}}:\overrightarrow{\Delta}^{op}\subset \Delta^{op},\iota_{{}_p\Lambda}:\overrightarrow{{}_p\Lambda}\subset {}_p\Lambda$, etc.
       \item ($j$): Recall from section 3.1 that there is a functor $j: {}_p\Lambda\rightarrow \Lambda$. This induces a functor $\overrightarrow{j}:\overrightarrow{{}_p\Lambda}\rightarrow \overrightarrow{\Lambda}$ and dg functors
$\overrightarrow{j_k}: (\overrightarrow{{}_p\Lambda})_k\rightarrow \overrightarrow{\Lambda}_k$ and $\overrightarrow{j_{\mathcal{A}^{dg}_{\infty}}}: \overrightarrow{{}_p\Lambda}\rtimes\mathcal{A}_{\infty}^{dg}\rightarrow \overrightarrow{\Lambda}\rtimes \mathcal{A}_{\infty}^{dg}$, etc.
       \item ($\mathfrak{i},\mathfrak{i}_p,\mathfrak{j}$): Finally, we let $\mathfrak{i}:*\rightarrow B\mathbb{T}$ and $\mathfrak{i}_p:*\rightarrow B\mathbb{Z}/p$ denote the inclusion of a basepoint, and $\mathfrak{j}:\mathbb{Z}/p\rightarrow \mathbb{T}$ the standard inclusion as the $p$-th roots of unity.
\end{itemize}
To build on certain results regarding Hochschild and cyclic homology in \cite{Hoy}, we need to use the language of $\infty$-categories (see Appendix D for a brief review). On the other hand, in section 3.2 we
set up the theory of Hochschild and cyclic homology in the context of $k$-linear dg/$\mathcal{A}_{\infty}$ categories. In particular, given a non-unital $\mathcal{A}_{\infty}$ $k$-cyclic module $X$, we would like to obtain from it a non-unital cyclic $k$-module in the $\infty$-categorical sense, i.e. a functor
of $\infty$-categories $N\overrightarrow{\Lambda}\rightarrow Mod_k$. This can be done in two steps:\par\indent
Step 1. We `pushforward' the dg functor $X: \overrightarrow{\Lambda}\rtimes \mathcal{A}_{\infty}^{dg}\rightarrow \mathrm{Mod}_k$ along $P_{\overrightarrow{\Lambda}}:\overrightarrow{\Lambda}\rtimes \mathcal{A}_{\infty}^{dg}\rightarrow \overrightarrow{\Lambda}_k$. More precisely, since $P_{\overrightarrow{\Lambda}}$ is a strictly unital dg quasi-equivalence (in particular an $\mathcal{A}_{\infty}$ quasi-equivalence),
it has a strictly unital $\mathcal{A}_{\infty}$ quasi-inverse $P_{\overrightarrow{\Lambda}}^{-1}$, which is unique up to $\mathcal{A}_{\infty}$-homotopy. Let $\mathrm{Fun}^{\mathcal{A}_{\infty},u}(\overrightarrow{\Lambda}\rtimes \mathcal{A}_{\infty}{dg},\mathrm{Mod}_k)$ denote the category of strictly unital $\mathcal{A}_{\infty}$-functors from $\overrightarrow{\Lambda}\rtimes \mathcal{A}_{\infty}^{dg}$ to $\mathrm{Mod}_k$; equivalently, it is the category of strictly unital left $\mathcal{A}_{\infty}$-modules over $\overrightarrow{\Lambda}\rtimes \mathcal{A}_{\infty}^{dg}$, cf. \cite[section (2f)]{Sei1}.
As a result, pulling back along $P_{\overrightarrow{\Lambda}}$ and $P_{\overrightarrow{\Lambda}}^{-1}$ gives rise to
an equivalence of categories
\begin{equation}
(P_{\overrightarrow{\Lambda}}^{-1})^*:\mathrm{Fun}^{\mathcal{A}_{\infty},u}(\overrightarrow{\Lambda}\rtimes \mathcal{A}_{\infty}^{dg}, \mathrm{Mod}_k)\rightleftarrows \mathrm{Fun}^{\mathcal{A}_{\infty},u}(\overrightarrow{\Lambda}_k, \mathrm{Mod}_k): P_{\overrightarrow{\Lambda}}^*.
\end{equation}
Given $X\in \mathrm{Fun}^{dg}(\overrightarrow{\Lambda}\rtimes \mathcal{A}_{\infty}^{dg}, \mathrm{Mod}_k)$, we simply denote by $(P^{-1}_{\overrightarrow{\Lambda}})^*X$ the strictly unital $\mathcal{A}_{\infty}$-functor given by the image of $X$ under the composition
\begin{equation}
\mathrm{Fun}^{dg}(\overrightarrow{\Lambda}\rtimes \mathcal{A}_{\infty}^{dg}, \mathrm{Mod}_k)\subset\mathrm{Fun}^{\mathcal{A}_{\infty},u}(\overrightarrow{\Lambda}\rtimes \mathcal{A}_{\infty}^{dg}, \mathrm{Mod}_k)\xrightarrow{(P^{-1}_{\overrightarrow{\Lambda}})^*} \mathrm{Fun}^{\mathcal{A}_{\infty},u}(\overrightarrow{\Lambda}_k, \mathrm{Mod}_k).
\end{equation}
Step 2. Apply Proposition D.3 to $(P^{-1}_{\overrightarrow{\Lambda}})^*X$, we obtain a functor of $\infty$-categories
\begin{equation}
((P^{-1}_{\overrightarrow{\Lambda}})^*X)_{\Delta}:N\overrightarrow{\Lambda}\rightarrow N^{dg}\mathrm{Mod}_k\rightarrow Mod_k,
\end{equation}
where $Mod_k$ denote the derived $\infty$-category of chain complexes over $k$, cf. Definition D.4 (note the second arrow of (3.19) is an equivalence since $k$ is a field).
To simply notation, we make this identification implicit and write $(P^{-1}_{\overrightarrow{\Lambda}})^*X$ for $((P^{-1}_{\overrightarrow{\Lambda}})^*X)_{\Delta}\in \mathrm{sSet}(N\overrightarrow{\Lambda},Mod_k)$. If $\mathcal{C},\mathcal{D}$ are $\infty$-categories, we write $\mathrm{Fun}(\mathcal{C},\mathcal{D})$ for their mapping $\infty$-category. In particular, $\mathrm{Fun}(\mathcal{C},\mathcal{D})_0=\mathrm{sSet}(\mathcal{C},\mathcal{D})$.\par\indent
Let $\mathbb{T}=\mathrm{Sing}_*(S^1)$ be the simplicial circle. Up to homotopy equivalence, it is uniquely characterized as a $K(\mathbb{Z},1)$. Let $N\Lambda\rightarrow \tilde{N\Lambda}$ be the $\infty$-groupoid completion of $N\Lambda$; in other words, $\tilde{N\Lambda}$ is up to homotopy equivalence the unique Kan complex equipped with a weak equivalence from $N\Lambda$. By \cite[Corollary 1.2]{Hoy}, there is a homotopy
equivalence $\tilde{N\Lambda}\simeq B\mathbb{T}$ and an equivalence of $\infty$-groups $\mathbb{T}\simeq \mathrm{Aut}_{[0]}(\tilde{N\Lambda})$.
Let $\mathfrak{i}:*\rightarrow B\mathbb{T}$ denote the inclusion of a basepoint and $q:B\mathbb{T}\rightarrow *$ the unique map to the terminal object. If $\mathcal{E}$ is any presentable symmetric monoidal $\infty$-category (e.g. $Mod_k$), the left and right Kan extension (cf. Appendix D.2) along $q$  will be denoted by $(-)_{h\mathbb{T}}, (-)^{h\mathbb{T}}: \mathrm{Fun}(B\mathbb{T},\mathcal{E})\rightarrow \mathcal{E}$, and called the \emph{$\mathbb{T}$-homotopy orbit} and \emph{$\mathbb{T}$-homotopy fixed point} functor, respectively.
\begin{lemma}
The diagram of $\infty$-categories
\begin{equation}
\begin{tikzcd}[row sep=1.2cm, column sep=0.8cm]
N\Delta^{op}\arrow[r,"i"]\arrow[d]&N\Lambda\arrow[d]\\
*\arrow[r,"\mathfrak{i}"]& B\mathbb{T}
\end{tikzcd}
\end{equation}
is a homotopy exact square (cf. Definition D.13).
\end{lemma}
\noindent\emph{Proof}. This follows from Proposition 1.1 of \cite{Hoy}.\qed\par\indent
Let $\int_{\mathbb{T}}: \mathrm{Fun}(N\Lambda,Mod_k)\rightarrow  \mathrm{Fun}(B\mathbb{T},Mod_k)$ denote the left Kan extension along $N\Lambda\rightarrow B\mathbb{T}$; let $\int^{\mathbb{T}}$ denote the right Kan extension along the same functor. By Lemma 3.10, for a functor
$Y:N\Lambda\rightarrow Mod_k$, the underlying chain complex of $\int_{\mathbb{T}}Y$ (i.e. $\mathfrak{i}^*\int_{\mathbb{T}} Y$) is quasi-isomorphic to $\mathrm{colim}_{N\Delta^{op}}i^*Y$.   \par\indent
We now state the main result of this subsection.
\begin{prop}
Let $X$ be a non-unital $\mathcal{A}_{\infty}$-cyclic $k$-module. Then there is a natural quasi-isomorphism
\begin{equation}
CC^{S^1}(X)\simeq \big(\int_{\mathbb{T}}(\iota_{\Lambda})_!(P^{-1}_{\overrightarrow{\Lambda}})^*X\big)^{h\mathbb{T}},
\end{equation}
where $(-)_!$ denotes the left Kan extension, cf. Appendix D.2.
\end{prop}
Let $k[\mathbb{T}]$ denote the $\mathcal{A}_{\infty}$-ring spectrum $k\otimes \Sigma^{\infty}_{+}\mathbb{T}$ over $k$ (more precisely, the Eilenberg-Maclane spectrum $Hk$). Let $k[\epsilon]:=k[\epsilon]/\epsilon^2$ denote the graded $k$-algebra with $|\epsilon|=-1$, which can be naturally viewed as an $\mathcal{A}_{\infty}$-ring spectrum over $k$. Then there is an equivalence of $\mathcal{A}_{\infty}$-$k$-ring spectra $\gamma:k[\epsilon]\simeq k[\mathbb{T}]$, which induces an equivalence of $\infty$-categories $\gamma^*:Mod_{k[\mathbb{T}]}\simeq Mod_{k[\epsilon]}$, where $Mod_{k[\mathbb{T}]}$ denotes the $\infty$-category of modules over the $k$-linear spectrum $k[\mathbb{T}]$. On the other hand, there is an equivalence of $\infty$-categories $Mod_{k[\mathbb{T}]}\simeq \mathrm{Fun}(B\mathbb{T}, Mod_k)$.  \par\indent
Let $CC: \mathrm{Fun}(\Lambda,\mathrm{Ch}_k)\rightarrow \mathrm{Ch}_{k[\epsilon]}$ denotes the ordinary functor that sends $X$ to its cyclic bar complex $CC(X)$, equipped with the unital $k[\epsilon]$-action where $\epsilon$ acts by Connes' differential $B$. $CC$ clearly preserves quasi-isomorphisms, and thus induces an $\infty$-functor $CC: \mathrm{Fun}(N\Lambda,Mod_k)\rightarrow Mod_{k[\epsilon]}$.
Proposition 3.11 is a consequence of the following result of \cite{Hoy}:
\begin{thm}
(Theorem 2.3 \cite{Hoy}) The following diagram
of $\infty$-categories commutes
\begin{equation}
\begin{tikzcd}[row sep=1.2cm, column sep=0.8cm]
\mathrm{Fun}(N\Lambda,Mod_k)\arrow[r,"{\int_{\mathbb{T}}}"]\arrow[dr,"CC"]& \mathrm{Fun}(B\mathbb{T},Mod_k)\simeq Mod_{k[\mathbb{T}]}\arrow[d,"{\gamma^*}"]\\
&Mod_{k[\epsilon]}
\end{tikzcd}.
\end{equation}\qed
\end{thm}

\noindent\emph{Proof of Proposition 3.11 given Theorem 3.12}.  We consider the following analogues of $CC$ in the non-unital settings. Define
\begin{equation}
CC^{nu}: \mathrm{Fun}(\overrightarrow{\Lambda},\mathrm{Ch}_k)\rightarrow \mathrm{Ch}_{k[\epsilon]}
\end{equation}
by
\begin{equation}
CC^{nu}(X):=CC(X)\oplus e^+CC(X),
\end{equation}
where the differential is given by the matrix (cf. Appendix B.5 for the definition of $b,b'$)
\begin{equation}
\begin{pmatrix}
b_X & \tau-1\\
0 & b'_X
\end{pmatrix}.
\end{equation}
Furthermore, this chain complex is equipped with the $k[\epsilon]$-action where $\epsilon$ acts by the matrix
\begin{equation}
\begin{pmatrix}
0 & 0\\
N & 0
\end{pmatrix}.
\end{equation}
This can be generalized to the non-unital $\mathcal{A}_{\infty}$-setting, where we define
\begin{equation}
CC^{nu}: \mathrm{Fun}^{dg}(\overrightarrow{\Lambda}\rtimes \mathcal{A}_{\infty}^{dg}, \mathrm{Mod}_k)\rightarrow \mathrm{Mod}_{k[\epsilon]}
\end{equation}
by
\begin{equation}
CC^{nu}(X):=CC(X)\oplus e^+CC(X),
\end{equation}
where the differential and the $k[\epsilon]$-action is given by the same formula as in (3.25) and (3.26), but using the $\mathcal{A}_{\infty}$-version of bar and cyclic bar differential; cf Definition 3.4 and 3.5.
When $X=A^{\sharp}$ is the Hochschild functor of a non-unital $\mathcal{A}_{\infty}$-algebra, $CC^{nu}(X)$ is the non-unital Hochschild chain complex, equipped with the non-unital Connes' differential. Given Theorem 3.12, to prove Proposition 3.11 it suffices to show there are natural equivalences:
\begin{equation}
CC^{nu}(X)\simeq CC((\iota_{\Lambda})_! (P^{-1}_{\overrightarrow{\Lambda}})^*X)
\end{equation}
of $k[\epsilon]$-modules, for $X\in \mathrm{Fun}^{dg}(\overrightarrow{\Lambda}\rtimes \mathcal{A}_{\infty}^{dg},\mathrm{Mod}_k)$.\par\indent
Consider the diagonal bimodule
\begin{equation}
k[\overrightarrow{\Lambda}\rtimes \mathcal{A}_{\infty}^{dg}]:\overrightarrow{\Lambda}^{op}\rtimes \mathcal{A}_{\infty}^{dg}\times \overrightarrow{\Lambda}\rtimes \mathcal{A}_{\infty}^{dg}\rightarrow \mathrm{Mod}_k
\end{equation}
given by
\begin{equation}
([m],[n])\mapsto \mathrm{Hom}_{\overrightarrow{\Lambda}\rtimes \mathcal{A}_{\infty}^{dg}}([m],[n]).
\end{equation}
Similarly, let
\begin{equation}
k[\Lambda]: \Lambda^{op}\times \Lambda\rightarrow \mathrm{Ch}_k
\end{equation}
denote the diagonal bimodule over $\Lambda$.
For $X\in \mathrm{Fun}^{dg}(\overrightarrow{\Lambda}\rtimes \mathcal{A}_{\infty}^{dg},\mathrm{Mod}_k)$, we have
\begin{equation}
CC^{nu}(X)\simeq CC^{nu}(k[\overrightarrow{\Lambda}\rtimes \mathcal{A}_{\infty}^{dg}])\otimes_{\overrightarrow{\Lambda}\rtimes \mathcal{A}_{\infty}^{dg}} X,
\end{equation}
where $\otimes_{\overrightarrow{\Lambda}\rtimes \mathcal{A}_{\infty}^{dg}}$ denotes coend pairing over $\overrightarrow{\Lambda}\rtimes \mathcal{A}_{\infty}^{dg}$. Similarly, for $X\in \mathrm{Fun}(\Lambda,\mathrm{Ch}_k)$, we have
\begin{equation}
CC(X)\simeq CC(k[\Lambda])\otimes_{\Lambda} X.
\end{equation}
There is a chain of weak equivalences
\begin{align}
CC(\iota_! (P^{-1})^*X)&\simeq CC(k[\Lambda])\otimes (\iota_{\Lambda})_!\circ (P^{-1}_{\overrightarrow{\Lambda}})^*(X)\\
&\simeq  (\iota_{\Lambda}\circ P_{\overrightarrow{\Lambda}})^*CC(k[\Lambda])\otimes X\\
&\simeq CC^{nu}(k[\overrightarrow{\Lambda}\rtimes \mathcal{A}_{\infty}^{dg}])\otimes X.
\end{align}
The first equivalence follows from (3.34). The second equivalence follows from the coend formula for Kan extension, cf Lemma D.6.
The third equivalence is induced by a composition
\begin{equation}
CC^{nu}\big(k[\overrightarrow{\Lambda}\rtimes \mathcal{A}_{\infty}^{dg}]\big)\rightarrow P_{\overrightarrow{\Lambda}}^*CC^{nu}(k[\overrightarrow{\Lambda}])\rightarrow (\iota_{\Lambda}\circ P_{\overrightarrow{\Lambda}})^*CC(k[\Lambda]).
\end{equation}
We now describe each map in the composition and show that they are weak equivalences. The first map in the composition
\begin{equation}
CC^{nu}\big(k[\overrightarrow{\Lambda}\rtimes \mathcal{A}_{\infty}^{dg}]\big)\rightarrow P_{\overrightarrow{\Lambda}}^*CC^{nu}(k[\overrightarrow{\Lambda}])
\end{equation}
is induced by the augmentation map on morphism spaces
\begin{equation}
\bigoplus_{f\in hom_{\overline{\Lambda}}([n],[m])} \bigotimes_{i=0}^m C^{cell}_{|f_i|}((\mathcal{A}_{\infty})_{|f^{-1}(i)|})\rightarrow \bigoplus_{f\in hom_{\overrightarrow{\Lambda}}([n],[m])} k.
\end{equation}
Since $(\mathcal{A}_{\infty})_d=\overline{\mathcal{R}}^{d+1}$ is contractible, (3.40) is a weak equivalence. To describe the second map,
\begin{equation}
P_{\overrightarrow{\Lambda}}^*CC^{nu}(k[\overrightarrow{\Lambda}])\rightarrow (\iota_{\Lambda}\circ P_{\overrightarrow{\Lambda}})^*CC(k[\Lambda]),
\end{equation}
we replace $CC(k[\Lambda])$ by a quasi-isomorphic cocyclic $k[\epsilon]$-module as follows.\par\indent
Fix $[m]$. The functors
\begin{equation}
\overrightarrow{\Lambda}([m],-)\;,\;\Lambda([m],-): \overrightarrow{\Delta}^{op}\rightarrow \mathrm{Sets}
\end{equation}
define semisimplicial sets, and $\Lambda([m],-)$ is in fact a genuine simplicial set. Let $\overline{CC}(\Lambda([m],-))$ denote the normalized cyclic bar complex obtained from $CC(\Lambda([m],-))$ by quotienting
out the degenerate subcomplex. The quotient map $CC(\Lambda([m],-))\rightarrow \overline{CC}(\Lambda([m],-))$ is a quasi-isomorphism, and $\overline{CC}(\Lambda([m],-))$ inherits a normalized
Connes' differential $\overline{B}$, giving it the structure of a $k[\epsilon]$-module. \par\indent
From the explicit description of $\Lambda$ in section 3.1, the elements of $\Lambda([m],[n])$ not in the images of the degeneracy maps $s_i: \Lambda([m],[n-1])\rightarrow \Lambda([m],[n])$
are given by (homotopy classes of) $f:S^1\rightarrow S^1$ that are surjective onto all except possibly the $0$-th marked point.
Thus, as a graded vector space, we can split $\overline{CC}(\Lambda([m],-))$ into two summands, the first is generated by surjective maps on marked points, and the second are maps that misses exactly the $0$-th marked point.
This gives a natural isomorphism
\begin{equation}
CC^{nu}(\overrightarrow{\Lambda}([m],-))=CC(\overrightarrow{\Lambda}([m],-))\otimes e^+CC(\overrightarrow{\Lambda}([m],-))\xrightarrow{\cong} \overline{CC}(\Lambda([m],-))
\end{equation}
which intertwines the differentials
\[
\begin{pmatrix}
b_{CC} & \tau-1\\
0 & b'_{CC}
\end{pmatrix}
\]
and $\overline{b}_{CC}$. Moreover, this isomorphism intertwines the $\epsilon$-actions given by
\[
\begin{pmatrix}
0&0\\
N&0
\end{pmatrix}
\]
with $\overline{B}$. This proves the third equivalence.\par\indent
Finally, using the equivalence (3.33) we conclude that $CC^{nu}(X)\simeq CC((\iota_{\Lambda})_! (P^{-1}_{\overrightarrow{\Lambda}})^*X)$. \qed\par\indent
We now discuss a consequence of Proposition 3.11.
\begin{cor}
Consider the diagram of $\infty$-categories
\begin{equation}
\begin{tikzcd}[row sep=1.2cm, column sep=0.8cm]
N\overrightarrow{\Delta}^{op}\arrow[d,"{\iota_{\Delta^{op}}}"]\arrow[r,"\overrightarrow{i}"]& N\overrightarrow{\Lambda}\arrow[d,"{\iota_{\Lambda}}"]\\
N\Delta^{op}\arrow[r,"i"]\arrow[d]&N\Lambda\arrow[d]\\
*\arrow[r,"\mathfrak{i}"]& B\mathbb{T}
\end{tikzcd}
\end{equation}
Suppose $X\in \mathrm{Fun}(\overrightarrow{\Lambda},\mathrm{Ch}_k)$ is $H$-unital, then the Beck-Chevalley transform (cf. Appendix D.3)
\begin{equation}
|(\iota_{\Delta^{op}})_!(\overrightarrow{i})^*X|\rightarrow \mathfrak{i}^*\int_{\mathbb{T}} (\iota_{\Lambda})_!X
\end{equation}
is a quasi-isomorphism. Here, $|\cdot|$ denotes the left Kan extension along $\Delta^{op}\rightarrow *$.
\end{cor}
\noindent\emph{Proof}. $|(\iota_{\Delta^{op}})_!(\overrightarrow{i})^*X|$ is computed by the cyclic bar complex $CC(X)$. On the other hand, by the argument of Proposition 3.11, there is a quasi-isomorphism
$\mathfrak{i}^*\int_{\mathbb{T}} (\iota_{\Lambda})_!X\simeq CC^{nu}(X)$. Under this identification, the natural map (3.45) is the inclusion $CC(X)\hookrightarrow CC^{nu}(X)$, which has cokernel the bar complex of
$X$. Since $X$ is $H$-unital, its bar complex is acyclic, and therefore (3.45) is a quasi-isomorphism.\qed
\begin{rmk}
By Lemma 3.10, the bottom square in (3.44) is homotopy exact. In contrast, neither the top square nor the outer square are homotopy exact squares.
Corollary 3.13 can be seen as a remedy for this failure when the $H$-unitality condition is satisfied.
\end{rmk}

\textbf{Negative finite cyclic homology as homotopy fixed points}. There is counterpart to the above story if one replaces $\Lambda$ by ${}_p\Lambda$, where $p$ is a prime.
\begin{mydef}
A dg functor $Q\in\mathrm{Fun}^{dg}(\overrightarrow{{}_p\Lambda}\rtimes \mathcal{A}_{\infty}^{dg}, \mathrm{Mod}_k)$ is called a \emph{non-unital} $\mathcal{A}_{\infty}$-\emph{finite $p$-cyclic $k$-module}.
\end{mydef}
We write $Q_{k_1,\cdots,k_p}$ for the value of $Q$ on the object $[k_1,\cdots,k_p]$.
\begin{mydef}
Let $Q$ be a non-unital $\mathcal{A}_{\infty}$-finite $p$-cyclic $k$-module. As a graded $k$-vector space, the $p$-\emph{fold cyclic bar complex} of $Q$ is
\begin{equation}
{}_pCC(Q):=\bigoplus_{k_1,\cdots,k_p\geq 0} Q_{k_1,\cdots,k_p}[k_1+\cdots+k_p].
\end{equation}
Its differential $b^p_Q$ is defined by, for an element $x\in Q_{k_1,\cdots,k_p}[k_1+\cdots+k_p]\subset {}_pCC(Q)$,
\begin{equation}
b^p_Q(x):=d_{Q_{\underline{k}}}(x)+\sum_{\underline{k'}<\underline{k}}\;\sum_{\substack{f\in \overrightarrow{{}_p\Lambda}(\underline{k},\underline{k'}):\exists!\;j\;s.t.\;|f^{-1}(j)|>1\\|f^{-1}(j)\cap\{p\;\textrm{distinguished pts}\}|\leq 1}}
\pm Q(\mathrm{id}\otimes \mathrm{id}\otimes \cdots \mu_{|f^{-1}(j)|}\otimes \cdots\otimes \mathrm{id})(x),
\end{equation}
where $\underline{k},\underline{k'}$ are the short hand notation for $[k_1,\cdots,k_p],[k_1',\cdots,k_p']$. The chain complex $({}_pCC(Q),b^p_Q)$ is called the $p$-\emph{fold cyclic bar complex} of $Q$. When $Q=i^*A^{\sharp}$ for an $\mathcal{A}_{\infty}$-algebra $A$, then
${}_pCC(Q)$ agrees with the $p$-fold Hochshild chain complex ${}_pCC(A)$ defined in section 2.3.
\end{mydef}
Let $\tau\in \overrightarrow{{}_p\Lambda}([k_1,k_2,\cdots,k_p],[k_p,k_1,\cdots,k_{p-1}])$ be the isomorphism that rotates `$p$ blocks' of marked points. By abuse of notation, we also denote by $\tau$
the morphism $(\tau, \mathrm{id}\otimes\cdots\otimes\mathrm{id})\in \overrightarrow{{}_p\Lambda}\rtimes \mathcal{A}_{\infty}^{dg}([k_1,k_2,\cdots,k_p],[k_p,k_1,\cdots,k_{p-1}])$. This induces a
$\mathbb{Z}/p$-action on ${}_pCC(Q)$.
\begin{mydef}
Define the \emph{negative} $\mathbb{Z}/p$-\emph{equivariant complex} of $Q$, denoted $CC^{\mathbb{Z}/p}(Q)$, to be the chain complex
\begin{equation}
CC^{\mathbb{Z}/p}(Q):={}_pCC(Q)[[t,\theta]],
\end{equation}
where $|t|=2,|\theta|=1, \theta^2=0$, equipped with a $t$-linear differential
\begin{equation}
\begin{cases}
x\mapsto b^p_Q(x)+(-1)^{|x|}(\tau-1)x\theta\\
x\theta\mapsto b^p_Q(x)\theta+(-1)^{|x|}(1+\tau+\cdots+\tau^{p-1})xt.
\end{cases}
\end{equation}
\end{mydef}

From the explicit description of ${}_p\Lambda$, every morphisms $f\in{}_p\Lambda([k_1,\cdots,k_p],[k_1',\cdots,k_p'])$ can be uniquely decomposed as
$f=s\circ f'$ where $f'\in(\Delta^{op})^p([k_1,\cdots,k_p],[k_1',\cdots,k_p'])$ and $s\in\mathbb{Z}/p\cong \mathrm{Aut}_{{}_p\Lambda}([k_1',\cdots,k_p'])$. In particular, since $N(\Delta^{op})^p$ is contractible,
the map $N{}_p\Lambda\rightarrow B\mathbb{Z}/p$ of simplicial sets, induced by the 1-categorical functor sending objects of ${}_p\Lambda$ to the unique object of $B\mathbb{Z}/p$ and morphisms
$f=s\circ f'\mapsto s$, is a weak equivalence. Since $B\mathbb{Z}/p$ is a Kan complex, the map $N{}_p\Lambda\rightarrow B\mathbb{Z}/p$ can be identified with the $\infty$-groupoid completion
$N{}_p\Lambda\rightarrow \tilde{N\Lambda}_p$ and $\mathrm{Aut}_{\tilde{N\Lambda}_p}([0,\cdots,0])\simeq \mathrm{Aut}_{N{}_p\Lambda}([0,\cdots,0])\simeq \mathbb{Z}/p$ as $\infty$-groups. In particular, we have a diagram
of $\infty$-categories
\begin{equation}
\begin{tikzcd}[row sep=1.2cm, column sep=0.8cm]
N(\Delta^{op})^p\arrow[d]\arrow[r,"{i_p}"]& N{}_p\Lambda\arrow[d]\\
*\arrow[r]& B\mathbb{Z}/p\simeq \tilde{N\Lambda}_p
\end{tikzcd}.
\end{equation}
\begin{lemma}
Diagram (3.50) is a homotopy exact square.
\end{lemma}
\noindent\emph{Proof}. By Theorem D.14, it suffices to show that for $[k_1,\cdots,k_p]\in {}_p\Lambda$ and $\varphi\in\mathbb{Z}/p$, viewed as an automorphism of the unique
object of $B\mathbb{Z}/p$, the double comma category $({}_{[k_1,\cdots,k_p]/}(\Delta^{op})^p_{/*})_{\varphi}$ is weakly contractible. \par\indent
The category $({}_{[k_1,\cdots,k_p]/}(\Delta^{op})^p_{/*})_{\varphi}$ has objects given by pairs $([k'_1,\cdots,k'_p], f)$, where $[k'_1,\cdots,k'_p]\in(\Delta^{op})^p$ and $f: [k_1,\cdots,k_p]\rightarrow [k'_1,\cdots,k_p']$
is a morphism in ${}_p\Lambda$ that decomposes uniquely as $f=\varphi\circ f'$ for some morphism $f'$ in $(\Delta^{op})^p$. Therefore, this category is equivalent to the overcategory
$(\Delta^{op})^p_{[k_1,\cdots,k_p]/}$, which is weakly contractible since it has an initial object. \qed\par\indent
Let $\int_{\mathbb{Z}/p}$ denote the left Kan extension along
$N{}_p\Lambda\rightarrow B\mathbb{Z}/p$. We can now prove an analogue of Proposition 3.11 in the finite cyclic setting.
\begin{prop}
Let $X$ be a non-unital $\mathcal{A}_{\infty}$-finite $p$-cyclic $k$-module. Then there is a natural quasi-isomorphism
\begin{equation}
CC^{\mathbb{Z}/p}(X)\simeq \big(\int_{\mathbb{Z}/p}(\iota_{{}_p\Lambda})_! (P^{-1}_{\overrightarrow{{}_p\Lambda}})^*X\big)^{h\mathbb{Z}/p}.
\end{equation}
\end{prop}
It suffices to prove the following theorem, and the deduction of Proposition 3.19 will follow in the same way that Proposition 3.11 follows from Theorem 3.12.
\begin{thm}
The following diagram of $\infty$-categories commutes
\begin{equation}
\begin{tikzcd}[row sep=1.2cm, column sep=0.8cm]
\mathrm{Fun}(N{}_p\Lambda,Mod_k)\arrow[r,"{\int_{\mathbb{Z}/p}}"]\arrow[dr,"{{}_pCC}"]& \mathrm{Fun}(B\mathbb{Z}/p, Mod_k)\arrow[d,"{\simeq}"]\\
&Mod_{k[\mathbb{Z}/p]}
\end{tikzcd},
\end{equation}
where ${}_pCC$ is $p$-fold cyclic bar construction in Definition 3.17, equipped with the $\mathbb{Z}/p$ action generated by $\tau$.
\end{thm}
\noindent\emph{Proof}. The idea of proof follows closely that of \cite[Theorem 2.3]{Hoy}. Let
\begin{equation}
k[{}_p\Lambda(-,-)]: {}_p\Lambda^{op}\times {}_p\Lambda\rightarrow \mathrm{Ch}_k
\end{equation}
be the diagonal bimodule. Fixing $[k_1,\cdots,k_p]\in{}_p\Lambda$, we have a dg $k[\mathbb{Z}/p]$-module ${}_pCC(k[{}_p\Lambda([k_1,\cdots,k_p],-)])$, which we simply denote ${}_pCC({}_p\Lambda([k_1,\cdots,k_p],-))$. As $[k_1,\cdots,k_p]$ varies over ${}_p\Lambda$,
this gives rise to a functor ${}_pCC({}_p\Lambda(-,-)): {}_p\Lambda^{op}\rightarrow \mathrm{Ch}_{k[\mathbb{Z}/p]}$. For $X\in\mathrm{Fun}({}_p\Lambda,\mathrm{Ch}_k)$, there is an isomorphism of dg-$k[\mathbb{Z}/p]$-module
\begin{equation}
{}_pCC(X)\xrightarrow{\cong}   {}_pCC({}_p\Lambda(-,-))\otimes_{{}_p\Lambda} X,
\end{equation}
where $\otimes$ denotes the coend pairing. Similarly, we have
\begin{equation}
\int_{\mathbb{Z}/p}X\simeq \int_{\mathbb{Z}/p}k[{}_p\Lambda(-,-)]\otimes_{{}_p\Lambda} X.
\end{equation}
Therefore, it suffices to prove the `universal' equivalence
\begin{equation}
\int_{\mathbb{Z}/p}k[{}_p\Lambda(-,-)]\simeq {}_pCC({}_p\Lambda(-,-))\in \mathrm{Fun}(N{}_p\Lambda^{op},Mod_{k[\mathbb{Z}/p]}).
\end{equation}
Consider the coalgebra $k\otimes_{k[\mathbb{Z}/p]}k$. Explicitly, it is given by the polynomial ring $k[\tilde{t},\tilde{\theta}], |\tilde{t}|=-2,|\tilde{\theta}|=-1,\tilde{\theta}^2=0$ with the
coalgebra structure
\begin{equation}
\begin{cases}
\tilde{t}^k\mapsto \sum_{k_1+k_2=k}\tilde{t}^{k_1}\otimes \tilde{t}^{k_2}\\
\tilde{t}^k\tilde{\theta}\mapsto \sum_{k_1+k_2=k}(\tilde{t}^{k_1}\tilde{\theta}\otimes \tilde{t}^{k_2}+\tilde{t}^{k_1}\otimes \tilde{t}^{k_2}\tilde{\theta})
\end{cases}.
\end{equation}
$-\otimes_{k[\mathbb{Z}/p]}k$ induces a fully faithful functor $Mod_{k[\mathbb{Z}/p]}\hookrightarrow Comod_{k[\tilde{t},\tilde{\theta}]}$. Therefore,
to prove (3.56), it suffices to show that
\begin{equation}
\int_{\mathbb{Z}/p}k[{}_p\Lambda(-,-)\otimes_{k[\mathbb{Z}/p]}k\simeq {}_pCC({}_p\Lambda(-,-))\otimes_{k[\mathbb{Z}/p]}k \in \mathrm{Fun}({}_p\Lambda^{op}, Comod_{k[\tilde{t},\tilde{\theta}]}).
\end{equation}
Note that both sides send morphisms in $N{}_p\Lambda^{op}$ to equivalences in $Comod_{k[\tilde{t},\tilde{\theta}]}$, and thus can be viewed as functors from $B\mathbb{Z}/p$ to $Comod_{k[\tilde{t},\tilde{\theta}]}$. Explicitly, a
$k[\tilde{t},\tilde{\theta}]$-comodule structure on $M\in Mod_k$ is classified by its $\tilde{t}$ and $\tilde{\theta}$ coefficients, i.e. a map $u: M\rightarrow M[2]$ and a map $\eta: M\rightarrow M[1]$. \par\indent
Let's first consider the left hand side of (3.58). By Lemma 3.18, given a functor $Y:{}_p\Lambda\rightarrow \mathrm{Sets}$, the underlying space of $\int_{\mathbb{Z}/p}Y$ is weakly equivalent
to $\mathrm{colim}_{N(\Delta^{op})^p}i_p^*Y=|i^*_pY|$, where $|\cdot|$ denote the $p$-fold geometric realization.\par\indent
Under the identification $\mathrm{Fun}(B\mathbb{Z}/p,Mod_k)\simeq Mod_{k[\mathbb{Z}/p]}$, we have
\begin{equation}
\int_{\mathbb{Z}/p}k[{}_p\Lambda(-,-)\otimes_{k[\mathbb{Z}/p]}k\simeq (\int_{\mathbb{Z}/p}k[{}_p\Lambda(-,-)])_{h\mathbb{Z}/p}\simeq k[(\int_{\mathbb{Z}/p}{}_p\Lambda(-,-))_{h\mathbb{Z}/p}].
\end{equation}
For each $\overrightarrow{k}\in{}_p\Lambda$, the geometric realization $|{}_p\Lambda(\overrightarrow{k},-)|$ is homeomorphic to $\mathbb{Z}/p\times \Delta^{\overrightarrow{k}}$, equipped with
the obvious $\mathbb{Z}/p$-action that cyclically permutes the first component and is trivial on the second component. Therefore, $(\int_{\mathbb{Z}/p}{}_p\Lambda(-,-))_{h\mathbb{Z}/p}$ is equivalent to the constant ${}_p\Lambda^{op}$-module $\underline{k}$ with (non-constant) $k[\tilde{t},\tilde{\theta}]$-comodule structure
classified by the standard generators of
\begin{equation}
[\underline{k},\underline{k}[2]]\simeq H^2(B\mathbb{Z}/p,k)=k\;\;,\;\;[\underline{k},\underline{k}[1]]\simeq H^1(B\mathbb{Z}/p,k)=k.
\end{equation}
Therefore, we need to show that the right hand side of (3.56) is also equivalent to $\underline{k}$ with the comodule structure classified by the standard generators of (3.60).\par\indent
The key observation is that
\begin{equation}
{}_pCC({}_p\Lambda(-,-))\otimes_{k[\mathbb{Z}/p]}k: {}_p\Lambda^{op}\rightarrow \mathrm{Ch}_k
\end{equation}
is a projective resolution of the constant ${}_p\Lambda^{op}$-module $\underline{k}$. This can be checked point-wise, so we fixed some $\overrightarrow{k}=[k_1,\cdots,k_p]\in {}_p\Lambda$. There is a bicomplex whose total complex computes $${}_pCC({}_p\Lambda([k_1,\cdots,k_p],-))\otimes_{k[\mathbb{Z}/p]}k.$$
Namely, consider
\begin{equation}
\begin{tikzcd}[row sep=1.2cm, column sep=1.5cm]
\vdots\arrow[d]&\vdots\arrow[d]&\vdots\arrow[d]&  \\
\bigoplus_{|\overrightarrow{k}'|=2} k[{}_p\Lambda(\overrightarrow{k},\overrightarrow{k}')]\arrow[d,"{b^p}"] &\bigoplus_{|\overrightarrow{k}'|=2} k[{}_p\Lambda(\overrightarrow{k},\overrightarrow{k}')]\arrow[d,"{b^p}"]\arrow[l,"{\tau-1}"]&\bigoplus_{|\overrightarrow{k}'|=2} k[{}_p\Lambda(\overrightarrow{k},\overrightarrow{k}')]\arrow[d,"{b^p}"]\arrow[l,"{1+\tau+\cdots+\tau^{p-1}}"]&\arrow[l] \cdots\\
\bigoplus_{|\overrightarrow{k}'|=1} k[{}_p\Lambda(\overrightarrow{k},\overrightarrow{k}')]\arrow[d,"{b^p}"] &\bigoplus_{|\overrightarrow{k}'|=2} k[{}_p\Lambda(\overrightarrow{k},\overrightarrow{k}')]\arrow[d,"{b^p}"]\arrow[l,"{\tau-1}"]&\bigoplus_{|\overrightarrow{k}'|=2} k[{}_p\Lambda(\overrightarrow{k},\overrightarrow{k}')]\arrow[d,"{b^p}"]\arrow[l,"{1+\tau+\cdots+\tau^{p-1}}"]&\arrow[l] \cdots\\
k[{}_p\Lambda(\overrightarrow{k},\overrightarrow{0})] &k[{}_p\Lambda(\overrightarrow{k},\overrightarrow{0})]\arrow[l,"{\tau-1}"]&k[{}_p\Lambda(\overrightarrow{k},\overrightarrow{0})]\arrow[l,"{1+\tau+\cdots+\tau^{p-1}}"]&\arrow[l] \cdots
\end{tikzcd}
\end{equation}
where we have written $|\overrightarrow{k}'|:=k_1'+\cdots+k_p'$. Projectivity is clear as each term in the bicomplex is a direct sum of representable functors. To see that this is a resolution
of $k$, we first look at each horizontal complex.\par\indent
There is a decomposition ${}_p\Lambda(\overrightarrow{k},\overrightarrow{k}')=(\Delta^{op})^p(\overrightarrow{k},\overrightarrow{k}')\times \mathbb{Z}/p$ compatible with $\mathbb{Z}/p$-action, where
the $\mathbb{Z}/p$-action on the right hand side is trivial on the first component, and the regular action on the second component. Therefore, the $i$-th horizontal complex is the tensor product of the
free $k$-module $\bigoplus_{|\overrightarrow{k}'|=i}k[(\Delta^{op})^p(\overrightarrow{k},\overrightarrow{k}')]$ with the complex
\begin{equation}
k[\mathbb{Z}/p]\xleftarrow{\tau-1}k[\mathbb{Z}/p]\xleftarrow{1+\tau+\cdots+\tau^{p-1}}k[\mathbb{Z}/p]\xleftarrow{\tau-1}\cdots,
\end{equation}
which is acyclic except in degree zero, where the homology is $k$. Therefore, the total complex of (3.62) is quasi-isomorphic to the vertical complex
\begin{equation}
\begin{tikzcd}[row sep=1.2cm, column sep=0.8cm]
\vdots\arrow[d]\\
\bigoplus_{|\overrightarrow{k}'|=2}(\Delta^{op})^p(\overrightarrow{k},\overrightarrow{k}')\arrow[d,"{b^p}"]\\
\bigoplus_{|\overrightarrow{k}'|=1}(\Delta^{op})^p(\overrightarrow{k},\overrightarrow{k}')\arrow[d,"{b^p}"]\\
(\Delta^{op})^p(\overrightarrow{k},\overrightarrow{0}).
\end{tikzcd}
\end{equation}
This complex computes the homology of $\Delta^{k_1}\times \Delta^{k_2}\times \cdots\times \Delta^{k_p}$ (which is the $p$-simplicial geometric realization of $(\Delta^{op})^p([k_1,\cdots,k_p],-)$), which is contractible. In particular, it is a
resolution of $k$.\par\indent
Writing the total complex of (3.62) as ${}_pCC({}_p\Lambda(-,-))[\tilde{t},\tilde{\theta}]$, then its $k[\tilde{t},\tilde{\theta}]$-comodule structure is classified by the degree $2$ map
\begin{equation}
u:
\begin{cases}
x\tilde{t}^k\mapsto x\tilde{t}^{k-1}\\
x\tilde{t}^k\tilde{\theta}\mapsto x\tilde{t}^{k-1}\tilde{\theta}
\end{cases}
\end{equation}
and the degree $1$ map
\begin{equation}
\eta:
\begin{cases}
x\tilde{t}^k\mapsto (-1)^{|x|}(\tau-1)^{p-2}x\tilde{t}^{k-1}\tilde{\theta}\\
x\tilde{t}^k\tilde{\theta}\mapsto x\tilde{t}^k
\end{cases}
\end{equation}
with the understanding that if the exponents become negative, we set the term zero. \par\indent
By the Yoneda lemma, $\mathrm{Hom}_{{}_p\Lambda^{op}}(\Lambda(-,[k_1,\cdots,k_p]),\underline{k})=\underline{k}_{k_1,\cdots,k_p}=k$. Thus, the mapping complex from ${}_pCC({}_p\Lambda(-,-))[\tilde{t},\tilde{\theta}]$ to $\underline{k}$ is quasi-isomorphic to $k[\tilde{t},\tilde{\theta}]$ (with trivial differential).
It is then easy to see that the compositions
\begin{equation}
{}_pCC({}_p\Lambda(-,-))[\tilde{t},\tilde{\theta}]\xrightarrow{u}{}_pCC({}_p\Lambda(-,-))[\tilde{t},\tilde{\theta}][2]\xrightarrow{\simeq}\underline{k}[2]
\end{equation}
and
\begin{equation}
{}_pCC({}_p\Lambda(-,-))[\tilde{t},\tilde{\theta}]\xrightarrow{\eta}{}_pCC({}_p\Lambda(-,-))[\tilde{t},\tilde{\theta}][1]\xrightarrow{\simeq}\underline{k}[1]
\end{equation}
correspond exactly to the standard generators $\tilde{t}\in[\underline{k},\underline{k}[2]]$ and $\tilde{\theta}\in[\underline{k},\underline{k}[1]]$.\qed\par\indent
3.4. \textbf{Cocyclic and finite cocyclic objects}. In this subsection, we briefly discuss the dual constructions associated to contravariant dg functors out of $\overrightarrow{\Lambda}\rtimes \mathcal{A}_{\infty}^{dg}$ and $\overrightarrow{{}_p\Lambda}\rtimes \mathcal{A}_{\infty}^{dg}$.
\begin{mydef}
A dg functor $Q: \overrightarrow{\Lambda}^{op}\rtimes \mathcal{A}^{dg}_{\infty}\rightarrow \mathrm{Mod}_k$ is called a \emph{non-unital} $\mathcal{A}_{\infty}$-\emph{cocyclic $k$-module}.
\end{mydef}
We denote $Q([n])$ as $Q_n$ and $d_{Q_n}$ its differential. For $d\geq 2$, let $\mu_d\in C_{d-2}^{cell}(\overline{R}^{d+1})$ be the top dimensional cell.
\begin{mydef}
Let $Q$ be a non-unital $A_{\infty}$-cocyclic $k$-module. As a graded $k$-vector space, the \emph{cocyclic cobar complex} of $Q$ is
\begin{equation}
CC^{\vee}(Q):=\prod_{n\geq 0} Q_n[-n],
\end{equation}
where $[-n]$ denotes shifting the degree by $-n$. We define the following operations on $CC^{\vee}(Q)$:
\begin{itemize}
    \item The \emph{cobar differential} $b'_Q$ is the degree $1$ differential on $CC^{\vee}(Q)$ defined by: for $x\in Q_n[-n]\subset CC^{\vee}(Q)$,
\begin{equation}
b'_Q(x):=d_{Q_n}(x)+\sum_{m>n}\sum_{\substack{f\in \overrightarrow{\Lambda}([m],[n]): \min f^{-1}(0)=0\\\textrm{and}\;\exists!\,j\in[n]\;\textrm{s.t.}\;|f^{-1}(j)|>1}}  \pm Q(\mathrm{id}\otimes\mathrm{id}\otimes\cdots\otimes \mu_{|f^{-1}(j)|}\otimes \cdots\otimes\mathrm{id})(x)
\end{equation}
and in general for $x=(x_0,x_1,x_2,\cdots)\in \prod_{n\geq 0} Q_n[-n]=CC^{\vee}(Q)$, we define $b'_Q(x)=b'_Q(x_0)+b'_Q(x_1)+b'_Q(x_2)+\cdots$. Note that for a fixed $m$, the $Q_m[-m]$-component of this expression
is a finite sum, and thus gives a well defined element of $CC^{\vee}(Q)$.
   \item The \emph{cocyclic cobar differential} $b_Q$ is the degree $1$ differential on $CC^{\vee}(Q)$
\begin{equation}
b_Q:=b'_Q+w,
\end{equation}
where $w$ is the \emph{wrapped around terms} defined by: for $x\in Q_n[-n]\subset CC^{\vee}(Q)$,
\begin{equation}
w(x):=\sum_{m>n}\sum_{\substack{f\in\overrightarrow{\Lambda}([m],[n]): 0\in f^{-1}(0)\backslash\{\min f^{-1}(0)\}\\\textrm{and}\;|f^{-1}(j)|=1\;\textrm{for all}\;j\neq 0}} \pm Q(\mu_{|f^{-1}(0)|}\otimes \mathrm{id}\otimes\cdots\otimes \mathrm{id})(x).
\end{equation}
$(CC(Q),b_Q)$ is called the \emph{cocyclic cobar complex} of $Q$.
\end{itemize}
The \emph{positive cocyclic complex} $CC^{S^1,\vee}(Q)$ is the totalization of
\begin{equation}
\cdots\xrightarrow{N}(CC^{\vee}(Q),b_Q)\xrightarrow{\tau-1} (CC^{\vee}(Q),b'_Q)\xrightarrow{N} (CC(Q)^{\vee},b_Q)\rightarrow 0,
\end{equation}
where $N$ acts as $1+\tau+\cdots+\tau^{d}$ on $Q_d[-d]\subset CC(Q)$. Alternatively, we can write the positive cocyclic complex as
\begin{equation}
CC^{\vee}(Q)[\tilde{u},e^+],
\end{equation}
where $|u|=-2, |e^+|=-1, (e^+)^2=0$, and the differential is given by
\begin{equation}
\begin{cases}
x\mapsto b_Q(x)\\
xu^k\mapsto b_Q(x)u^k+(-1)^{|x|}(\tau-1)xu^{k-1}e^+\,,\,k\geq 1\\
xu^ke^+\mapsto b'_Q(x)u^ke^++Nxu^k\,,\,k\geq 0.
\end{cases}
\end{equation}
\end{mydef}
\begin{mydef}
$Q\in \mathrm{Fun}^{dg}(\overrightarrow{\Lambda}^{op}\rtimes \mathcal{A}^{dg}_{\infty},\mathrm{Mod}_k)$ is called \emph{H-counital} if the cobar complex $(CC^{\vee}(Q),b'_Q)$ is acyclic.
\end{mydef}
Let $\int^{\mathbb{T}}$ denote the right Kan extension along $N\Lambda^{op}\rightarrow B\mathbb{T}$, the following proposition is the dual version of Proposition 3.11.
\begin{prop}
Let $X$ be a non-unital $\mathcal{A}_{\infty}$-cocyclic $k$-module. Then there is a natural equivalence
\begin{equation}
CC^{S^1,\vee}(X)\simeq \big(\int^{\mathbb{T}}(\iota^{op}_{\Lambda})_* ((P^{op}_{\overrightarrow{\Lambda}})^{-1})^*X\big)_{h\mathbb{T}}.
\end{equation}
\qed
\end{prop}
The following corollary is the dual to Corollary 3.13.
\begin{cor}
If $X$ is $H$-counital, the dual Beck-Chevalley transform
\begin{equation}
|(\iota^{op}_{\Delta})_*(\overrightarrow{i}^{op})^*X|\leftarrow (\mathfrak{i}^{op})^*\int^{\mathbb{T}}(\iota^{op}_{\Lambda})_*X
\end{equation}
associated to (the opposite of) diagram (3.44) is a quasi-isomorphism.
\end{cor}
We now consider finite $p$-cocyclic objects.
\begin{mydef}
A dg functor $Q\in\mathrm{Fun}^{dg}(\overrightarrow{{}_p\Lambda}^{op}\rtimes \mathcal{A}_{\infty}^{dg}, \mathrm{Mod}_k)$ is called a \emph{non-unital} $\mathcal{A}_{\infty}$-\emph{finite $p$-cocyclic $k$-module}.
\end{mydef}
We write $Q_{k_1,\cdots,k_p}$ for the value of $Q$ on the object $[k_1,\cdots,k_p]$.
\begin{mydef}
Let $Q$ be a non-unital $A_{\infty}$-finite $p$-cocyclic $k$-module. As a graded $k$-vector space, the $p$-\emph{fold cocyclic cobar complex} of $Q$ is
\begin{equation}
{}_pCC^{\vee}(Q):=\prod_{k_1,\cdots,k_p\geq 0} \pm Q_{k_1,\cdots,k_p}[-k_1-\cdots-k_p].
\end{equation}
Its differential $b^p_Q$ is defined by, for an element $x\in Q_{k_1,\cdots,k_p}[-k_1-\cdots-k_p]\subset {}_pCC^{\vee}(Q)$,
\begin{equation}
b^p_Q(x):=d_{Q_{\underline{k}}}(x)+\sum_{\underline{k'}>\underline{k}}\;\sum_{\substack{f\in \overrightarrow{{}_p\Lambda}(\underline{k'},\underline{k}):\exists!\;j\;s.t.\;|f^{-1}(j)|>1\\|f^{-1}(j)\cap\{p\;\textrm{distinguished pts}\}|\leq 1}}
\pm Q(\mathrm{id}\otimes \mathrm{id}\otimes \cdots \otimes\mu_{|f^{-1}(j)|}\otimes \cdots\otimes \mathrm{id})(x),
\end{equation}
where $\underline{k},\underline{k'}$ are the short hand notation for $[k_1,\cdots,k_p],[k_1',\cdots,k_p']$. For a general $x=(x_{k_1,\cdots,k_p})_{k_1,\cdots,k_p\geq 0}$, we define
$b^p_Q(x):=\sum_{k_1,\cdots,k_p\geq 0}b^p_Q(x_{k_1,\cdots,k_p})$, which is a finite sum in each $(k_1',\cdots,k_p')$-component, and thus gives a well defined element of ${}_pCC^{\vee}(Q)$.
The chain complex $({}_pCC^{\vee}(Q),b^p_Q)$ is called the $p$-\emph{fold cocyclic cobar complex} of $Q$.\par\indent
Let $\tau\in \overrightarrow{{}_p\Lambda}([k_1,k_2,\cdots,k_p],[k_p,k_1,\cdots,k_{p-1}])$ be the isomorphism that rotates `$p$ blocks' of marked points. This induces a
$\mathbb{Z}/p$-action on ${}_pCC^{\vee}(Q)$, and we define the \emph{positive} $\mathbb{Z}/p$-\emph{equivariant complex} of $Q$, denoted $CC^{\mathbb{Z}/p,\vee}(Q)$, as the chain complex
\begin{equation}
CC^{\mathbb{Z}/p}(Q):={}_pCC^{\vee}(Q)[\tilde{t},\tilde{\theta}],
\end{equation}
where $|\tilde{t}|=-2,|\theta|=-1, \tilde{\theta}^2=0$, equipped with the differential
\begin{equation}
\begin{cases}
x\mapsto b^p_Q(x)\\
x\tilde{t}^k\mapsto b^p_Q(x)\tilde{t}^k+(-1)^{|x|}(1+\tau+\cdot+\tau^{p-1})x\tilde{t}^{k-1}\tilde{\theta}\,,\,k\geq 1\\
x\tilde{t}^k\tilde{\theta}\mapsto b^p_Q(x)\tilde{t}^k\tilde{\theta}+(-1)^{|x|}(\tau-1)x\tilde{t}^k\,,\,k\geq 0.
\end{cases}
\end{equation}
\end{mydef}
Let $\int^{\mathbb{Z}/p}$ denote the right Kan extension along $N{}_p\Lambda\rightarrow B\mathbb{Z}/p$. The following proposition is the dual version of Proposition 3.17.
\begin{prop}
Let $X$ be a non-unital $\mathcal{A}_{\infty}$-finite $p$-cocyclic $k$-module. Then there is a natural quasi-isomorphism
\begin{equation}
CC^{\mathbb{Z}/p,\vee}(X)\simeq \big(\int^{\mathbb{Z}/p}(\iota^{op}_{{}_p\Lambda})_* ((P^{op}_{\overrightarrow{{}_p\Lambda}})^{-1})^*X\big)_{h\mathbb{Z}/p}.
\end{equation}
\qed
\end{prop}

\renewcommand{\theequation}{4.\arabic{equation}}
\setcounter{equation}{0}

\section{Operadic open-closed maps}
In classical Floer theory, one often considers a \emph{parametrized moduli problem} of counting solutions to a pseudo-holomorphic curve equation with varying domains, equipped with appropriate Floer data. For instance, in Lagrangian Floer theory,
the $\mathcal{A}_{\infty}$-structure map $\mu^d$ of the Fukaya category is defined out of a parametrized moduli problem where one allows the domain to vary over the top dimensional cell of $\overline{\mathcal{R}}^{d+1}$. Another example is the open-closed map from Hochschild homology of the Fukaya category to quantum cohomology, in which case the domain varies over the top dimension cell of $\overline{\mathcal{R}}^1_{d+1}$. Both of these constructions rely on particular cellular structures of the parameter space of domains. While having the advantage of being explicit,
these constructions have the drawback that they are not sufficiently functorial, as a lot of interesting maps among the space of domains might not be cellular with respect to an a priori chosen cell structure. Therefore,
it would be convenient to have a theory of parametrized moduli problems based on \emph{singular chains} instead of \emph{cellular chains} on the space of domains. \par\indent
The technical framework that addresses this issue is Abouzaid, Groman and Varolgunes' operadic Floer theory, \cite{AGV}, which uses a model for homology based on singular symmetric cubical chains. The goal of this section
is to introduce their methods, based on which we will define the \emph{operadic open-closed maps}.\par\indent
4.1. \textbf{Parameter spaces of domains}.
Let $X$ be a moduli space of disks with interior and boundary marked points. The primary examples for us are $X=\mathcal{R}^{d+1}$, $X=\mathcal{R}^1_{d+1}$ and $X=\mathcal{R}^1_{k_1,\cdots,k_p}$. All of the above $X$ have the property that:
\begin{itemize}
    \item $X$ has a Deligne-Mumford compactification $\overline{X}$ consisting of nodal disks, such that set theoretically $\overline{X}=\bigsqcup_{T\in \mathcal{T}} X^T$, where $\mathcal{T}$ (which depends on $X$) is a collection of decorated finite planar trees and each $X^T$ is a product of parameter spaces of unbroken disks.
    \item The topology of $X$ near each boundary strata is compatible with the gluing of disks near the nodes.
\end{itemize}
We now discuss the three examples, and in particular, specify the collection of tree types $\mathcal{T}$ in each case.
\begin{enumerate}[label=\arabic*)]
    \item The spaces $\overline{X}=\overline{\mathcal{R}}^{d+1}$ was discussed in section 2.1. In this case, $\mathcal{T}$ consists of equivalence classes of rooted planar trees with $d\geq 2$ leaves, such that each vertex is adjacent to $\geq 3$ edges. The boundary component of $\overline{\mathcal{R}}^{d+1}$
corresponding to $T$ is
\begin{equation}
\mathcal{R}^T:=\prod_{v\in T} \mathcal{R}^{|v|},
\end{equation}
where $|v|$ is the valency of $v$ (cf. (2.7)).
   \item The spaces $\overline{X}=\overline{\mathcal{R}}^1_d$ are special cases of the spaces ${}_k\overline{\check{\mathcal{R}}}^1_d$ defined in section 2.2. In this case, $\mathcal{T}$ consists of equivalence classes of rooted planar trees with $d\geq 0$ leaves, with one internal vertex marked as \emph{main}, such that
each non-main vertex is adjacent to $\geq 3$ edges. The boundary component of $\overline{\mathcal{R}}^1_{d+1}$ corresponding to $T$ is
\begin{equation}
\mathcal{R}^T:=\mathcal{R}^1_{|main|}\times \prod_{v\in T\backslash\{main\}} \mathcal{R}^{|v|}.
\end{equation}
\begin{figure}[H]
 \centering
 \includegraphics[width=1.0\textwidth]{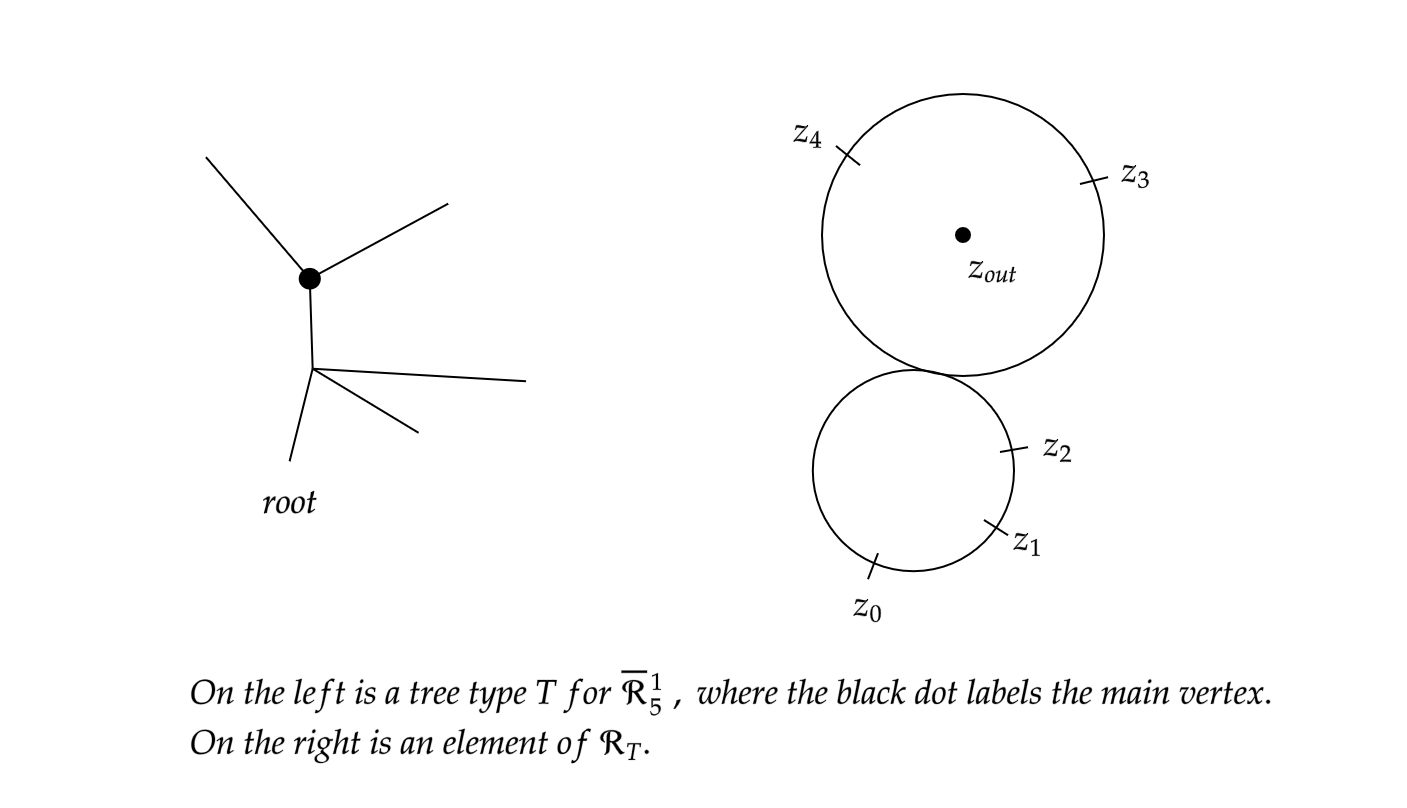}
 \caption{}
\end{figure}
 \item The spaces $\overline{X}=\overline{\mathcal{R}}^1_{k_1,\cdots,k_p}$ are defined in section 2.3. In this case, $\mathcal{T}$ consists of equivalence classes of rooted planar trees with $k_1+\cdots+k_p+p-1$ leaves, with one internal vertex
marked as \emph{main} and $p-1$ leaves marked as \emph{distinguished}, see Figure 5. Those $p-1$ leaves, together with the root, are called the $p$ distinguished semi-infinite edges. Moreover, we require that no two distinguished semi-infinite edges
are adjacent to the same non-main vertex and each non-main vertex is adjacent to $\geq 3$ edges. In particular, one has $|main|\geq p$. The boundary component of $\overline{\mathcal{R}}^1_{k_1,\cdots,k_p}$
corresponding to $T$ is
\begin{equation}
\mathcal{R}_T:=\mathcal{R}^1_{k_1^T,\cdots,k_p^T}\times \prod_{v\in \mathcal{T}\backslash\{main\}} \mathcal{R}^{|v|},
\end{equation}
where the $k_i^T$'s are defined as follows: among the edges adjacent to the main vertex, there are $p$ special ones that are contained in the respective paths from the main vertex to the $p$ distinguished semi-infinite edges. $k_i^T$ is defined as the number of edges in between the $i$-th and $i+1$-th special edge in counter-clockwise order.
\end{enumerate}
4.2. \textbf{Operadic Floer theory}. Take $X$ to be one of the spaces of domains from section 4.1. We now describe the symmetric cubical set $\mathcal{F}_{\bullet}(\overline{X})$ of Floer data on $\overline{X}$, following \cite{AGV}. An element of the $0$-cubes $\mathcal{F}_0(\overline{X})$ consists of the following data:
\begin{enumerate}[label=\arabic*)]
    \item A stable disk $\Sigma\in \overline{X}$, and a labeling of each interior and boundary marked points as input/output.
    \item For each component $\Sigma_v$ of $\Sigma$, and each boundary marked point $p$ of $\Sigma_v$, a choice of strip-like ends at $p$
\begin{equation}
\epsilon_p^{+}:[0,\infty)\times [0,1]\rightarrow \Sigma_v\;\;\mathrm{or}\;\;\epsilon_p^-:(-\infty,0]\times[0,1]\rightarrow \Sigma_v
\end{equation}
depending on whether $p$ is an output or input.
  \item For each component $\Sigma_v$, a pair $(K_v,J_v)$ where $K_v\in \Omega^1(\Sigma_v, \mathcal{H}), J_v\in C^{\infty}(\Sigma_v,\mathcal{J})$ satisfying
\begin{equation}
(\epsilon^{\pm}_p)^*K_v=H_tdt\;\;,\;\;(\epsilon^{\pm}_p)^*J_v=J_t,
\end{equation}
where $H_t, J_t$ are the pre-chosen time-dependent Hamiltonians and almost complex structures,
\end{enumerate}
modulo the following equivalence relation: two such data are equivalent if there exists an isomorphism of the underlying stable disks which intertwines the strip-like ends, and intertwines the
perturbation data $(K,J)$ up to rescaling of a constant, cf. Definition 2.3 \cite{AGV}. We let $\mathcal{F}_{0,T}(\overline{X})$ be the subset of $\mathcal{F}_0(\overline{X})$ where the underlying stable disk has combinatorial type $T$.
\begin{figure}[H]
 \centering
 \includegraphics[width=1.0\textwidth]{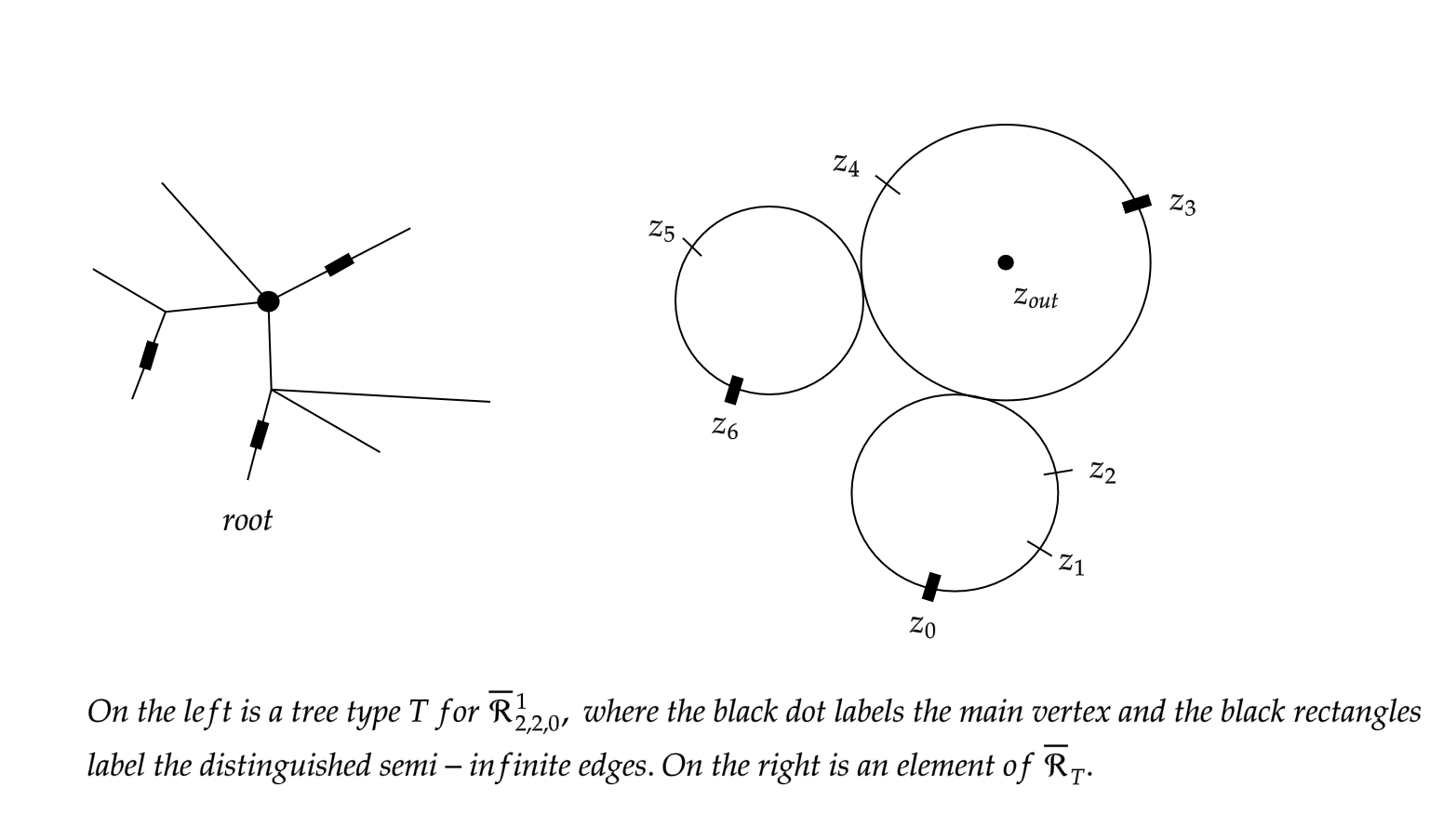}
 \caption{}
\end{figure}
\begin{rmk}
For the three spaces discussed in section 4.1, the convention for input/output is the following. For $\overline{X}=\overline{\mathcal{R}}^{d+1}$, the root is the output, and all other marked points are inputs. For $\overline{X}=\overline{\mathcal{R}}^1_{d+1}$ and $\overline{X}=\overline{\mathcal{R}}^1_{k_1,\cdots,k_p}$, the unique interior marked point is the output, and all boundary marked points are inputs.
\end{rmk}
More generally, an element $\mathfrak{o}\in\mathcal{F}_n(\overline{X})$ consists of the following data:
\begin{enumerate}[label=\arabic*)]
    \item For each face $f$ of the standard topological $n$-cube $\Box^n$, a choice of a tree $T_f\in \mathcal{T}$.
    \item A map of sets $b: \Box^n\rightarrow \mathcal{F}_0(\overline{X})$, such that the stable disks underlying $b(f)$ has combinatorial type $T_f$.
    \item For each face $f$ of $\Box^n$, a parametrization of pre-glued disk with Floer data by a corner chart
\begin{equation}
b_f: f\times [0,\epsilon)^{Codim(f)}\rightarrow \mathcal{F}_{0,T_f}(\overline{X}),
\end{equation}
and a choice of gluing lengths parametrized by the same chart
\begin{equation}
g_f: f\times [0,\epsilon)^{Codim(f)}\rightarrow [0,1)^{|E_{int}(T_f)|}.
\end{equation}
The triple $(T_f,b_f,g_f)$ is called a \emph{local model for the face} $f$, and the collection of triples $(T_f,b_f,g_f)_{f\in\mathrm{face}(\Box^n)}$ is called a \emph{gluing atlas}.
\end{enumerate}
These are required to satisfy:
\begin{itemize}
    \item For each $f$ and $p\in f\times [0,\epsilon)^{E_{int}(T_f)}$,
\begin{equation}
b(p)=\Gamma_{g_f(p)}(b_f(p)),
\end{equation}
where $\Gamma$ stands for gluing of the nodal Riemann surface equipped with Floer data given by $b_f(p)$, with gluing lengths $g_f(p)$, cf. \cite[Definition A.4]{AGV},\cite[section (9e)]{Sei1}.
     \item For an inclusion of faces $f_0\subset f_1$, let $E$ be the set of edges of $T_{f_0}$ that get collapsed under $T_{f_0}\rightarrow T_{f_1}$. Let $g_{f_0}^E: f_0\times [0,\epsilon)^{Codim(f_0)}\rightarrow [0,1)^E$ be the composition of $g_{f_0}$ with $[0,1)^{E_{int}(T_{f_0})}\rightarrow [0,1)^E$. Then we have
\begin{equation}
b_{f_1}=\Gamma_{g_{f_0}^E}(b_{f_0}).
\end{equation}
\end{itemize}
The face maps (cf. Definition C.1) $d^{\pm}_{n,i}$ are defined on the underlying map $b$ by $b\mapsto b\circ \iota^{pm}_{n,i}$, the degeneracy maps $\sigma_i$ are given by $b\mapsto b\circ \pi_i$ and the transposition maps are defined
by $b\mapsto b\circ \tau_{n,i}$; these induce obvious maps on the gluing atlases. Thus, we have defined the symmetric cubical set $\mathcal{F}_{\bullet}(\overrightarrow{X})$. The following proposition
guarantees that this is a correct model for homotopy theory on $\overline{X}$.
\begin{prop}
The natural projection
\begin{equation}
\pi: \mathcal{F}_{\bullet}(\overline{X})\rightarrow \overline{X}^{\Box}_{\bullet}
\end{equation}
is a homotopy equivalence of cubical sets, where $(-)^{\Box}_{\bullet}$ denotes the singular cubical set of a space.
\end{prop}
\noindent\emph{Proof}. This is an immediate adaptation of Proposition 2.20, Lemma 2.21 and Lemma 2.22 of \cite{AGV}.\qed\par\indent
\textbf{Transversality}. Having defined $\mathcal{F}_{\bullet}(\overline{X})$, we recall \cite{AGV}'s construction of a `Floer functor'. Informally, for each cube $\mathfrak{o}$ in $\mathcal{F}_{\bullet}(\overline{X})$, we would like to count rigid solutions to a parametrized moduli problem specified by $\mathfrak{o}$, and that would give rise to certain Floer operations. Moreover, to ensure that these operations satisfy the desired algebraic properties, we need to consider the parametrized moduli space of virtual dimension $1$ and its boundary.\par\indent
Hence, we need to equip the parametrized moduli space associated to each $\mathfrak{o}$ the structure of a topological manifold with boundary (at least in virtual dimension $0$ and $1$). To achieve this, an extra piece of data is added to the construction $\mathcal{F}_{\bullet}$, namely a perturbation of Floer data in the thick part of a glued surface. This is used to achieve transversality for the moduli spaces, because in general, the result of gluing disks equipped with regular Floer data might not be regular. We define a new cubical set (cf \cite[Def 3.8]{AGV}) $\tilde{\mathcal{F}}_{\bullet}(\overline{X})$ where an $n$-cube in $\tilde{\mathcal{F}}_{n}(\overline{X})$ is given by a triple $(B,(T_f,b_f,g_f)_{f\in\mathrm{face}(\Box^n)})$, where $B$ is a map from $[0,1]\times [0,1]^n$ to $\mathcal{F}_0(\overline{X})$ and $(T_f,b_f,g_f)$ is a local model for $f$. Writing $b=B|_{\{0\}\times [0,1]^n}$, these data are required to satisfy:
\begin{enumerate}[label=\arabic*)]
    \item The projection of $B$ onto $\overline{X}$ is independent of the first coordinate.
    \item $(b,(T_f,b_f,g_f)_f)$ defines an element of $\mathcal{F}_n(\overline{X})$.
    \item $B$ is supported in the thick part of each domain, with uniformly bounded $C^k$ norms for all $k$.
\end{enumerate}
\begin{lemma}
The map
\begin{equation}
\pi: \tilde{\mathcal{F}}_{\bullet}(\overline{X})\rightarrow \mathcal{F}_{\bullet}(\overline{X})
\end{equation}
forgetting the datum of $B$ is a homotopy equivalence.
\end{lemma}
\noindent\emph{Proof}. The proof is the same as \cite[Lemma 3.11]{AGV}.\qed\par\indent
For a cube $\mathfrak{o}\in \tilde{\mathcal{F}}_{n}(\overline{X})$, let $\overline{\mathcal{M}}(\mathfrak{o})$ be the parametrized moduli space of solutions to the Floer equation associated to $B|_{\{1\}\times [0,1]^n}$. For a point $u\in\overline{\mathcal{M}}(\mathfrak{o})$, if the underlying disk of $u$ is smooth, we define the \emph{extended Cauchy-Riemann operator} at $u$ to be the sum of the linearized Cauchy-Riemann operator with the deformation of the operator corresponding to moving the domain (with Floer data) within the cube $\mathfrak{o}$. In general, we take it to be the product of the extended Cauchy-Riemann operators of each component.
\begin{mydef}
A cube $\mathfrak{o}\in \tilde{\mathcal{F}}_n(\overline{X})$ is \emph{regular} if the extended linearized Cauchy-Riemann operator is surjective for all elements in $\overline{\mathcal{M}}(\mathfrak{o})$.
\end{mydef}
\noindent Since regularity is closed under taking boundaries and degeneracies, we can define:
\begin{mydef}
The cubical set of regular Floer data $\tilde{\mathcal{F}}^{reg}_{\bullet}(\overline{X})$ is the cubical subset of $\tilde{\mathcal{F}}_{\bullet}(\overline{X})$ consisting of regular cubes.
\end{mydef}
\begin{prop}
The inclusion $\tilde{\mathcal{F}}_{\bullet}^{reg}(\overline{X})\subset \tilde{\mathcal{F}}_{\bullet}(\overline{X})$ is a homotopy equivalence of cubical sets.
\end{prop}
\noindent \emph{Proof}. This follows from the proof of \cite[Prop 3.4]{AGV}. In fact, there is a deformation retraction $r: \tilde{\mathcal{F}}_{\bullet}(\overline{X})\rightarrow \tilde{\mathcal{F}}^{reg}_{\bullet}(\overline{X})$. \qed\par\indent
\textbf{The operadic cubical and dg associahedron}. Consider the topological associahedron $\mathcal{A}_{\infty}$ from section 3; for $d\geq 2$, we have $(\mathcal{A}_{\infty})_d=\overline{\mathcal{R}}^{d+1}$, whose
strata are described in section 4.1 1). There is an induced operadic structure on the collection of symmetric cubical sets $\{\tilde{\mathcal{F}}^{reg}((\mathcal{A}_{\infty})_d)\}_{d\geq 1}$, given by:\par\indent
Given a pair $(\mathfrak{o}_1,\mathfrak{o}_2)\in \tilde{\mathcal{F}}^{reg}_{n_1}((\mathcal{A}_{\infty})_{d_1})\times\tilde{\mathcal{F}}^{reg}_{n_2}((\mathcal{A}_{\infty})_{d_2}) $, we define $\mathfrak{o}_2\circ_i\mathfrak{o}_1$ to be the product $n_1+n_2$-cube which is pointwise given by sending
$(x,y)$ to $\mathfrak{o}_2\circ_i\mathfrak{o}_2$, where $\circ_i$ is the operadic structure map on $\mathcal{A}_{\infty}$, i.e. concatenation of disks (equipped with Floer data). This defines a map of symmetric cubical sets (cf. Appendix C for the definition of the symmetric monoidal product)
\begin{equation}
\circ_i: \tilde{\mathcal{F}}^{reg}_{\bullet}((\mathcal{A}_{\infty})_{d_1})\otimes\tilde{\mathcal{F}}^{reg}_{\bullet}((\mathcal{A}_{\infty})_{d_2})\rightarrow \tilde{\mathcal{F}}^{reg}_{\bullet}((\mathcal{A}_{\infty})_{d_1+d_2-1}),
\end{equation}
and makes the collection $\{\tilde{\mathcal{F}}^{reg}((\mathcal{A}_{\infty})_d)\}_{d\geq 1}$ into an operad valued in symmetric cubical sets. We call this the \emph{operadic cubical associahedron}, denoted $\mathcal{A}_{\infty}^{oper,\Box}$.\par\indent
Now apply the symmetric monoidal functor $C_*$ of normalized cubical chains (cf Appendix C), we obtain a dg operad $\{C_*((\mathcal{A}_{\infty}^{oper,\Box})_d)\}_{d\geq 2}$. We call this the
\emph{operadic dg associahedron}, denoted $\mathcal{A}_{\infty}^{oper,dg}$.
\begin{lemma}
Let $L$ be a monotone Lagrangian brane, then there is a map of dg operads
\begin{equation}
\mathcal{A}_{\infty}^{oper,dg}\rightarrow End(CF^*(L,L)).
\end{equation}
\end{lemma}
\noindent\emph{Proof}. The proof is analogous to \cite[Proposition 3.4]{AGV}, so we only describe how the map goes.\par\indent
For an $n$-cube $\mathfrak{o}\in (\mathcal{A}^{oper,dg}_{\infty})_d=C_*(\tilde{\mathcal{F}}^{reg}(\overline{\mathcal{R}}^{d+1}))$, we assign the multilinear map $CF(L,L)^{\otimes d}\rightarrow CF(L,L)$ of degree $-n$ obtained by counting rigid elements
\begin{equation}
(x\in [0,1]^n, u: \Sigma_{(\mathfrak{o}(x))}\rightarrow M)
\end{equation}
that satisfy Floer's equation
\begin{equation}
(du-{Y_K}_{\mathfrak{o}(x)})^{0,1}_{J_{\mathfrak{o}(x)}}=0,
\end{equation}
whose boundary lies on $L$, and with asymptotics given by elements of $CF^*(L,L)$. Since $\mathfrak{o}$ is regular cube, this parametrized moduli space of rigid elements is a $0$-manifold. To see that this assignment defines a map of dg operads, we consider the $1$-dimensional parametrized moduli space. By transversality and gluing, the boundary of this moduli spaces is the union of the parametrized moduli space
associated to $\partial \mathfrak{o}$, together with the moduli space involving semi-stable strip breaking (i.e. the differential of $CF^*(L,L)$). \qed\par\indent
4.3. \textbf{Operadic open-closed maps}. In this subsection, we define the \emph{operadic open-closed map} and its equivariant versions. \par\indent
For simplicity, we pretend the category $\mathrm{Fuk}(M)_{\lambda}$ contains a single Lagrangian $L$. We will explain in Appendix A the generalization to multiple objects. We denote $A_L:=CF^*(L,L)$. Let $\overrightarrow{\Lambda}\rtimes\mathcal{A}_{\infty}^{oper,dg}$ be the dg category described in Definition 3.1. By Lemma 4.7, $A_L$ is an algebra over the dg operad $\mathcal{A}_{\infty}^{oper,dg}$. We can form the \emph{operadic Hochschild functor},
which is the dg functor
\begin{equation}
A^{\sharp}_L: \overrightarrow{\Lambda}\rtimes \mathcal{A}_{\infty}^{oper,dg}\rightarrow \mathrm{Mod}_k
\end{equation}
sending $[n]\mapsto (A_L,\mu^1)^{\otimes n+1}$, just as in Definition 3.2.
Let $QH:=CM^*(f)$ be the Morse chain complex serving as a chain model for quantum cohomology (cf. section 2.2). Postcomposing with
the mapping chain complex functor $\mathrm{Map}(-,QH)$ with $A^{\sharp}_L$ one obtains a dg functor
\begin{equation}
\mathrm{Map}(A^{\sharp}_L,QH): \overrightarrow{\Lambda}^{op}\rtimes \mathcal{A}_{\infty}^{oper,dg}\rightarrow \mathrm{Mod}_k.
\end{equation}
That is, on objects, (4.17) sends $[n]$ to the chain complex $\mathrm{Map}(A_L^{\otimes n+1}, QH)$.
We call a dg functor from $\overrightarrow{\Lambda}\rtimes\mathcal{A}_{\infty}^{oper,dg}$ to $\mathrm{Mod}_k$ a \emph{non-unital $\mathcal{A}_{\infty}^{oper,dg}$-cocyclic $k$-module}. \par\indent
\textbf{The (non-unital) $\mathcal{A}_{\infty}^{oper,dg}$-cocyclic structure on $C_*(\tilde{\mathcal{F}}^{reg}(\overline{\mathcal{R}}^1_d)$}. For each $d\geq 1$, there is an action of $\frac{\mathbb{Z}}{d+1}$ on $\overline{\mathcal{R}}^1_{d+1}$
by cyclically permuting the $d+1$ boundary marked points. When the underlying Riemann surfaces are equipped with Floer data, one can package these actions into the structure of a non-unital $\mathcal{A}_{\infty}^{oper,dg}$-cocyclic $k$-module.
\begin{mydef}
We define $\overline{\mathfrak{R}}^1\in\mathrm{Fun}^{dg}(\overrightarrow{\Lambda}^{op}\rtimes \mathcal{A}_{\infty}^{oper,dg},\mathrm{Mod}_k)$ as follows. On objects, $\overline{\mathfrak{R}}^1$ is given by
\begin{equation}
[d]\mapsto C_*(\tilde{\mathcal{F}}^{reg}(\overline{\mathcal{R}}^1_{d+1})).
\end{equation}
Since every morphism $\tilde{h}\in\overrightarrow{\Lambda}\rtimes \mathcal{A}_{\infty}^{oper,dg}([m],[n])$ can be uniquely decomposed as $\tilde{h}=s\circ \tilde{f}$, where
$\tilde{f}\in \overrightarrow{\Delta}^{op}\rtimes\mathcal{A}_{\infty}^{oper,dg}([m],[n])$ and $s\in\frac{\mathbb{Z}}{n+1}$, it suffices to describe $\overline{\mathfrak{R}}^1(\tilde{f})$ and $\overline{\mathfrak{R}}^1(s)$ separately.
\begin{itemize}
    \item Let $f\in \Delta^{op}([m],[n])$. Recall from section 3.1 that viewing $f$ as a homotopy class of maps on $S^1$ (with marked points), we have $f(0)=0$. As $f^{-1}(0)=\{i_1<\cdots<i_{|f^{-1}(0)|}\}$ comes with a canonical ordering, we can define $l_0$ to be the unique index such that $i_{l_0}=0$. Now, given $f$ and stable disks $S_{main}\in \overline{\mathcal{R}}^1_{n+1}, S_i\in \overline{\mathcal{R}}^{|f^{-1}(i)|+1}, 0\leq i\leq n$, we can define the pre-glued stable disk
\begin{equation}
S_{main}\cup_f\{S_0,\cdots,S_n\}
\end{equation}
by concatenating the output of $S_i$ at the $i$-th marked point of $S_{main}$, and the $l_0$-th boundary marked point of $S_0$ becomes the $0$-th boundary marked point of $S_{main}\cup_f\{S_0,\cdots,S_n\}$.
By section 4.1 2), this
defines an element of $\overline{\mathcal{R}}^1_{m+1}$. If $S_{main}$ and $S_i,0\leq i\leq n$ are equipped with Floer data (on each component), then $S_{main}\cup_f\{S_1,\cdots,S_n\}$ is naturally equipped with a pre-glued Floer data.
In particular, given $f\in\overrightarrow{\Delta}^{op}([m],[n])$ and cubes $c_i\in \tilde{\mathcal{F}}^{reg}_{l_i}(\overline{R}^{|f^{-1}(i)|+1})$ so that
\begin{equation}
\tilde{f}=\Big(f\in\overrightarrow{\Delta}^{op}([m],[n]), \otimes_{i=0}^n c_i\in \bigotimes_{i=0}^n C_{l_i}(\tilde{\mathcal{F}}^{reg}(\overline{\mathcal{R}}_{|f^{-1}(i)|+1}))\Big)
\end{equation}
represents a morphism from $[m]$ to $[n]$ in  $\overrightarrow{\Delta}^{op}\rtimes \mathcal{A}_{\infty}^{oper,dg}$, together with a cube $c\in C_{l}(\tilde{\mathcal{F}}^{reg}(\overline{R}^1_{n+1}))$, we define $\overline{\mathfrak{R}}^1(\tilde{f})(c)$
to be the product cube in $C_{l+l_0+\cdots+l_n}(\tilde{\mathcal{F}}^{reg}(\overline{R}^1_{m+1}))$ that assigns a point $(x,x_0,\cdots,x_n)$ to $c(x)\cup_f\{c_0(x_0),\cdots,c_n(x_n)\}$ (while taking the
product gluing atlas and perturbation). It is clear that regularity is preserved under this process.
   \item For $s\in \frac{\mathbb{Z}}{n+1}$, viewed as an automorphism of $[n]$, we define
\begin{equation}
\overline{\mathfrak{R}}^1(s): C_*(\tilde{\mathcal{F}}^{reg}(\overline{\mathcal{R}}^1_{n+1}))\rightarrow C_*(\tilde{\mathcal{F}}^{reg}(\overline{\mathcal{R}}^1_{n+1}))
\end{equation}
by cyclically permuting the boundary marked points of the underlying surfaces in $\overline{\mathcal{R}}^1_{n+1}$, equipped with the pulled back Floer data.
\end{itemize}
\end{mydef}
Let $\mathfrak{o}$ be a cube in $C_l(\tilde{\mathcal{F}}^{reg}(\overline{\mathcal{R}}^1_{d+1}))$, $\mathbf{a}=a_0\otimes a_1\otimes a_n\in A_L^{\otimes n+1}$ and $y_{out}\in\mathrm{crit(f)}$. We define $\overline{M}(\mathfrak{o}, \mathbf{a},y_{out})$ to be the moduli space of pairs
\begin{equation}
\overline{M}(\mathfrak{o}, \mathbf{a},y_{out}):=\{(x\in[0,1]^l, u:\Sigma_{\mathfrak{o}(x)}\rightarrow M)\}
\end{equation}
that satisfy Floer's equation $(du-{Y_K}_{\mathfrak{o}(x)})^{0,1}_{J_{\mathfrak{o}(x)}}=0$, has boundary compoents constrained on $L$, has boundary marked points asymptotic to $a_0,\cdots,a_n$ and interior
marked point incident to the unstable manifold of $y_{out}$. By a standard boundary analysis for parametrized moduli space, the assignment
\begin{equation}
C_*(\tilde{\mathcal{F}}^{reg}(\overline{\mathcal{R}}^1_{d+1}))\rightarrow Map(A^{\otimes d+1}_L,QH)
\end{equation}
sending a cube $\mathfrak{o}\in C_l(\tilde{\mathcal{F}}^{reg}(\overline{\mathcal{R}}^1_{d+1}))$ to the map $A^{\otimes d+1}_L\rightarrow QH$ defined by counting rigid elements of (4.22)
defines a map of non-unital $\mathcal{A}_{\infty}^{oper,dg}$-cocyclic $k$-module of degree $\frac{\dim_{\mathbb{R}} M}{2}$, denoted as
\begin{equation}
OC^{oper}:\overline{\mathfrak{R}}^1\rightarrow\mathrm{Map}(A^{\sharp}_L,QH).
\end{equation}
\begin{mydef}
We call the map (4.24) the \emph{operadic open-closed map}.
\end{mydef}
We now define the operadic negative cyclic open-closed map.\par\indent
There is a projection $\tilde{P}^{op}_{\overrightarrow{\Lambda}}: \overrightarrow{\Lambda}^{op}\rtimes \mathcal{A}_{\infty}^{oper,dg}\rightarrow \overrightarrow{\Lambda}^{op}_k$ induced
by the augmentation map $C_*(\tilde{\mathcal{F}}^{reg}(\overline{\mathcal{R}}^{d+1}))\rightarrow k$. This is a unital quasi-equivalence of dg categories, and we let $(\tilde{P}^{op}_{\overrightarrow{\Lambda}})^{-1}$ denote a choice (unique up to homotopy) of $\mathcal{A}_{\infty}$-quasi-inverse.
Given a dg functor $X\in \mathrm{Fun}^{dg}(\overrightarrow{\Lambda}^{op}\rtimes \mathcal{A}_{\infty}^{oper,dg}, \mathrm{Mod}_k)$, similar to (3.18) we can obtain an $\infty$-functor $(\tilde{P}^{op}_{\overrightarrow{\Lambda}})^{-1})^*X\in \mathrm{Fun}(N\overrightarrow{\Lambda}^{op},Mod_k)_0$. We then right Kan extend along $\iota^{op}_{\Lambda}$ to obtain an $\infty$-functor
\begin{equation}
(\iota^{op}_{\Lambda})_*(\tilde{P}^{op}_{\overrightarrow{\Lambda}})^{-1})^*X\in \mathrm{Fun}(N\Lambda^{op},Mod_k)_0.
\end{equation}
Then, apply $\int^{\mathbb{T}}$ to (4.25), we obtain
\begin{equation}
\int^{\mathbb{T}}(\iota^{op}_{\Lambda})_*(\tilde{P}^{op}_{\overrightarrow{\Lambda}})^{-1})^*X\in \mathrm{Fun}(B\mathbb{T},Mod_k)_0.
\end{equation}
Since $\int^{\mathbb{T}}(\iota^{op}_{\Lambda})_*(\tilde{P}^{op}_{\overrightarrow{\Lambda}})^{-1})^*$ is clearly functorial in $X$, we can apply it to (4.24) and obtain a map of $\mathbb{T}$-modules
\begin{equation}
\int^{\mathbb{T}}(\iota^{op}_{\Lambda})_* ((\tilde{P}^{op}_{\overrightarrow{\Lambda}})^{-1})^*OC^{oper}: \int^{\mathbb{T}}(\iota^{op}_{\Lambda})_*((\tilde{P}^{op}_{\overrightarrow{\Lambda}})^{-1})^*\overline{\mathfrak{R}}^1\rightarrow \int^{\mathbb{T}}(\iota^{op}_{\Lambda})_*((\tilde{P}^{op}_{\overrightarrow{\Lambda}})^{-1})^*\mathrm{Map}(A^{\sharp}_L,QH).
\end{equation}
For any functors $f: \mathcal{C}\rightarrow \mathcal{D}, Y: \mathcal{C}\rightarrow N^{dg}\mathrm{Mod}_k$ and chain complex $B\in \mathrm{Mod}_k$, there is a canonical map $\mathrm{Map}(f_!Y,B)\rightarrow f^{op}_*\mathrm{Map}(Y,B)$ in $\mathrm{Fun}(\mathcal{D},N^{dg}\mathrm{Mod}_k)$ by the universal property of Kan extensions. Thus, there is a canonical map of $\mathbb{T}$-modules
\begin{equation}
\int^{\mathbb{T}}(\iota^{op}_{\Lambda})_* ((\tilde{P}^{op}_{\overrightarrow{\Lambda}})^{-1})^*\mathrm{Map}(A^{\sharp}_L,QH)\xleftarrow{\simeq} \mathrm{Map}(\int_{\mathbb{T}}(\iota_{\overrightarrow{\Lambda}})_! (\tilde{P}_{\overrightarrow{\Lambda}}^{-1})^*A^{\sharp}_L,QH).
\end{equation}
We claim that (4.28) is an equivalence, which can be checked by pulling back along $\mathfrak{i}: *\rightarrow B\mathbb{T}$ (an equivalence of $\mathbb{T}$-modules can be checked on the underlying chain complexes). Since $A_L^{\sharp}$ is $H$-unital, by Corollary 3.13 the pullback of (4.28) along $\mathfrak{i}$ is just the natural map $\mathrm{lim}_{\overrightarrow{\Lambda}^{op}}\mathrm{Map}((\tilde{P}_{\overrightarrow{\Lambda}}^{-1})^*A_L^{\sharp},QH)\leftarrow\mathrm{Map}(\mathrm{colim}_{\overrightarrow{\Lambda}}((\tilde{P}_{\overrightarrow{\Lambda}}^{-1})^*A_L^{\sharp},QH)$, which is an equivalence by the universal property of (co)limits. \par\indent
Postcomposing (4.27) with the inverse of (4.28), one obtains a map of $\mathbb{T}$-modules
\begin{equation}
\int^{\mathbb{T}}(\iota^{op}_{\Lambda})_* ((\tilde{P}^{op}_{\overrightarrow{\Lambda}})^{-1})^*\overline{\mathfrak{R}}^1\rightarrow \mathrm{Map}(\int_{\mathbb{T}}(\iota_{\overrightarrow{\Lambda}})_! (\tilde{P}_{\overrightarrow{\Lambda}}^{-1})^*A^{\sharp}_L,QH).
\end{equation}
Now, apply $(-)_{h\mathbb{T}}$ to (4.29). Since $(-)^{h\mathbb{T}}$ ($(-)_{h\mathbb{T}}$) is nothing but the (co)limit over $B\mathbb{T}$, for any $\mathbb{T}$-module $A$ and chain complex $B$ over $k$, there is a natural map (not an equivalence in general) $\mathrm{Map}(A,B)_{h\mathbb{T}}\rightarrow \mathrm{Map}(A^{h\mathbb{T}},B)$ of chain complexes. Thus, we obtain a map of chain complexes
\begin{equation}
\Big(\int^{\mathbb{T}}(\iota^{op}_{\Lambda})_* ((\tilde{P}^{op}_{\overrightarrow{\Lambda}})^{-1})^*\overline{\mathfrak{R}}^1\Big)_{h\mathbb{T}}\rightarrow \mathrm{Map}(\Big(\int_{\mathbb{T}}(\iota_{\overrightarrow{\Lambda}})_! (\tilde{P}_{\overrightarrow{\Lambda}})^{-1})^*A^{\sharp}_L\Big)^{h\mathbb{T}},QH).
\end{equation}
Let $\underline{k}$ denote the constant non-unital $\mathcal{A}_{\infty}^{oper,dg}$-cocyclic $k$-module with value $k$. There is a map $\overline{\mathfrak{R}}^1\rightarrow \underline{k}$ of non-unital $\mathcal{A}_{\infty}^{oper,dg}$-cocyclic $k$-modules induced by the augmentation maps
$C_*(\tilde{\mathcal{F}}^{reg}(\overline{\mathcal{R}}^1_{d+1})\rightarrow k$. This map is an equivalence since $\overline{\mathcal{R}}^1_{d+1}$ is contractible. Therefore, by Proposition 3.24,
\begin{equation}
H^*(\Big(\int^{\mathbb{T}}(\iota^{op}_{\Lambda})_* ((\tilde{P}^{op}_{\overrightarrow{\Lambda}})^{-1})^*\overline{\mathfrak{R}}^1\Big)_{h\mathbb{T}})\cong H^*(CC^{S^1,\vee}(\underline{k}))\cong k[\tilde{t}], |\tilde{t}|=-2.
\end{equation}
Thus, for each $m$, the image of $\tilde{t}^m\in k[\tilde{t}]=H^*(\Big(\int^{\mathbb{T}}(\iota^{op}_{\Lambda})_* ((\tilde{P}^{op}_{\overrightarrow{\Lambda}})^{-1})^*\overline{\mathfrak{R}}^1\Big)_{h\mathbb{T}})$ under (4.30) defines a chain map of degree $-2m+\frac{1}{2}\dim_{\mathbb{R}} M$
\begin{equation}
OC^{S^1,oper}_m:CC^{S^1,oper}(A_L):=\Big(\int_{\mathbb{T}}(\iota_{\overrightarrow{\Lambda}})_! ((\tilde{P}_{\overrightarrow{\Lambda}})^{-1})^*A^{\sharp}_L\Big)^{h\mathbb{T}}\rightarrow QH(M).
\end{equation}
\begin{mydef}
The \emph{operadic negative cyclic open-closed map} is the degree $\frac{1}{2}\dim_{\mathbb{R}} M$ chain map
\begin{equation}
OC^{S^1,oper}:=\sum_{m\geq 0} OC^{S^1,oper}_mu^m: CC^{S^1,oper}(A)\rightarrow QH(M)[[u]], |u|=2.
\end{equation}
\end{mydef}
One can define the operadic $\mathbb{Z}/p$-equivariant open-closed map analogously. First, we pullback the operadic open-closed map along $\overrightarrow{j^{op}_{\mathcal{A}_{\infty}^{oper,dg}}}: \overrightarrow{\Lambda}^{op}_p\rtimes\mathcal{A}_{\infty}^{oper,dg}\rightarrow\overrightarrow{\Lambda}^{op}\rtimes \mathcal{A}_{\infty}^{oper,dg}$ to obtain
a map
\begin{equation}
(\overrightarrow{j^{op}_{\mathcal{A}_{\infty}^{oper,dg}}})^*OC^{oper}:(\overrightarrow{j^{op}_{\mathcal{A}_{\infty}^{oper,dg}}})^*\overline{\mathfrak{R}}^1\rightarrow (\overrightarrow{j^{op}_{\mathcal{A}_{\infty}^{oper,dg}}})^*\mathrm{Map}(A^{\sharp}_L,QH)=\mathrm{Map}((\overrightarrow{j_{\mathcal{A}_{\infty}^{oper,dg}}})^*A^{\sharp}_L,QH).
\end{equation}
As in the $\mathbb{T}$-equivariant case, we first apply $\int^{\mathbb{Z}/p}(\iota^{op}_{{}_p\Lambda})_* ((\tilde{P}^{op}_{\overrightarrow{{}_p\Lambda}})^{-1})^*(-)$ and then $(-)_{h\mathbb{Z}/p}$ to obtain a map
\begin{equation}
\Big(\int^{\mathbb{Z}/p}(\iota^{op}_{{}_p\Lambda})_*((\tilde{P}^{op}_{\overrightarrow{{}_p\Lambda}})^{-1})^*(\overrightarrow{j^{op}_{\mathcal{A}_{\infty}^{oper,dg}}})^*\overline{\mathfrak{R}}^1\Big)_{h\mathbb{Z}/p}\rightarrow \mathrm{Map}(\Big(\int_{\mathbb{Z}/p}(\iota_{\overrightarrow{{}_p\Lambda}})_!(\tilde{P}_{\overrightarrow{{}_p\Lambda}}^{-1})^*(\overrightarrow{j_{\mathcal{A}_{\infty}^{oper,dg}}})^*A^{\sharp}_L\Big)^{h\mathbb{Z}/p},QH).
\end{equation}
Since $(\overrightarrow{j^{op}_{\mathcal{A}_{\infty}^{oper,dg}}})^*\overline{\mathfrak{R}}^1\rightarrow \underline{k}$ is an equivalence of dg-$\overrightarrow{{}_p\Lambda}\rtimes \mathcal{A}_{\infty}^{oper,dg}$-modules,
we have, by Proposition 3.28,
\begin{equation}
H^*(\Big(\int^{\mathbb{Z}/p}(\iota^{op}_{{}_p\Lambda})_* ((\tilde{P}^{op}_{\overrightarrow{{}_p\Lambda}})^{-1})^*(\overrightarrow{j^{op}_{\mathcal{A}_{\infty}^{oper,dg}}}^*)^*\overline{\mathfrak{R}}^1\Big)_{h\mathbb{Z}/p})\cong H^*(CC^{\mathbb{Z}/p,\vee}(\underline{k}))\cong k[\tilde{t},\tilde{\theta}],
\end{equation}
where $|\tilde{t}|=-2, |\tilde{\theta}|=-1, \tilde{\theta}^2=0$.
Therefore, the images of $\tilde{t}^m$ (resp. $\tilde{t}^m\tilde{\theta}$) under (4.35) define chain maps
\begin{equation}
OC^{\mathbb{Z}/p,oper}_{2m}\,(\textrm{resp}.\;\;OC^{\mathbb{Z}/p,oper}_{2m+1}): CC^{\mathbb{Z}/p,oper}(A):=\Big(\int_{\mathbb{Z}/p}(\iota_{\overrightarrow{{}_p\Lambda}})_!(\tilde{P}_{\overrightarrow{{}_p\Lambda}}^{-1})^*(\overrightarrow{j_{\mathcal{A}_{\infty}^{oper,dg}}})^*A^{\sharp}_L\Big)^{h\mathbb{Z}/p}\rightarrow QH
\end{equation}
of degree $-2m+\frac{1}{2}\dim_{\mathbb{R}} M$ (resp. $-2m-1+\frac{1}{2}\dim_{\mathbb{R}} M$).
\begin{mydef}
The \emph{operadic negative $\mathbb{Z}/p$-equivariant open-closed map} is the degree $\frac{1}{2}\dim_{\mathbb{R}} M$ chain map
\begin{equation}
OC^{\mathbb{Z}/p,oper}:=\sum_{m\geq 0} (OC^{\mathbb{Z}/p,oper}_{2m}+OC^{\mathbb{Z}/p,oper}_{2m+1}\theta)t^k: CC^{\mathbb{Z}/p,oper}(A_L)\rightarrow QH(M)[[t,\theta]],
\end{equation}
where $|t|=2, |\theta|=1,\theta^2=0$.
\end{mydef}
To conclude this section, we note the following property of $\overline{\mathfrak{R}}^1$.
\begin{lemma}
$\overline{\mathfrak{R}}^1$ is $H$-counital.
\end{lemma}
\emph{Proof}. By definition, we need to show the cobar complex $(CC^{\vee}(\overline{\mathfrak{R}}^1),b')$ is acyclic. Since $\overline{\mathcal{R}}^1_{d+1}$ is contractible, there is a quasi-equivalence of non-unital $\mathcal{A}_{\infty}$-cocyclic
$k$-modules $\overline{\mathfrak{R}}^1\rightarrow \underline{k}$. Here, $\underline{k}$ denotes the non-unital $\mathcal{A}_{\infty}$-cocyclic
$k$-module with constant value $k$, and the structure map $k\rightarrow k$ associated to a morphism
$g\in \bigotimes_{i=0}^m (\mathcal{A}_{\infty})_{|f^{-1}(i)|}, f\in \overrightarrow{\Lambda}([n],[m])$ is given by
multiplication by the image of $g$ under the (tensor product of) augmentation map
\begin{equation}
\bigotimes_{i=0}^m (\mathcal{A}_{\infty})_{|f^{-1}(i)|}\rightarrow \bigotimes_{i=0}^m k\cong k.
\end{equation}
As a result, there is a quasi-isomorphism $(CC^{\vee}(\overline{\mathfrak{R}}^1),b')\simeq (CC^{\vee}(\underline{k}),b')$. Since $(CC^{\vee}(\underline{k}),b')$ is dual to the bar complex $(CC(\underline{k}),b')$, it suffices to show that the latter is acyclic. But $(CC(\underline{k}),b')$ can be explicitly written as (the $2$-periodification of)
\begin{equation}
0\rightarrow k\xrightarrow{1} k\xrightarrow{0} k\xrightarrow{1}\cdots,
\end{equation}
which is clearly acyclic.
\qed

\renewcommand{\theequation}{5.\arabic{equation}}
\setcounter{equation}{0}

\section{$\mathbb{Z}/p$-Gysin comparison map}
Given a topological space $X$ with an $S^1$-action, restriction along $p$-th roots of unity gives an induced $\mathbb{Z}/p$-action on $X$. In classical topology, the $\mathbb{Z}/p$-Gysin sequence states that there is an isomorphism $H^*_{\mathbb{Z}/p}(X;\mathbb{F}_p)=H^*_{S^1}(X;\mathbb{F}_p)\oplus H^{*-1}_{S^1}(X;\mathbb{F}_p)$ on equivariant homologies. This isomorphism easily follows from the Gysin long exact sequence associated to the $S^1$-bundle $X\times_{\mathbb{Z}/p}ES^1\rightarrow X\times_{S^1}ES^1$. In section 5.1, we formulate and prove a version of the $\mathbb{Z}/p$-Gysin comparison for chain complexes over $k$, where chain level $S^1$-actions are modeled using Connes' cyclic category $\Lambda$ and the induced $\mathbb{Z}/p$-action is modeled by considering the restriction along $j: {}_p\Lambda\rightarrow \Lambda$ (cf. section 3.1).
In section 5.2, we prove that under the $\mathbb{Z}/p$-Gysin comparison (for the operadic Hochschild complex), the operadic negative cyclic and finite cyclic open-closed maps are compatible. \par\indent
5.1. \textbf{$\mathbb{Z}/p$-Gysin comparison for cyclic modules}. As in section 3, here we also work in the context of $\infty$-categories. Taking $\infty$-groupoid completion of the functor $j: N{}_p\Lambda\rightarrow N\Lambda$, we get a commutative square
of $\infty$-categories
\begin{equation}
\begin{tikzcd}[row sep=1.2cm, column sep=0.8cm]
N{}_p\Lambda\arrow[r,"j"]\arrow[d]& N\Lambda\arrow[d]\\
\tilde{N\Lambda}_p\arrow[r,"{\tilde{j}}"]& \tilde{N\Lambda}
\end{tikzcd}
\end{equation}
In section 3, we saw that $\tilde{N\Lambda}_p\simeq B\mathbb{Z}/p$ and $\tilde{N\Lambda}\simeq B\mathbb{T}$ via identifying $\mathbb{Z}/p$ (resp. $\mathbb{T}$) with the automorphism $\infty$-group
of an object of $\tilde{N\Lambda}_p$ (resp. $\tilde{N\Lambda}$).
\begin{lemma}
Under the above identifications, $\tilde{j}$ is homotopic to the map induced by the standard inclusion $\mathfrak{j}:\mathbb{Z}/p\subset \mathbb{T}$.
\end{lemma}
\noindent\emph{Proof}. Consider the object $[0,\cdots,0]\in{}_p\Lambda$, then $\tilde{\Lambda}_p([0,\cdots,0],-)\simeq\int_{\mathbb{Z}/p}{}_p\Lambda([0,\cdots,0],-)$, where we recall that $\int_{\mathbb{Z}/p}$ denotes the left Kan extension along $N{}_p\Lambda\rightarrow B\mathbb{Z}/p$. By definition, the underlying space (forgetting the $\mathbb{Z}/p$-action) of $\tilde{\Lambda}_p([0,\cdots,0],-)$ is weakly equivalent to $\mathrm{Aut}_{\tilde{N\Lambda}}([0,\cdots,0])$, and by Lemma 3.18 the underlying space of $\int_{\mathbb{Z}/p}{}_p\Lambda([0,\cdots,0],-)$ is weakly equivalent to
the geometric realization (denoted $|\cdot|$) of the underlying $p$-fold simplicial set of ${}_p\Lambda([0,\cdots,0],-)$. The only nondegenerate simplices of ${}_p\Lambda([0,\cdots,0],-)$ are its $0$-simplicies, and thus
its geometric realization is the discrete space $\mathbb{Z}/p$. The image of $[0,\cdots,0]\in{}_p\Lambda$ under $j$ is $[p-1]\in \Lambda$, and the underlying space of $\tilde{\Lambda}([p-1],-)\simeq \int_{\mathbb{Z}/p}\tilde{\Lambda}([p-1],-)$ is equivalent to
the geometric realization of the underlying $p$-fold simplicial set of $\Lambda([p-1],-)$ (via $(\Delta^{op})^p\rightarrow {}_p\Lambda\rightarrow \Lambda$). Hence, it suffices to show that (the geometric realizations of)
the map of $p$-fold simplicial sets ${}_p\Lambda([0,\cdots,0],-)\rightarrow j^*\Lambda(j([0,\cdots,0],-)=j^*\Lambda([p-1],-)$ is homotopic to the inclusion $\mathbb{Z}/p\subset S^1$. \par\indent
Consider the cyclic set $C$ given by $C_n=\frac{\mathbb{Z}}{n+1}=\mathrm{Aut}_{\Lambda}([n])$ (cf. \cite[6.1.10]{Lod} for a full description), whose geometric realization is $S^1$. For each $[m],[n]\in\Lambda$, recall that a morphism $f\in \Lambda([m],[n])$ canonically factors as $f=g(f)\circ s(f), g(f)\in\Delta^{op}([m],[n]), s(f)\in C_n$. Therefore, there is a map of cyclic sets
$\Lambda([m],-)\rightarrow C$ that sends $f\in\Lambda([m],[n])$ to $s(f)\in C_n$. The geometric realization of $\Lambda([m],-)$ is homeomorphic to $S^1\times \Delta^m$ (cf. \cite[E.7.2.1]{Lod}), and $\Lambda([m],-)\rightarrow C$ induces the projection $S^1\times \Delta^m\rightarrow S^1$ on geometric realizations. Since the geometric realization of $C$ is homotopy equivalent to the ($p$-fold) geometric realization of $j^*C$ (cf. Lemma B.10), it suffices to show that
the composition ${}_p\Lambda([0,\cdots,0],-)\rightarrow j^*\Lambda([p-1],-)\rightarrow j^*C$ induces the inclusion $\mathbb{Z}/p\subset S^1$ on ($p$-fold) geometric realizations. \par\indent
The $(k_1,\cdots,k_p)$-simplices of $j^*C$ are given by the set $\mathbb{Z}/(k_1+\cdots+k_p+p)$. But upon a more careful observation, one sees that $j^*C$ has only
$2p$ non-degenerate simplices:
\begin{itemize}
    \item All of its $(0,0,\cdots,0)$-simplices, of which there are $p$.
    \item For each $1\leq i\leq p$, the $(0,0,\cdots,1,\cdots,0)$-simplex (where $1$ is in position $i$) corresponding to
$i\in \{0,1,\cdots,p\}\cong \mathbb{Z}/(p+1)=(j^*C)_{0,0,\cdots,1,\cdots,0}$. One can think of this as corresponding to the unique (up to homotopy) circle with $p+1$ marked points, where the $i$-th marked point is non-distinguished,
and all the other $p$ marked points are distinguished.
\end{itemize}
Therefore,
\begin{align*}
|j^*C|&=\Big(\bigsqcup_{k_1,\cdots,k_p}\frac{\mathbb{Z}}{k_1+\cdots+k_p+p}\times\Delta^{k_1}\times\cdots\times \Delta^{k_p}\Big)\,/\sim\\
&=\Big(\mathbb{Z}/p\times \Delta^0\times\cdots\times\Delta^0\sqcup\bigsqcup_{i=1}^p \Delta^0\times\cdots\Delta^1\times\cdots\times \Delta^0\Big)\,/\sim,
\end{align*}
where the equivalence relation in the last equation is exactly identifying the endpoints of the $p$ intervals cyclically with the $p$ $(0,0,\cdots,0)$-simplices, indicated in Figure 6.
Thus, the geometric realization of $j^*C$ is a copy of $S^1$.
\begin{figure}[H]
 \centering
 \includegraphics[width=1.0\textwidth]{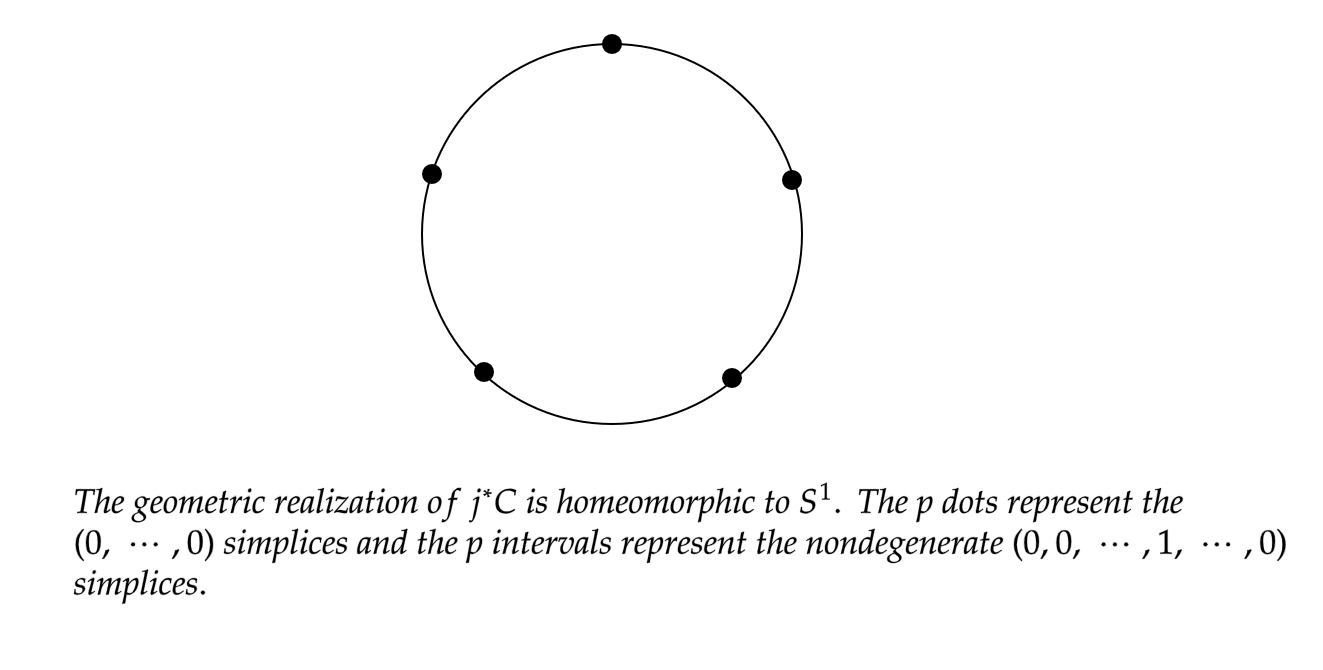}
 \caption{}
\end{figure}
Since the map of $p$-fold simplicial sets ${}_p\Lambda([0,\cdots,0],-)\rightarrow j^*C$ gives a bijection on $(0,\cdots,0)$-simplices, the induced map on geometric realization is the inclusion
$\mathbb{Z}/p\subset S^1$. \qed
\begin{lemma}
All the squares in
\begin{equation}
\begin{tikzcd}[row sep=1.2cm, column sep=0.8cm]
N(\Delta^{op})^p\arrow[r,"{i_p}"]\arrow[d]& N{}_p\Lambda\arrow[r,"j"]\arrow[d]& N\Lambda\arrow[d]\\
*\arrow[r,"\mathfrak{i}_p"]& B\mathbb{Z}/p\arrow[r,"\mathfrak{j}"]& B\mathbb{T}
\end{tikzcd}
\end{equation}
are homotopy exact.
\end{lemma}
\emph{Proof}. The left square of (5.2) is a homotopy exact square by Lemma 3.18. Therefore, it suffices to show that the outer square is homotopy exact. The outer square of (5.2) has another decomposition given by
\begin{equation}
\begin{tikzcd}[row sep=1.2cm, column sep=0.8cm]
N(\Delta^{op})^p\arrow[r,"o"]\arrow[d]& N\Delta^{op}\arrow[r,"i"]\arrow[d]& N\Lambda\arrow[d]\\
*\arrow[r,"\mathrm{id}"]& *\arrow[r,"\mathfrak{i}"]& B\mathbb{T}.
\end{tikzcd}
\end{equation}
The right square of (5.3) is homotopy exact by \cite[Proposition 1.1]{Hoy}. Therefore it suffices to show that the left square of (5.3) is homotopy exact, or equivalently that the functor
$o: N(\Delta^{op})^p\rightarrow N\Delta^{op}$ is cofinal. We precompose this with the diagonal functor $\mathrm{diag}:N\Delta^{op}\rightarrow N(\Delta^{op})^p$, which is cofinal by \cite[Lemma 5.5.8.4]{Lur1}.
It thus suffices to show that the composition $o\circ \mathrm{diag}: N\Delta^{op}\rightarrow N\Delta^{op}$ is cofinal. This follows from \cite[Proposition 2.1]{Bar} and \cite[Theorem 4.6]{Vel}.\qed\par\indent
We now prove the main result of this subsection.
\begin{prop}($\mathbb{Z}/p$-Gysin comparison)
Let $k$ be a field of characteristic $p$. Let $X: N\Lambda\rightarrow Mod_k$ be a cyclic $k$-module and $X': N\Lambda^{op}\rightarrow Mod_k$ a cocyclic $k$-module.
Then there are quasi-isomorphisms
\begin{equation}
(\int_{\mathbb{Z}/p}j^*X)^{h\mathbb{Z}/p}\xleftarrow{\sim} (\int_{\mathbb{T}}X)^{h\mathbb{T}}\langle 1,\theta\rangle:\phi_p
\end{equation}
and
\begin{equation}
\tilde{\phi}_p: (\int^{\mathbb{Z}/p}(j^{op})^*X')_{h\mathbb{Z}/p}\xrightarrow{\sim} (\int^{\mathbb{T}}X')_{h\mathbb{T}}\langle 1,\tilde{\theta}\rangle,
\end{equation}
where $\theta$ and $\tilde{\theta}$ are formal variables of degree $1$ and $-1$, respectively. Moreover, $\phi_p$ and $\tilde{\phi}_p$. satisfy the following properties:
\begin{enumerate}[label=\arabic*)]
    \item  $\phi_p$ and $\tilde{\phi}_p$ are natural in $X$ and $X'$, respectively,
    \item for a cyclic $k$-module $X$ and a chain complex $K$, the diagram
\begin{equation}
\begin{tikzcd}[row sep=1.2cm, column sep=0.8cm]
(\int^{\mathbb{Z}/p}Map(j^*M,K))_{h\mathbb{Z}/p}\arrow[r] \arrow[d,"{\tilde{\phi}_p}"]& Map((\int_{\mathbb{Z}/p}j^*M)^{h\mathbb{Z}/p},K) \arrow[d,"{Map({\phi_p},K)}"]\\
(\int^{\mathbb{T}}Map(M,K))_{h\mathbb{T}}\langle 1,\tilde{\theta}\rangle\arrow[r] & Map((\int_{\mathbb{T}}M)^{h\mathbb{T}}\langle 1,\theta\rangle,K)
\end{tikzcd}
\end{equation}
is commutes up to homotopy.
     \item let $\underline{k}:N\Lambda^{op}\rightarrow Mod_k$ be the constant cocyclic module with value $k$. Under the canonical identifications
\begin{equation}
H^*((\int^{\mathbb{T}}\underline{k})_{h\mathbb{T}}) \cong H_{-*}(BS^1)=k[\tilde{t}], |t|=-2
\end{equation}
and
\begin{equation}
H^*((\int^{\mathbb{Z}/p}\underline{k})_{h\mathbb{Z}/p})\cong H_{-*}(B\mathbb{Z}/p)=k[\tilde{t},\tilde{\theta}], |\tilde{t}|=-2, |\tilde{\theta}|=-1, \tilde{\theta}^2=0,
\end{equation}
the Gysin comparison map $\tilde{\phi}_p$ for $\underline{k}$ is the identity on $k[\tilde{t},\tilde{\theta}]$.
\end{enumerate}
\end{prop}

\noindent\emph{Proof}.  By Lemma 5.2, we have the following commutative diagram of $\infty$-categories (and its opposite version)
\begin{equation}
\begin{tikzcd}[row sep=1.2cm, column sep=0.8cm]
\mathrm{Fun}(N\Lambda,Mod_k)\arrow[r,"{j^*}"]\arrow[d,"{\int_{\mathbb{T}}}"]&\mathrm{Fun}(N{}_p\Lambda, Mod_k)\arrow[d,"{\int_{\mathbb{Z}/p}}"]\\
\mathrm{Fun}(B\mathbb{T},Mod_k)\arrow[r,"{\mathfrak{j}^*}"]&\mathrm{Fun}(B\mathbb{Z}/p,Mod_k)
\end{tikzcd}
\end{equation}
where $\mathfrak{j}:\mathbb{Z}/p\subset\mathbb{T}$ denotes the standard inclusion. The canonical equivalence of $\infty$-categories
\begin{equation}
\mathrm{Fun}(B\mathbb{T}, Mod_k)\simeq Mod_{k[\mathbb{T}]}
\end{equation}
gives rise to commutative diagrams
\begin{equation}
\begin{tikzcd}[row sep=1.2cm, column sep=0.8cm]
\mathrm{Fun}(B\mathbb{T},Mod_k)\arrow[r,"{(-)_{h\mathbb{T}}}"]\arrow[d,"\cong"]&Mod_k\arrow[d,"="]\\
Mod_{k[\mathbb{T}]}\arrow[r,"{k\otimes_{k[\mathbb{T}]}-}"]&Mod_k,
\end{tikzcd}
\begin{tikzcd}[row sep=1.2cm, column sep=0.8cm]
\mathrm{Fun}(B\mathbb{T},Mod_k)\arrow[r,"{(-)^{h\mathbb{T}}}"]\arrow[d,"\cong"]&Mod_k\arrow[d,"="]\\
Mod_{k[\mathbb{T}]}\arrow[r,"{\mathrm{Hom}_{k[\mathbb{T}]}(k,-)}"]&Mod_k.
\end{tikzcd}
\end{equation}
Similarly, the canonical equivalence
\begin{equation}
\mathrm{Fun}(B\mathbb{Z}/p,Mod_k)\simeq Mod_{k[\mathbb{Z}/p]}
\end{equation}
gives rise to commutative diagrams
\begin{equation}
\begin{tikzcd}[row sep=1.2cm, column sep=0.8cm]
\mathrm{Fun}(B\mathbb{Z}/p,Mod_k)\arrow[r,"{(-)_{h\mathbb{Z}/p}}"]\arrow[d,"\simeq"]&Mod_k\arrow[d,"="]\\
Mod_{k[\mathbb{Z}/p]}\arrow[r,"{k\otimes_{k[\mathbb{Z}/p]}-}"]&Mod_k,
\end{tikzcd}
\begin{tikzcd}[row sep=1.2cm, column sep=0.8cm]
\mathrm{Fun}(B\mathbb{Z}/p,Mod_k)\arrow[r,"{(-)^{h\mathbb{Z}/p}}"]\arrow[d,"\simeq"]&Mod_k\arrow[d,"="]\\
Mod_{k[\mathbb{Z}/p]}\arrow[r,"{\mathrm{Hom}_{k[\mathbb{Z}/p]}(k,-)}"]&Mod_k.
\end{tikzcd}
\end{equation}
Moreover, these are related by the diagram
\begin{equation}
\begin{tikzcd}[row sep=1.2cm, column sep=0.8cm]
\mathrm{Fun}(B\mathbb{T},Mod_k)\arrow[d,"\simeq"]\arrow[r,"\mathfrak{j}^*"] &\mathrm{Fun}(B\mathbb{Z}/p,Mod_k)\arrow[d,"\simeq"]\\
Mod_{k[\mathbb{T}]}\arrow[r,"\mathfrak{j}^*"]&Mod_{k[\mathbb{Z}/p]}
\end{tikzcd}
\end{equation}
We will show that for $\mathfrak{X}\in Mod_{k[\mathbb{T}]}$, there exist quasi-isomorphisms
\begin{equation}
\mathrm{Hom}_{k[\mathbb{Z}/p]}(k,\mathfrak{j}^*\mathfrak{X})\xleftarrow{\simeq} \mathrm{Hom}_{k[\mathbb{T}]}(k,\mathfrak{X})\langle 1,\theta\rangle: \varphi_p
\end{equation}
and
\begin{equation}
\tilde{\varphi}_p: \mathfrak{j}^*\mathfrak{X}\otimes_{k[\mathbb{Z}/p]}k\xrightarrow{\simeq} \mathfrak{X}\otimes_{k[\mathbb{T}]}k\langle 1,\tilde{\theta}\rangle.
\end{equation}
Then, we define the $\mathbb{Z}/p$-Gysin comparison maps in (5.4) (resp. (5.5)) as
\begin{equation}
\phi_p:= \varphi_p\circ\int_{\mathbb{T}}\quad(\textrm{resp.}\;\tilde{\phi}_p:= \tilde{\varphi}_p\circ\int^{\mathbb{T}}).
\end{equation}
To construct the maps in (5.15) and (5.16), we write down explicit complexes computing the corresponding homotopy orbits/fixed points. Let $\tau\in\mathbb{Z}/p$ be the standard generator, then there is a free $k[\mathbb{Z}/p]$-resolution of $k$ given by
\begin{equation}
\cdots\xrightarrow{1+\tau+\cdots+\tau^{p-1}} k[\mathbb{Z}/p]\xrightarrow{\tau-1}k[\mathbb{Z}/p]\xrightarrow{} k\xrightarrow{} 0.
\end{equation}
Therefore, for $\mathfrak{Y}\in Mod_{k[\mathbb{Z}/p]}$, an explicit complex computing $\mathrm{Hom}_{k[\mathbb{Z}/p]}(k,\mathfrak{Y})$ is $\mathfrak{Y}[[t,\theta]]$ with $t$-linear differential
\begin{equation}
\begin{cases}
y\mapsto d_{\mathfrak{Y}}y+(-1)^{|y|}(\tau-1)y\theta\\
y\theta\mapsto d_{\mathfrak{Y}}y\theta+(-1)^{|y|}(1+\tau+\cdots+\tau^{p-1})yt.
\end{cases}
\end{equation}
We now write down an explicit complex for $\mathrm{Hom}_{k[\mathbb{T}]}(k,\mathfrak{X})$. Consider the dg-algebra
\begin{equation}
k[\tau,\sigma]:=k[\tau,\sigma]/(\tau^p=1,d\sigma=\tau-1,\sigma^2=0),
\end{equation}
where $\tau$ is in degree $0$ and $\sigma$ is in degree $-1$. There is a canonical equivalence of $\mathcal{A}_{\infty}$-$k$-ring spectra $k[\tau,\sigma]\cong k[\mathbb{T}]$,
coming from the cellular structure on $S^1$ with $p$ 0-cells and $p$ 1-cells, which induces a diagram of $\infty$-categories
\begin{equation}
\begin{tikzcd}[row sep=1.2cm, column sep=0.8cm]
k[\mathbb{Z}/p]\arrow[r,"{\tau\mapsto \tau}"]\arrow[dr,"\mathfrak{i}"]& k[\tau,\sigma]\arrow[d,"\cong"]\\
&k[\mathbb{T}].
\end{tikzcd}
\end{equation}
Moreover, under this identification, $\mathrm{Hom}_{k[\mathbb{T}]}(k,-)$ is identified with $\mathrm{Hom}_{k[\tau,\sigma]}(k,-)$. There is a free $k[\tau,\sigma]$-resolution of $k$ given by
\begin{equation}
\cdots\xrightarrow{(1+\tau+\cdots+\tau^{p-1})\sigma} k[\tau,\sigma]\xrightarrow{(1+\tau+\cdots+\tau^{p-1})\sigma} k[\tau,\sigma]\rightarrow k\rightarrow 0.
\end{equation}
Therefore, an explicit complex computing $\mathrm{Hom}_{k[\tau,\sigma]}(k,\mathfrak{X})\cong \mathrm{Hom}_{k[\mathbb{T}]}(k,\mathfrak{X})$ is $\mathfrak{X}[[t]]$ with $t$-linear differential given by
\begin{equation}
x\mapsto d_{\mathfrak{X}}x+(1+\tau+\cdots+\tau^{p-1})\sigma xt.
\end{equation}
Using the explicit descriptions above, one finds a $t$-linear map
\begin{equation}
\varphi_p:\mathrm{Hom}_{k[\tau,\sigma]}(k,\mathfrak{X})\langle 1,\theta\rangle\rightarrow \mathrm{Hom}_{k[\mathbb{Z}/p]}(k,\mathfrak{j}^*\mathfrak{X})
\end{equation}
given by
\begin{equation}
\begin{cases}
x\mapsto x+(-1)^{|x|}\sigma x\theta\\
x\theta\mapsto x\theta+(-1)^{|x|}(\tau-1)^{p-2}\sigma xt.
\end{cases}
\end{equation}
$\varphi_p$ is a quasi-isomorphism because its constant term is the identity. \par\indent
Similarly, an explicit chain complex computing $\mathfrak{j}^*\mathfrak{X}\otimes_{k[\mathbb{Z}/p]} k$ is $\mathfrak{X}[\tilde{t},\tilde{\theta}], |\tilde{t}|=-2,|\tilde{\theta}|=-1,\tilde{\theta}^2=0$, with differential given by
\begin{equation}
\begin{cases}
x\mapsto d_{\mathfrak{X}}(x)\\
x\tilde{t}^k\mapsto d_{\mathfrak{X}}(x)\tilde{t}^k+(-1)^{|x|}(1+\tau+\cdots+\tau^{p-1})x\tilde{t}^{k-1}\tilde{\theta},\;\;k>0\\
x\tilde{t}^k\tilde{\theta}\mapsto d_{\mathfrak{X}}(x)\tilde{t}^k\tilde{\theta}+(-1)^{|x|}(\tau-1)\tilde{t}^k,\;\;k\geq 0.
\end{cases}
\end{equation}
An explicit chain complex computing $\mathfrak{X}\otimes_{k[\mathbb{T}]}k$ is given by $\mathfrak{X}[\tilde{t}], |\tilde{t}|=-2$, with differential given by
\begin{equation}
\begin{cases}
x\mapsto d_{\mathfrak{X}}(x)\\
x\tilde{t}^k\mapsto d_{\mathfrak{X}}(x)\tilde{t}^k+(1+\tau+\cdots+\tau^{p-1})\sigma \tilde{t}^{k-1},\;\;k>0
\end{cases}
\end{equation}
In this case, we define the $\mathbb{Z}/p$-Gysin comparison map $\tilde{\varphi}_p: \mathfrak{j}^*\mathfrak{X}\otimes_{k[\mathbb{Z}/p]}k\xrightarrow{\simeq} \mathfrak{X}\otimes_{k[\mathbb{T}]}k\langle 1,\tilde{\theta}\rangle$ by
\begin{equation}
\begin{cases}
x\mapsto x\\
x\tilde{t}^k\mapsto x\tilde{t}^k+(-1)^{|x|}(\tau-1)^{p-2}\sigma x\tilde{t}^{k-1}\tilde{\theta},\;\;k>0\\
x\tilde{t}^k\tilde{\theta}\mapsto x\tilde{t}^k\tilde{\theta}+(-1)^{|x|}\sigma x\tilde{t}^k.
\end{cases}
\end{equation}
Property 1) is clear from the definition. \par\indent
For Property 2): the explicit formulas for $\phi_p, \tilde{\phi}_p$ show that for $X\in Mod_{k[\mathbb{T}]}$, there is a commutative diagram of chain complexes
\begin{equation}
\begin{tikzcd}[row sep=1.2cm, column sep=0.8cm]
Map(\mathfrak{i}^*X,K)_{h\mathbb{Z}/p}\arrow[r]\arrow[d,"{\tilde{\phi}_p}"]& Map(\mathfrak{i}^*X^{h\mathbb{Z}/p},K)\arrow[d,"{Map({\phi_p},K)}"]\\
Map(X,K)_{h\mathbb{T}}\langle 1,\tilde{\theta}\rangle \arrow[r]& Map(X^{h\mathbb{T}}\langle 1,\theta\rangle,K)
\end{tikzcd}.
\end{equation}
On the other hand, for $X\in \mathrm{Fun}(N\Lambda,Mod_k)$, there is a natural equivalence of $\mathbb{T}$-modules
\begin{equation}
\int^{\mathbb{T}}Map(X,K)\xleftarrow{\simeq} Map(\int_{\mathbb{T}}X,K).
\end{equation}
By Lemma 3.10, this is a weak equivalence (since weak equivalences can be detected by the forgetful functor to $\mathrm{Mod}_k$). Combining (5.29) and (5.30) we conclude Property 2).\par\indent
3): It follows from Proposition 3.24 that as a $\mathbb{T}$-module, $\int^{\mathbb{T}}\underline{k}$ is equivalent to the chain complex $k$ equipped with the trivial $\mathbb{T}$-action.
In particular, under the identification $k[\mathbb{T}]\simeq k[\tau,\sigma]$, it is equivalent to the $k[\tau,\sigma]$-module $k$ where $\tau$ acts by identity and $\sigma$ acts by zero.
Therefore, $(\int^{\mathbb{T}}\underline{k})_{h\mathbb{T}}$ is computed by the total complex of
\begin{equation}
\big(\cdots\xrightarrow{(1+\tau+\cdots+\tau^{p-1})\sigma} k[\tau,\sigma]\xrightarrow{(1+\tau+\cdots+\tau^{p-1})\sigma} k[\tau,\sigma]\rightarrow 0\big)\otimes_{k[\tau,\sigma]} k,
\end{equation}
whose cohomology is naturally isomorphic to $k[\tilde{t}]$. On the other hand, $(\int^{\mathbb{Z}/p}\underline{k})_{h\mathbb{Z}/p}$ is computed by the total complex of
\begin{equation}
(\cdots\xrightarrow{1+\tau+\cdots+\tau^{p-1}} k[\mathbb{Z}/p]\xrightarrow{\tau-1}k[\mathbb{Z}/p]\xrightarrow{} 0)\otimes_{k[\mathbb{Z}/p]}k,
\end{equation}
whose cohomology (in characteristic $p$) is naturally isomorphic to $k[\tilde{t},\tilde{\theta}]$. By the explicit formula defining $\varphi_p$, one easily sees that the comparison map $\tilde{\phi}_p$ is the identity on $k[\tilde{t},\tilde{\theta}]$. \qed\par\indent
5.2. \textbf{Compatibility of operadic (negative) cyclic and $\mathbb{Z}/p$-equivariant open-closed maps}.
\begin{thm}
Let $A_L=CF^*(L,L)$. The following diagram of chain complexes is homotopy commutative
\begin{equation}
\begin{tikzcd}[row sep=1.2cm, column sep=0.8cm]
\int_{\mathbb{Z}/p}j^*(\iota_{\Lambda})_!(\tilde{P}^{-1}_{\overrightarrow{\Lambda}})^*A^{\sharp}_L\simeq CC^{\mathbb{Z}/p,oper}(A_L)\arrow[rrr,"{OC^{\mathbb{Z}/p,oper}}"]& & &QH(M)[[t,\theta]]  \\
(\int_{\mathbb{T}}(\iota_{\Lambda})_!(\tilde{P}^{-1}_{\overrightarrow{\Lambda}})^*A^{\sharp}_L)\langle 1,\theta\rangle=CC^{S^1,oper}(A_L)\langle 1,\theta\rangle \arrow[rrr,"{OC^{S^1,oper}\langle 1,\theta\rangle}"]\arrow[u,"{\Phi_p^{oper}}"]& & & QH(M)[[t,\theta]]\arrow[u,"="]
\end{tikzcd},
\end{equation}
where $\Phi_p^{oper}$ is the $\mathbb{Z}/p$-Gysin comparison map from Proposition 5.3 applied to the cyclic $k$-module $(\iota_{\Lambda})_!(\tilde{P}^{-1}_{\overrightarrow{\Lambda}})^*A^{\sharp}_L$.
\end{thm}

\noindent\emph{Proof}.
By contractibility of $\overline{\mathcal{R}}^1_{d+1}$, there is an equivalence $(\iota^{op}_{\Lambda})_*((\tilde{P}^{op}_{\overrightarrow{\Lambda}})^{-1})^*\overline{\mathfrak{R}}^1\simeq \underline{k}$ of cocyclic $k$-modules. Thus by Proposition 5.3 3), there is
a commutative diagram
\begin{equation}
\begin{tikzcd}[row sep=1.2cm, column sep=0.8cm]
k[\tilde{t},\tilde{\theta}]\arrow[r,"\simeq"]\arrow[d,"="]& H^*(\Big(\int^{\mathbb{Z}/p}(j^{op})^*(\iota^{op}_{\Lambda})_*((\tilde{P}^{op}_{\overrightarrow{\Lambda}})^{-1})^*\overline{\mathfrak{R}}^1\Big)_{h\mathbb{Z}/p})\arrow[d,"{H^*(\tilde{\phi}_p)}"] \\
k[\tilde{t},\tilde{\theta}]\arrow[r,"\simeq"]& H^*(\Big(\int^{\mathbb{T}}(\iota^{op}_{\Lambda})_*((\tilde{P}^{op}_{\overrightarrow{\Lambda}})^{-1})^*\overline{\mathfrak{R}}^1\Big)_{h\mathbb{T}})\langle 1,\tilde{\theta}\rangle
\end{tikzcd}.
\end{equation}
For any $X\in \mathrm{Fun}^{dg}(\overrightarrow{\Lambda}^{op}\rtimes \mathcal{A}^{oper,dg}_{\infty},\mathrm{Mod}_k)$,
the Beck-Chevalley transform (cf. Appendix D.3) gives rise to a morphism
\begin{equation}
(\iota^{op}_{{}_p\Lambda})_*((\tilde{P}^{op}_{\overrightarrow{{}_p\Lambda}})^{-1})^*(\overrightarrow{j^{op}_{\mathcal{A}_{\infty}^{oper,dg}}})^*X\leftarrow (j^{op})^*(\iota^{op}_{\Lambda})_*((\tilde{P}^{op}_{\overrightarrow{\Lambda}})^{-1})^*X
\end{equation}
in $\mathrm{Fun}(N{}_p\Lambda,\mathrm{Mod}_k)$. Apply the right Kan extension $\int^{\mathbb{Z}/p}$ we obtain a map of $\mathbb{Z}/p$-modules
\begin{equation}
\int^{\mathbb{Z}/p}(\iota^{op}_{{}_p\Lambda})_*((\tilde{P}^{op}_{\overrightarrow{{}_p\Lambda}})^{-1})^*(\overrightarrow{j^{op}_{\mathcal{A}_{\infty}^{oper,dg}}})^*X\leftarrow \int^{\mathbb{Z}/p}(j^{op})^*(\iota^{op}_{\Lambda})_*((\tilde{P}^{op}_{\overrightarrow{\Lambda}})^{-1})^*X.
\end{equation}
If $X$ is H-counital, analogous to Corollary 3.13 one can show that (5.36) is an equivalence. Apply this to $X=\mathrm{Map}(A^{\sharp}_L,QH)$ and $X=\overline{\mathfrak{R}}^1$, and combining with
Proposition 5.3 1) applied to the map of cocyclic $k$-modules
\begin{equation}
(\iota^{op}_{\Lambda})_*((\tilde{P}^{op}_{\overrightarrow{\Lambda}})^{-1})^*OC^{oper}: (\iota^{op}_{\Lambda})_*((\tilde{P}^{op}_{\overrightarrow{\Lambda}})^{-1})^*\overline{\mathfrak{R}}^1\rightarrow (\iota^{op}_{\Lambda})_*((\tilde{P}^{op}_{\overrightarrow{\Lambda}})^{-1})^*Map(A^{\sharp}_L,QH),
\end{equation}
we obtain a homotopy commutative diagram
\begin{equation}
\begin{tikzcd}[row sep=1.2cm, column sep=0.8cm]
\Big(\int^{\mathbb{Z}/p}(\iota^{op}_{{}_p\Lambda})_*((\tilde{P}^{op}_{\overrightarrow{{}_p\Lambda}})^{-1})^*(\overrightarrow{j^{op}_{\mathcal{A}_{\infty}^{oper,dg}}})^*\overline{\mathfrak{R}}^1\Big)_{h\mathbb{Z}/p}\arrow[rr]\arrow[d,"{\tilde{\phi}_p}"] & &\Big(\int^{\mathbb{Z}/p}(\iota^{op}_{{}_p\Lambda})_*((\tilde{P}^{op}_{\overrightarrow{{}_p\Lambda}})^{-1})^*(\overrightarrow{j^{op}_{\mathcal{A}_{\infty}^{oper,dg}}})^*Map(A^{\sharp}_L,QH)\Big)_{h\mathbb{Z}/p}\arrow[d,"{\tilde{\phi}_p}"]\\
\Big(\int^{\mathbb{T}}(\iota^{op}_{\Lambda})_*((\tilde{P}^{op}_{\overrightarrow{\Lambda}})^{-1})^*\overline{\mathfrak{R}}^1\Big)_{h\mathbb{T}}\langle 1,\tilde{\theta}\rangle\arrow[rr] & &
\Big(\int^{\mathbb{T}}(\iota^{op}_{\Lambda})_*((\tilde{P}^{op}_{\overrightarrow{\Lambda}})^{-1})^*Map(A^{\sharp}_L,QH)\Big)_{h\mathbb{T}}\langle 1,\tilde{\theta}\rangle,
\end{tikzcd}
\end{equation}
where the horizontal arrows are induced by $OC^{oper}$. \par\indent
On the other hand, there is a natural equivalence of $\mathbb{T}$-modules \begin{equation}
\int^{\mathbb{T}}(\iota^{op}_{\Lambda})_*((\tilde{P}^{op}_{\overrightarrow{\Lambda}})^{-1})^*Map(A^{\sharp}_L,QH)\simeq Map(\int_{\mathbb{T}}(\iota_{\Lambda})_!(\tilde{P}^{-1}_{\overrightarrow{\Lambda}})^*A^{\sharp}_L,QH).
\end{equation}
By Proposition 5.3 1) again, there is a homotopy commutative diagram
\begin{equation}
\begin{tikzcd}[row sep=1.2cm, column sep=0.8cm]
\Big(\int^{\mathbb{Z}/p}(\iota^{op}_{{}_p\Lambda})_*((\tilde{P}^{op}_{\overrightarrow{{}_p\Lambda}})^{-1})^*(\overrightarrow{j^{op}_{\mathcal{A}_{\infty}^{oper,dg}}})^*Map(A^{\sharp}_L,QH)\Big)_{h\mathbb{Z}/p}\arrow[r,"\simeq"]\arrow[d,"{\tilde{\phi}_p}"]&{Map(\int_{\mathbb{Z}/p}(\iota_{{}_p\Lambda})_!(\tilde{P}^{-1}_{\overrightarrow{{}_p\Lambda}})^*(\overrightarrow{j_{\mathcal{A}_{\infty}^{oper,dg}}})^*A^{\sharp}_L,QH)}_{h\mathbb{Z}/p}\arrow[d,"{\tilde{\phi}_p}"]\\
\Big(\int^{\mathbb{T}}(\iota^{op}_{\Lambda})_*((\tilde{P}^{op}_{\overrightarrow{\Lambda}})^{-1})^*Map(A^{\sharp}_L,QH)\Big)_{h\mathbb{T}}\langle 1,\tilde{\theta}\rangle\arrow[r,"\simeq"]&{Map(\int_{\mathbb{T}}(\iota_{\Lambda})_!(\tilde{P}^{-1}_{\overrightarrow{\Lambda}})^*A^{\sharp}_L,QH)}_{h\mathbb{T}}\langle 1,\tilde{\theta}\rangle.
\end{tikzcd}.
\end{equation}
By Proposition 5.3 2) there is a homotopy commutative diagram
\begin{equation}
\begin{tikzcd}[row sep=1.2cm, column sep=0.8cm]
{Map(\int_{\mathbb{Z}/p}(\iota_{{}_p\Lambda})_!(\tilde{P}^{-1}_{\overrightarrow{{}_p\Lambda}})^*(\overrightarrow{j_{\mathcal{A}_{\infty}^{oper,dg}}})^*A^{\sharp}_L,QH)}_{h\mathbb{Z}/p}\arrow[r]\arrow[d,"{\tilde{\phi}_p}"]& Map(\Big(\int_{\mathbb{Z}/p}(\iota_{{}_p\Lambda})_!(\tilde{P}^{-1}_{\overrightarrow{{}_p\Lambda}})^*(\overrightarrow{j_{\mathcal{A}_{\infty}^{oper,dg}}})^*A^{\sharp}_L\Big)^{h\mathbb{Z}/p},QH)\arrow[d,"{Map(\phi_p,QH)}"]\\
{Map(\int_{\mathbb{T}}(\iota_{\Lambda})_!(\tilde{P}^{-1}_{\overrightarrow{\Lambda}})^*A^{\sharp}_L,QH)}_{h\mathbb{T}}\langle 1,\tilde{\theta}\rangle\arrow[r] & Map(\Big(\int_{\mathbb{T}}(\iota_{\Lambda})_!(\tilde{P}^{-1}_{\overrightarrow{\Lambda}})^*A^{\sharp}_L\Big)^{h\mathbb{T}}\langle 1,\theta\rangle,QH).
\end{tikzcd}.
\end{equation}
Combining the diagrams (5.34),(5.38), (5.40) and (5.41), one concludes Theorem 5.4. To be more precise, take the element $\tilde{t}^m\in k[\tilde{t},\tilde{\theta}]$ and trace its image along the
 bottom rows of (5.34),(5.38), (5.40) and finally (5.41), we obtain the chain map $OC^{S^1,oper}_{m}$, cf. Definition 4.9. Similarly, take the element $\tilde{t}^m\;(\textrm{resp}.\;\tilde{t}^m\tilde{\theta})\in k[\tilde{t},\tilde{\theta}]$
 and trace its image along the top rows of (5.34), (5.38), (5.40) and finally (5.41), we obtain the chain map $OC^{\mathbb{Z}/p,oper}_{2m}$ (resp. $OC^{\mathbb{Z}/p,oper}_{2m+1}$), cf. Definition 4.10. \qed\par\indent
 5.3. \textbf{Proof of Proposition 1.5}. We construct the desired quasi-isomorphism $\Phi_p$ as follows. Let $A$ be an H-unital $\mathcal{A}_{\infty}$-algebra. By Proposition 3.11 and 3.19, we have natural quasi-isomorphisms
 \begin{equation}
CC^{S^1}(A)\xrightarrow[\simeq]{(3.21)} \big(\int_{\mathbb{T}}(\iota_{\Lambda})_!(P^{-1}_{\overrightarrow{\Lambda}})^*A^{\sharp})\big)^{h\mathbb{T}}
\end{equation}
and
\begin{equation}
CC^{\mathbb{Z}/p}(A)\xrightarrow[\simeq]{(3.51)} \big(\int_{\mathbb{Z}/p}(\iota_{{}_p\Lambda})_!(P^{-1}_{\overrightarrow{{}_p\Lambda}})^*(\overrightarrow{j_{\mathcal{A}_{\infty}^{dg}}})^*A^{\sharp})\big)^{h\mathbb{Z}/p}.
 \end{equation}
Since $A$ is cohomological unital, this implies that $A^{\sharp}$ is H-unital and hence the map induced by Beck-Chevalley transform
\begin{equation}
\big(\int_{\mathbb{Z}/p}(\iota_{{}_p\Lambda})_!(P^{-1}_{\overrightarrow{{}_p\Lambda}})^*(\overrightarrow{j_{\mathcal{A}_{\infty}^{dg}}})^*A^{\sharp})\big)^{h\mathbb{Z}/p}\xleftarrow{\simeq} \big(\int_{\mathbb{Z}/p}j^*(\iota_{\Lambda})_!(P^{-1}_{\overrightarrow{\Lambda}})^*A^{\sharp})\big)^{h\mathbb{Z}/p}
\end{equation}
is a quasi-isomorphism. Then, we define $\Phi_p$ to be the composition
\begin{equation*}
\Phi_p: CC^{S^1}(A)\langle 1,\theta\rangle\xrightarrow[\simeq]{(3.21)}  \big(\int_{\mathbb{T}}(\iota_{\Lambda})_!(P^{-1}_{\overrightarrow{\Lambda}})^*A^{\sharp})\big)^{h\mathbb{T}}\langle 1,\theta\rangle \xrightarrow[\simeq]{\phi_p} \big(\int_{\mathbb{Z}/p}j^*(\iota_{\Lambda})_!(P^{-1}_{\overrightarrow{\Lambda}})^*A^{\sharp})\big)^{h\mathbb{Z}/p}
\end{equation*}
\begin{equation}
\xrightarrow[\simeq]{(5.44)} \big(\int_{\mathbb{Z}/p}(\iota_{{}_p\Lambda})_!(P^{-1}_{\overrightarrow{{}_p\Lambda}})^*(\overrightarrow{j_{\mathcal{A}_{\infty}^{dg}}})^*A^{\sharp})\big)^{h\mathbb{Z}/p}\xrightarrow[\simeq]{(3.51)^{-1}}  CC^{\mathbb{Z}/p}(A),
\end{equation}
where $\phi_p$ denotes $\mathbb{Z}/p$-Gysin comparison map in Proposition 5.3 applied to the cyclic $k$-module $(\iota_{\Lambda})_!(P^{-1}_{\overrightarrow{\Lambda}})^*A^{\sharp}$. \qed

\renewcommand{\theequation}{6.\arabic{equation}}
\setcounter{equation}{0}

\section{Comparison with the classical construction}
In this section, we prove that the operadic cyclic and $\mathbb{Z}/p$-equivariant open-closed map agree with their classical counterparts.\par\indent
To be more precise, recall from Lemma 4.7 that $A_L=CF^*(L,L)$ is an algebra over the operad $\mathcal{A}_{\infty}^{oper,dg}$.
In section 4.3, we studied the operadic Hochschild functor
\begin{equation}
A^{\sharp}_L: \overrightarrow{\Lambda}\rtimes \mathcal{A}_{\infty}^{oper,dg}\rightarrow \mathrm{Mod}_k,
\end{equation}
and defined the operadic cyclic homology chain complex as
\begin{equation}
CC^{S^1,oper}(A_L):=\big(\int_{\mathbb{T}}(\iota_{\Lambda})_!(\tilde{P}^{-1}_{\overrightarrow{\Lambda}})^*A^{\sharp}_L\big)^{h\mathbb{T}}.
\end{equation}
In this section, we prove the following theorem.
\begin{thm}
1) There exists a quasi-isomorphism
\begin{equation}
\Xi^{S^1}:CC^{S^1,oper}(A_L)\simeq CC^{S^1}(A_L),
\end{equation}
where $CC^{S^1}(A_L)$ denotes the negative cyclic chain complex of (2.15) associated to $A_L$ (viewed as an $\mathcal{A}_{\infty}$-algebra).
Moreover, $\Xi^{S^1}$ induces a homotopy commutative diagram
\begin{equation}
\begin{tikzcd}[row sep=1.2cm, column sep=0.8cm]
CC^{S^1,oper}(A_L)\arrow[rr,"{OC^{S^1,oper}}"]\arrow[d,"{\Xi^{S^1}}"]& &QH(M)[[t]]\\
CC^{S^1}(A_L)\arrow[urr,"OC^{S^1}"]& &
\end{tikzcd},
\end{equation}
\end{thm}
where $OC^{S^1}$ is Ganatra's cyclic open-closed map (2.16) and $OC^{S^1,oper}$ is the operadic cyclic open-closed map (4.33). \par\indent
2) Similarly for the $\mathbb{Z}/p$-equivariant case, there exists a quasi-isomorphism
\begin{equation}
\Xi^{\mathbb{Z}/p}: CC^{\mathbb{Z}/p,oper}(A_L)\simeq CC^{\mathbb{Z}/p}(A_L)
\end{equation}
which makes the following diagram homotopy commute
\begin{equation}
\begin{tikzcd}[row sep=1.2cm, column sep=0.8cm]
CC^{\mathbb{Z}/p,oper}(A_L)\arrow[rr,"{OC^{\mathbb{Z}/p,oper}}"]\arrow[d,"{\Xi^{\mathbb{Z}/p}}"]& &QH(M)[[t,\theta]]\\
CC^{\mathbb{Z}/p}(A_L)\arrow[urr,"OC^{\mathbb{Z}/p}"]& &
\end{tikzcd},
\end{equation}
where $OC^{\mathbb{Z}/p}$ is the $\mathbb{Z}/p$-equivariant open-closed map (2.36) and $OC^{\mathbb{Z}/p,oper}$ is operadic $\mathbb{Z}/p$-equivariant open-closed map (4.38). \par\indent
3) There is a homotopy commutative diagram
\begin{equation}
\begin{tikzcd}[row sep=1.2cm, column sep=0.8cm]
CC^{S^1,oper}(A_L)\langle 1,\theta\rangle\arrow[rr,"\Phi_p^{oper}"]\arrow[d,"{\Xi^{S^1}\langle1,\theta\rangle}"]& &CC^{\mathbb{Z}/p,oper}(A_L)\arrow[d,"\Xi^{\mathbb{Z}/p}"]\\
CC^{S^1}(A_L)\langle 1,\theta\rangle\arrow[rr,"\Phi_p"]& & CC^{\mathbb{Z}/p}(A_L)
\end{tikzcd}
\end{equation}
\emph{Proof of Theorem 1.6}. Theorem 1.6 follows immediately from Theorem 6.1 and Theorem 5.4. Explicitly, Consider the following diagram
\begin{equation}
\begin{tikzcd}[row sep=1.2cm, column sep=0.8cm]
 &CC^{S^1,oper}(A_L)\langle 1,\theta\rangle\arrow[rr,"\Phi_p^{oper}"]\arrow[ddl,"{OC^{S^1,oper}\langle1, \theta\rangle}"']\arrow[d,"{\Xi^{S^1}\langle1,\theta\rangle}"]& &CC^{\mathbb{Z}/p,oper}(A_L)\arrow[d,"\Xi^{\mathbb{Z}/p}"]\arrow[ddr,"{OC^{\mathbb{Z}/p,oper}}"]&\\
&CC^{S^1}(A_L)\langle 1,\theta\rangle\arrow[rr,"\Phi_p"]\arrow[dl,"{OC^{S^1}\langle1,\theta\rangle}"]& & CC^{\mathbb{Z}/p}(A_L)\arrow[dr,"OC^{\mathbb{Z}/p}"'] & \\
QH[[t,\theta]]\arrow[rrrr,"="]& & & & QH[[t,\theta]]
\end{tikzcd}
\end{equation}
The top square is homotopy commutative by Theorem 6.1. 3). The left and right triangles are homotopy commutative by Theorem 6.1. 1) and 2), respectively. The outer square is homotopy commutative by Theorem 5.4. Since $\Xi^{S^1}$ and $\Xi^{\mathbb{Z}/p}$ are homotopy equivalences, this implies that the bottom square of (6.8) is also homotopy commutative. This conludes the proof for the one object scenario, and generalization to the multiple object case will be discussed in Appendix A.\qed\par\indent
6.1. \textbf{Comparing the operads $\mathcal{A}_{\infty}^{dg}$ and $\mathcal{A}_{\infty}^{oper,dg}$}. Classically, for an object $L\in \mathrm{Fuk}(M)_{\lambda}$, its Floer
complex $A_L=CF^*(L,L)$ can be endowed with the structure of an $\mathcal{A}_{\infty}$-algebra, or equivalently, an algebra over the operad $\mathcal{A}_{\infty}^{dg}$. This structure is obtained from
\emph{a choice} of consistent Floer data on the moduli spaces of disks $\{\overline{\mathcal{R}}^{d+1}\}_{d\geq 2}$, cf. section 2.1. \par\indent
On the other hand, the action of the operad $\mathcal{A}_{\infty}^{oper,dg}$ on $A=CF^*(L,L)$ is independent of choices, or rather, the space of all consistent choices of Floer data is already built
into the definition of $\mathcal{A}_{\infty}^{oper,dg}$. As we will see in the next lemma, a consistent choice of Floer data over $\{\overline{\mathcal{R}}^{d+1}\}_{d\geq 2}$ gives rise to
a homotopy equivalence of dg operads $\tilde{s}: \mathcal{A}_{\infty}^{dg}\rightarrow \mathcal{A}_{\infty}^{oper,dg}$.

\begin{lemma}
Let $\mathcal{A}_{\infty}^{\Box,dg}:=\{C^{\Box}_*(\mathcal{A}_{\infty})\}_{d\geq 1}$ denote the dg operad of normalized smooth symmetric cubical chains on topological associahedron $\mathcal{A}_{\infty}$,
which naturally receives a projection map $\pi: \mathcal{A}_{\infty}^{oper,dg}\rightarrow \mathcal{A}_{\infty}^{\Box,dg}$ induced by the projections $C_*(\tilde{\mathcal{F}}^{reg}(\overline{\mathcal{R}}^{d+1}))\rightarrow C^{\Box}_*(\overline{\mathcal{R}}^{d+1})$ (by results in section 4.2, $\pi$ is a homotopy equivalence of dg operads).\par\indent
Then there exist
\begin{enumerate}[label=\arabic*)]
    \item a universal choice of Floer data for the spaces $\{\overline{\mathcal{R}}^{d+1}\}_{d\geq 2}$, compatible with products and boundary structures;
    \item for each cell $C$ in the standard cellular structure of $\overline{\mathcal{R}}^{d+1}$, the assignment of a cubical subdivision $s_d(C)$ of $C$. Moreover, viewing $s(C)$ as an element
of $C^{\Box}_*(\overline{\mathcal{R}}^{d+1})$, these assignments fit together into a homotopy equivalence of dg operads
\begin{equation}
s: \mathcal{A}_{\infty}^{dg}\rightarrow \mathcal{A}_{\infty}^{\Box,dg};
\end{equation}
    \item a homotopy equivalence of dg operads $\tilde{s}: \mathcal{A}_{\infty}^{dg}\rightarrow \mathcal{A}_{\infty}^{oper,dg}$ fitting into the diagram of dg operads
\begin{equation}
\begin{tikzcd}[row sep=1.2cm, column sep=0.8cm]
 & \mathcal{A}_{\infty}^{oper,dg}\arrow[d,"{\pi}"]\\
\mathcal{A}_{\infty}^{dg}\arrow[r,"{s}"] \arrow[ur,"{\tilde{s}}"] &\mathcal{A}_{\infty}^{\Box,dg}
\end{tikzcd}
\end{equation}
such that the Floer data on each cube in $\tilde{s}(C)$ is the universal choice of Floer data in 1) restricted to the corresponding cube of Riemann surfaces in $s(C)$.
\end{enumerate}
\end{lemma}
\noindent\emph{Proof}. First we construct the cubical subdivision of 2)
\begin{equation}
s_d: C_*^{cell}(\overline{\mathcal{R}}^{d+1})\rightarrow C_*^{\Box}(\overline{\mathcal{R}}^{d+1})
\end{equation}
inductively, for $d\geq 2$. These will be chain maps and are moreover compatible with boundary and product structures of the $\overline{\mathcal{R}}^{d+1}$'s.
When $d=2$, we take $s_0$ to be the trivial cubical subdivision of $\overline{\mathcal{R}}^2=pt$. Suppose we have constructed $s$ up to $s_{d-1}$, then the map $s_d$
is determined in degree $*<d-2$ (i.e. all except the top dimensional cell $[\mathcal{R}^{d+1}$) by operadic structures since the boundary of $\overline{\mathcal{R}}^{d+1}$
consists of products of lower dimensional $\overline{\mathcal{R}}^k$'s, see section 4.1.
We then define $s_d([\mathcal{R}^{d+1}])$ to be a cubical subdivision of the top dimensional cell that extends the cubical subdivision of the boundary.
Moreover, we require that for each face $f$ of a cube in $s_d([\mathcal{R}^{d+1}])$, its interior $f^o$ is contained in a single stratum of $\overline{\mathcal{R}}^{d+1}$.
This can be achieved by first extending the boundary cubical subdivision to a small tubular neighborhood $\mathcal{U}_{\epsilon}=\partial\overline{\mathcal{R}}^{d+1}\times [0,\epsilon]\rightarrow \overline{\mathcal{R}}^{d+1}$
of the boundary, by taking the product of the boundary cubical subdivision with $[0,\epsilon)$; then, we take an arbitrary cubical subdivision of the polytope $\overline{\mathcal{R}}^{d+1}\backslash \mathcal{U}_{\epsilon}$ that is compatible with the prior chosen subdivision of $\partial(\overline{\mathcal{R}}^{d+1}\backslash \mathcal{U}_{\epsilon})$. \par\indent
The construction of 1) is classical, but we will reproduce it, along with an inductive construction of $\tilde{s}$ in 3).
In fact, we first construct a lift of $s$ to
$\tilde{s}': C^{cell}_*(\overline{\mathcal{R}}^{d+1})\rightarrow C_*(\mathcal{F}(\overline{\mathcal{R}}^{d+1}))$, and then consider the necessary perturbations to obtain $\tilde{s}$.
For the base case, $\tilde{s}'_2$ is defined by choosing an arbitrary Floer data on the disk with two boundary inputs and one boundary output.
We now proceed to the inductive step and suppose we have constructed a choice of Floer data over the spaces $\overline{\mathcal{R}}^{k+1}$ as well as maps $\tilde{s}'_k: C^{cell}_*(\overline{\mathcal{R}}^{k+1})\rightarrow C_*(\mathcal{F}(\overline{\mathcal{R}}^{d+1}))$ lifting $s_k$, for $k<d$.
By the product structures at the boundary, this again determines $\tilde{s}'_d$ except for its value on the top dimensional cell $[\mathcal{R}^{d+1}]$.
Recall from the previous paragraph that $s([\mathcal{R}^{d+1}])\in C_*^{\Box}(\overline{\mathcal{R}}^{d+1})$ is a sum
$\sum_{i} C_i\times [0,\epsilon]+\sum_j D_j$, where the $d-3$-cubes $C_i$'s form a cubical subdivision of $\partial \overline{\mathcal{R}}^{d+1}$ and the $d-2$-cubes
$D_j$'s are entirely contained in the interior $\mathcal{R}^{d+1}$. By inductive hypothesis,
we've chosen, for each $i$, a cube $\eta_i\in C_{d-3}(\mathcal{F}(\partial\overline{\mathcal{R}}^{d+1}))$ which lifts $C_i$. For each $x\in [0,1]^{d-3}$, let $T_x$ denote the tree type of $C_i(x)$.
Take $\epsilon$ small enough, then there exists a unique map
\begin{equation}
l: [0,\epsilon]\rightarrow  [0,1]^{|E_{int}(T_{x})|}
\end{equation}
such that
\begin{equation}
\gamma_{l(t)}(x)=(x,t)\in \partial\overline{\mathcal{R}}^{d+1}\times[0,\epsilon]\subset\mathcal{R}^{d+1},
\end{equation}
where $\gamma$ denotes the gluing map for Riemann surfaces (without Floer data; compare (4.8)). This follows because $\gamma$ is a local diffeomorphism near a boundary stratum.
Now, define the map $b_i$ (see section 4.2) associated to the cube
\begin{equation}
\mathfrak{o}_i\in C_{d-2}(\mathcal{F}(\overline{\mathcal{R}}^{d+1}))
\end{equation}
that lifts $C_i\times [0,\epsilon]$ by
\begin{equation}
b_i(x,t):=\Gamma_{l(t)}(\eta_i(x)),
\end{equation}
where $\Gamma$ denotes the gluing of surfaces equipped with Floer data. Now we specify the gluing atlas $\{T_f,b_f,g_f\}$ (see (4.6), (4.7)) associated to the $\mathfrak{o}_i$.
By our assumption that $f$ is contained in a single stratum, there are two cases:
\begin{itemize}
    \item $f$ is a face such that $f^o\subset \mathcal{R}^{d+1}$. In this case $T_f$ is the unique $d$-leafed tree with one internal node. $g_f$ vanishes identically,
i.e. no gluing occurs. Moreover, the restriction of $b_i$ to $W_f$ determines a smooth family of unbroken disks with Floer data, which we define to be $b_f$.
     \item $f$ is a face such that $f^o\subset \partial \overline{\mathcal{R}}^{d+1}$. In this case, $f$ is a product of (faces of) lower dimensional cubes
coming from previously defined cubical subdivisons of lower dimensional associahedrons. $T_f$ is then defined to be the product of the corresponding tree types.
$b_f$, when restricted to $W_f\cap [0,1]^{d-3}\times \{0\}$, is a product of the form $b_{f_1}\times b_{f_2}$, where $b_{f_i}, i=1,2$ are already defined by inductive hypothesis; it is
constant along the last coordinate $[0,1]$. ${g_f}$, when restricted to $W_f\cap [0,1]^{d-3}\times \{0\}$, is a product of the form $g_{f_1}\times g_{f_2}$,
where $g_{f_i}, i=1,2$ are already defined by inductive hypothesis. In general, it is defined for $(x,t)\in [0,1]^{d-3}\times [0,\epsilon]$ by
\begin{equation}
g_f(x,t)=(g_{f_1}\times g_{f_2})(x)+l(t).
\end{equation}
\end{itemize}
These define lifts $\mathfrak{o}_i$ for all the `boundary cubes' $C_i\times [0,\epsilon]$. The above construction in particular gives an extension of the choice of Floer data on $\partial\overline{\mathcal{R}}^{d+1}$ to the tubular neighborhood $\mathcal{U}_{\epsilon}$, i.e. a section of the fibration
\begin{equation}
\mathcal{F}_{0,T_0}(\overline{\mathcal{R}}^{d+1})|_{\mathcal{U}_{\epsilon}}\rightarrow \mathcal{U}_{\epsilon},
\end{equation}
where $T_0$ denotes the unique $d$-leafed tree with one internal node (cf. section 4.2 for $\mathcal{F}_{0,T_0}$). Since the fibers of (6.17) are contractible, this section can be extended to a global section of
\begin{equation}
\mathcal{F}_{0,T_0}(\overline{\mathcal{R}}^{d+1})\rightarrow \mathcal{R}^{d+1}.
\end{equation}
Thus, we have inductively constructed the universal choice of Floer data over $\overline{\mathcal{R}}^{d+1}$.
Furthermore, this determines a lift $\tilde{D}_j\in C_{d-2}^*(\mathcal{F}(\overline{\mathcal{R}}^{d+1}))$ of each $D_j$: we set the tree type to be $T_f:=T_0$ for all face $f$ of $D_j$;
the gluing parameters $g_f$ vanish identically for all $f$; $b_f$ is given by the restriction of $\sigma$ to the face $f$ of $D_j$.\par\indent
Coming from the cubical subdivision, $s$ is a homotopy equivalence of dg operads. By Proposition 4.1, $\tilde{s}'$ is also a homotopy equivalence. \par\indent
Finally, we construct $\tilde{s}$. Recall that an $n$-cube of $\tilde{\mathcal{F}}(\overline{\mathcal{R}}^{d+1})$ is an $n$-cube of
$\mathcal{F}(\overline{\mathcal{R}}^{d})$ together with a perturbation $B: [0,1]^n\times [0,1]\rightarrow \mathcal{F}_0(\overline{\mathcal{R}}^{d+1})$.
Given $\mathfrak{o}\in\mathcal{F}_n(\overline{\mathcal{R}}^{d+1})$, one can first set $B$ to be constant along
the second last $[0,1]$-coordinate and post-compose with the deformation retract $\tilde{\mathcal{F}}(\overline{\mathcal{R}}^{d+1})\rightarrow \tilde{\mathcal{F}}^{reg}(\overline{\mathcal{R}}^{d+1})$.
Apply this construction cube-wise to $\tilde{s}'$, one obtain the desired lift $\tilde{s}$.
Alternatively, since each $\overline{\mathcal{R}}^{d+1}$ has finitely many cells, each of which is subdivided into finitely many cubes,
and there are countably many such spaces ($d\geq 2$) in total, as a standard consequence of Sard's theorem one can choose the perturbations $B$ to make each cube regular. \qed \par\indent
6.2. \textbf{Proof of Theorem 6.1. 1)}. 
Pulling back $A_L=CF^*(L,L)$, viewed as an $\mathcal{A}_{\infty}^{oepr,dg}$-algebra, along $\tilde{s}: \mathcal{A}^{dg}_{\infty}\rightarrow \mathcal{A}_{\infty}^{oper,dg}$ of Lemma 6.2, we recover the classical construction of the Floer
$\mathcal{A}_{\infty}$-algebra $\tilde{s}^*A_L$ associated to $L\in \mathrm{Fuk}(M)_{\lambda}$.
Let $\tilde{s}_{\overrightarrow{\Lambda}}$ denote the induced homotopy equivalence of dg categories
\begin{equation}
\tilde{s}_{\overrightarrow{\Lambda}}: \overrightarrow{\Lambda}\rtimes \mathcal{A}^{dg}_{\infty}\rightarrow \overrightarrow{\Lambda}\rtimes \mathcal{A}_{\infty}^{oper,dg}.
\end{equation}
induced by $\tilde{s}$. There is a homotopy commutative diagram of dg categories
\begin{equation}
\begin{tikzcd}[row sep=1.2cm, column sep=0.8cm]
\overrightarrow{\Lambda}\rtimes \mathcal{A}^{dg}_{\infty}\arrow[r,"\tilde{s}_{\overrightarrow{\Lambda}}"]\arrow[d,"P_{\overrightarrow{\Lambda}}"]& \overrightarrow{\Lambda}\rtimes \mathcal{A}_{\infty}^{oper,dg}\arrow[dl,"\tilde{P}_{\overrightarrow{\Lambda}}"]\\
\overrightarrow{\Lambda}_k&
\end{tikzcd}.
\end{equation}
Therefore, there is a quasi-isomorphism
\begin{equation}
CC^{S^1,oper}(A_L):=\Big(\int_{\mathbb{T}}(\iota_{\Lambda})_!(\tilde{P}^{-1}_{\overrightarrow{\Lambda}})^*A^{\sharp}_L\Big)^{h\mathbb{T}}\simeq \Big(\int_{\mathbb{T}} (\iota_{\Lambda})_!(P^{-1}_{\overrightarrow{\Lambda}})^*\tilde{s}_{\overrightarrow{\Lambda}}^*A_L^{\sharp}\Big)^{h\mathbb{T}}.
\end{equation}
Composing with the quasi-isomorphism of Proposition 3.11, we obtain a quasi-isomorphism
\begin{equation}
(\Xi^{S^1})^{-1}: CC^{S^1}(A_L)=CC^{S^1}(\tilde{s}_{\overrightarrow{\Lambda}}^*A^{\sharp}_L)\xrightarrow[\simeq]{(3.21)} \Big(\int_{\mathbb{T}}(\iota_{\Lambda})_!(P^{-1}_{\overrightarrow{\Lambda}})^*\tilde{s}_{\overrightarrow{\Lambda}}^*A^{\sharp}_L\Big)^{h\mathbb{T}}\xrightarrow[\simeq]{(6.21)} \Big(\int_{\mathbb{T}} (\iota_{\Lambda})_!(\tilde{P}_{\overrightarrow{\Lambda}}^{-1})^*A^{\sharp}_L\Big)^{h\mathbb{T}}=CC^{S^1,oper}(A_L).
\end{equation}
We take $\Xi^{S^1}$ in Theorem 6.1. 1) to be a choice of homotopy inverse.
Pulling back $OC^{oper}$ (cf. (4.24)), viewed as a morphism of dg functors in
$\mathrm{Fun}^{dg}(\overrightarrow{\Lambda}^{op}\rtimes\mathcal{A}_{\infty}^{oper,dg},\mathrm{Mod}_k)$, along $\tilde{s}_{\overrightarrow{\Lambda}}$ gives rise to a map of $\mathcal{A}_{\infty}$-cocyclic $k$-modules
\begin{equation}
(\tilde{s}^{op}_{\overrightarrow{\Lambda}})^*\overline{\mathfrak{R}}^1\xrightarrow{(\tilde{s}^{op}_{\overrightarrow{\Lambda}})^*OC^{oper}} \mathrm{Map}((\tilde{s}^{op}_{\overrightarrow{\Lambda}})^*A^{\sharp}_L, QH).
\end{equation}
Apply $\Big(\int^{\mathbb{T}}(\iota^{op}_{\Lambda})_*((P^{op}_{\overrightarrow{\Lambda}})^{-1})^*(-)\Big)_{h\mathbb{T}}$ to both sides of the equation (6.23), and undergo the same process as in (4.27)-(4.30), we obtain a map of
chain complex
\begin{equation}
\Big(\int^{\mathbb{T}}(\iota^{op}_{\Lambda})_* (P^{op}_{\overrightarrow{\Lambda}})^{-1})^*(\tilde{s}^{op}_{\overrightarrow{\Lambda}})^*\overline{\mathfrak{R}}^1\Big)_{h\mathbb{T}}\rightarrow \mathrm{Map}(\Big(\int_{\mathbb{T}}(\iota_{\Lambda})_! (P^{-1}_{\overrightarrow{\Lambda}})^*\tilde{s}_{\overrightarrow{\Lambda}}^*A^{\sharp}_L\Big)^{h\mathbb{T}},QH).
\end{equation}
By commutativity of (6.20), we have a commutative diagram of chain complexes
\begin{equation}
\begin{tikzcd}[row sep=1.2cm, column sep=0.8cm]
\Big(\int^{\mathbb{T}}(\iota^{op}_{\Lambda})_* ((\tilde{P}^{op}_{\overrightarrow{\Lambda}})^{-1})^*\overline{\mathfrak{R}}^1\Big)_{h\mathbb{T}}\arrow[rr,"{(4.30)}"]\arrow[d,"{\simeq}"]& & \mathrm{Map}(\Big(\int_{\mathbb{T}}(\iota_{\Lambda})_! (\tilde{P}^{-1}_{\overrightarrow{\Lambda}})^*A^{\sharp}_L\Big)^{h\mathbb{T}},QH)\arrow[d,"{\simeq}"]\\
\Big(\int^{\mathbb{T}}(\iota^{op}_{\Lambda})_* ((P^{op}_{\overrightarrow{\Lambda}})^{-1})^*(\tilde{s}^{op}_{\overrightarrow{\Lambda}})^*\overline{\mathfrak{R}}^1\Big)_{h\mathbb{T}}\arrow[rr,"{(6.24)}"]& & \mathrm{Map}(\Big(\int_{\mathbb{T}}(\iota_{\Lambda})_! (P^{-1}_{\overrightarrow{\Lambda}})^*\tilde{s}_{\overrightarrow{\Lambda}}^*A^{\sharp}_L\Big)^{h\mathbb{T}},QH)
\end{tikzcd},
\end{equation}
where the vertical arrows are quasi-isomorphisms. There is also a chain of isomorphisms
\begin{equation}
H^*(\Big(\int^{\mathbb{T}}(\iota^{op}_{\Lambda})_* ((\tilde{P}^{op}_{\overrightarrow{\Lambda}})^{-1})^*\overline{\mathfrak{R}}^1\Big)_{h\mathbb{T}})\xrightarrow[\cong]{(*)} H^*(\Big(\int^{\mathbb{T}}(\iota^{op}_{\Lambda})_* ((P^{op}_{\overrightarrow{\Lambda}})^{-1})^*(\tilde{s}^{op}_{\overrightarrow{\Lambda}})^*\overline{\mathfrak{R}}^1\Big)_{h\mathbb{T}})\xrightarrow[\cong]{(**)} k[\tilde{t}],
\end{equation}
where (*) follows from commutativity of (6.20) and (**) follows from the contractibility of $\overline{\mathfrak{R}}^1$ together with Proposition 3.11.
Let
\begin{equation}
\tilde{OC}^{S^1,oper}_i:\Big(\int_{\mathbb{T}}(\iota_{\Lambda})_! (P^{-1}_{\overrightarrow{\Lambda}})^*\tilde{s}_{\overrightarrow{\Lambda}}^*A^{\sharp}_L\Big)^{h\mathbb{T}}\rightarrow QH(M)
\end{equation}
denote the chain map given by the image
of $\tilde{t}^i$ under (the cohomological level map of) (6.26) and then (6.24). Define
\begin{equation}
\tilde{OC}^{S^1,oper}:=\sum_{i\geq 0}\tilde{OC}^{S^1,oper}_it^i: \Big(\int_{\mathbb{T}}(\iota_{\Lambda})_! (P^{-1}_{\overrightarrow{\Lambda}})^*\tilde{s}_{\overrightarrow{\Lambda}}^*A^{\sharp}_L\Big)^{h\mathbb{T}}\rightarrow QH(M)[[t]].
\end{equation}
By commutativity of (6.25), to prove Theorem 6.1. 1) it thus suffices to show the following diagram commutes
\begin{equation}
\begin{tikzcd}[row sep=1.2cm, column sep=0.8cm]
\Big(\int_{\mathbb{T}}(\iota_{\Lambda})_! (P^{-1}_{\overrightarrow{\Lambda}})^*\tilde{s}_{\overrightarrow{\Lambda}}^*A^{\sharp}_L\Big)^{h\mathbb{T}}\arrow[rr,"{\tilde{OC}^{S^1,oper}}"]& &QH(M)[[t]]\\
CC^{S^1}(A)\arrow[urr,"OC^{S^1}"]\arrow[u,"\simeq"]& &
\end{tikzcd},
\end{equation}
where the vertical quasi-isomorphism is given by Proposition 3.11.\par\indent
The key ingredient to proving commutativity of (6.29) is the following Lemma, which allows one to give explicit formula for the generators
$\tilde{t}^i\in k[\tilde{t}]\cong H^*(\Big(\int^{\mathbb{T}}(\iota^{op}_{\Lambda})_* ((P^{op}_{\overrightarrow{\Lambda}})^{-1})^*(\tilde{s}^{op}_{\overrightarrow{\Lambda}})^*\overline{\mathfrak{R}}^1\Big)_{h\mathbb{T}})$, where we compute the latter using the explicit
positive cocyclic complex $CC^{S^1,\vee}((\tilde{s}^{op}_{\overrightarrow{\Lambda}})^*\overline{\mathfrak{R}}^1)$, cf. Proposition 3.24.
\begin{lemma}
There exists\par\indent
1) A choice of Floer data for the cyclic open-closed map (cf. section 2.2).\par\indent
2) For each cell
\begin{equation}
K\in\{{}_{k}\overline{\check{\mathcal{R}}^1}_{d+1},{}_{k}\overline{\check{\mathcal{R}}^1}_{d+1}^{S^1}, {}_k^{i,i+1}\overline{\check{\mathcal{R}}^1}_{d+1}, {}_{k}\overline{\hat{\mathcal{R}}^1}_{d+1},{}_{k}\overline{\hat{\mathcal{R}}^1}_{d+1}^{S^1},{}_k^{i,i+1}\overline{\hat{\mathcal{R}}^1}_{d+1}\}_{d,k\geq 0, 1\leq i\leq k-1},
\end{equation}
one can assign a cubical subdivision $s_K([K])$, compatibly with product and boundary structures. In particular,
they fit together into chain quasi-isomorphisms $s_K: C_*^{cell}(K)\rightarrow C_*^{\Box}(K)$. \par\indent
3) For each $K$ in (6.30), one can find chain maps $\tilde{s}_K: C_*^{cell}(K)\rightarrow C_*(\tilde{\mathcal{F}}^{reg}(\overline{\mathcal{R}}^1_{d+1}))$
compatible with product structures at the boundary such that
\begin{enumerate}[label=3\alph*)]
    \item  These fit into commutative diagrams
\begin{equation}
\begin{tikzcd}[row sep=1.2cm, column sep=0.8cm]
C_*^{cell}(K)\arrow[r,"{\tilde{s}_K}"]\arrow[d,"{s_K}"]& C_*(\tilde{\mathcal{F}}^{reg}(\overline{\mathcal{R}}^1_{d+1}))\arrow[d,"{\pi}"]\\
C_*^{\Box}(K)\arrow[r,"{\mathrm{forget}}"]& C_*^{\Box}(\overline{\mathcal{R}}^1_{d+1})
\end{tikzcd},
\end{equation}
where `forget' denotes the map that forgets the additional interior marked points $p_1,\cdots,p_k$, as well as the boundary marked point $z_f$ if $K$ is of hat($\,\hat{}\,$)-type.
   \item $\tilde{s}_K([K])$ has the same degree as the dimension of $K$, and
\begin{equation}
OC^{oper}(\tilde{s}_K([K])): A_L^{\otimes d+1}\rightarrow QH,
\end{equation}
agrees with the map induced by counting rigid solutions to the parametrized moduli problem associated to the universal choice of Floer data in 1) restricted to $K$.
   \item $\tilde{s}(K)=0$ if $K$ is of type $\{{}_k^{i,i+1}\overline{\check{\mathcal{R}}^1}_{d+1},{}_k^{i,i+1}\overline{\hat{\mathcal{R}}^1}_{d+1},{}_k\overline{\hat{\mathcal{R}}^1}^{S^1}_{d+1}\}$.
\end{enumerate}
\end{lemma}
When the context is clear, we will denote simply $s$ and $\tilde{s}$ and omit the subscript $K$.\par\indent
\emph{Proof of Theorem 6.1 1) given Lemma 6.3}. Recall that it remains to show commutativity of (6.29). By definition,
the $i$-th coefficient of $\tilde{OC}^{S^1,oper}$ is image of $\tilde{t}^i\in k[\tilde{t}]\cong H^*(\Big(\int^{\mathbb{T}}(\iota^{op}_{\Lambda})_* ((P^{op}_{\overrightarrow{\Lambda}})^{-1})^*(\tilde{s}^{op}_{\overrightarrow{\Lambda}})^*\overline{\mathfrak{R}}^1\Big)_{h\mathbb{T}})$
under (6.24). By Proposition 3.24, $H^*(\Big(\int^{\mathbb{T}}(\iota^{op}_{\Lambda})_* ((P^{op}_{\overrightarrow{\Lambda}})^{-1})^*(\tilde{s}^{op}_{\overrightarrow{\Lambda}})^*\overline{\mathfrak{R}}^1\Big)_{h\mathbb{T}})$ can be computed as the
cohomology of the chain complex
\begin{equation}
CC^{S^1,\vee}((\tilde{s}^{op}_{\overrightarrow{\Lambda}})^*\overline{\mathfrak{R}}^1)=
Tot^{\oplus}\Big(\cdots\xrightarrow{\tau-1}(CC^{\vee}((\tilde{s}^{op}_{\overrightarrow{\Lambda}})^*\overline{\mathfrak{R}}^1),b')\xrightarrow{N}(CC^{\vee}((\tilde{s}^{op}_{\overrightarrow{\Lambda}})^*\overline{\mathfrak{R}}^1),b)\xrightarrow{\tau-1}(CC^{\vee}((\tilde{s}^{op}_{\overrightarrow{\Lambda}})^*\overline{\mathfrak{R}}^1),b') \rightarrow 0\Big)
\end{equation}
or equivalently,
\begin{equation}
CC^{S^1,\vee}((\tilde{s}^{op}_{\overrightarrow{\Lambda}})^*\overline{\mathfrak{R}}^1)=CC^{\vee}((\tilde{s}^{op}_{\overrightarrow{\Lambda}})^*\overline{\mathfrak{R}}^1)[\tilde{t},e^+]=\prod_{d\geq0}{C_{-*}(\tilde{\mathcal{F}}^{reg}(\overline{\mathcal{R}}^1_{d+1}))[-d]}[\tilde{t},e^+],
\end{equation}
where $\tilde{t}$ is a formal variable of degree $-2$, $e^+$ has degree $-1$ and satisfy $(e^+)^2=0$, equipped with the differential
\begin{equation}
\begin{cases}
x\tilde{t}^k\mapsto b'x\tilde{t}^k+Nxe^+\tilde{t}^{k-1}\\
xe^+\tilde{t}^k\mapsto bx e^+\tilde{t}^k+(\tau-1)x\tilde{t}^k.
\end{cases}
\end{equation}
Lemma 6.3. 2) and 3c) together imply the relations (compare the boundary decompositions (2.22)-(2.24), (2.26)-(2.29))
\begin{equation}
\partial \tilde{s}({}_{k}\overline{\check{\mathcal{R}}^1}_{d+1})=\sum_{m,i} \pm \tilde{s}(\overline{\mathcal{R}}^{m+1})\times_i s({}_{k}\overline{\check{\mathcal{R}}^1}_{d-m+1})+(1+\tau+\cdots+\tau^d)\tilde{s}({}_{k-1}\overline{\hat{\mathcal{R}}^1}_{d+1})
\end{equation}
and
\begin{equation}
\partial \tilde{s}({}_{k}\overline{\hat{\mathcal{R}}^1}_{d+1})=\sum_{m,i} \pm \tilde{s}(\overline{\mathcal{R}}^{m+1})\times_i \tilde{s}({}_{k}\overline{\hat{\mathcal{R}}^1}_{d-m+1})+(\tau-1)\tilde{s}({}_k\overline{\check{\mathcal{R}}^1}_{d+1}).
\end{equation}
Therefore, the generator corresponding to $\tilde{t}^i\in k[\tilde{t}]\cong H^*(CC^{S^1,\vee}((\tilde{s}^{op}_{\overrightarrow{\Lambda}})^*\overline{\mathfrak{R}}^1))$ has a chain representative in $CC^{S^1,\vee}((\tilde{s}^{op}_{\overrightarrow{\Lambda}})^*\overline{\mathfrak{R}}^1)$
given by
\begin{equation}
\Big(\prod_d\tilde{s}({}_i\overline{\hat{\mathcal{R}}^1}_{d+1})+ \prod_d\tilde{s}({}_i\overline{\check{\mathcal{R}}^1}_{d+1})e^+\Big)+\cdots+\Big(\prod_d\tilde{s}({}_0\overline{\hat{\mathcal{R}}^1}_{d+1})\tilde{t}^i+ \prod_d\tilde{s}({}_0\overline{\check{\mathcal{R}}^1}_{d+1})\tilde{t}^ie^+\Big).
\end{equation}
As a result, the $t^i$-coefficient of the composition
\begin{equation}
CC^{S^1}(A_L)\simeq \Big(\int_{\mathbb{T}}(\iota_{\Lambda})_!(P^{-1}_{\overrightarrow{\Lambda}})^*\tilde{s}_{\overrightarrow{\Lambda}}^*A^{\sharp}_L\Big)^{h\mathbb{T}}\xrightarrow{\tilde{OC}^{S^1,oper}}QH(M)[[t]]
\end{equation}
is homotopic to the chain map that sends $(x_0+x_0'e^+)+(x_1u+x_1'te^+)+\cdots\in CC^{S^1}(A_L)$, where $x_i\in A_L^{\otimes d_i+1}, x_i'\in A_L^{\otimes d_i'+1}\subset CC(A_L)$, to
\begin{equation}
\sum_{0\leq j\leq i} \big(OC^{oper}(\tilde{s}([{}_j\overline{\check{\mathcal{R}}^1}_{d_j+1}]))(x_{i-j})+OC^{oper}(\tilde{s}([{}_j\overline{\hat{\mathcal{R}}^1}_{d_j'+1}]))(x_{i-j}')\big)\in QH(M).
\end{equation}
Lemma 6.3. 3b) implies that this is exactly the $t^i$-coefficient of Ganatra's cyclic open-closed map, which proves that (6.29) commutes. \qed\par\indent
\emph{Proof of Lemma 6.3}. The existence of a universal choice of Floer data with the specified requirements follows from \cite[Proposition. 10]{Gan2}.
Below we reproduce the proof and in the meantime construct the desired cubical subdivisions $s([K])\in C_*^{\Box}(K)$ and a lift $\tilde{s}': C_*^{cell}(K)\rightarrow C_*(\mathcal{F}(\overline{\mathcal{R}}^1_{d+1})$
of the composition $\mathrm{forget}\circ s$, cf. (6.29). Then, the lifts $\tilde{s}(K)$ to $C_*(\tilde{\mathcal{F}}^{reg}(\overline{\mathcal{R}}^1_{d+1}))$ is obtained by a standard perturbation argument as in Lemma 6.2.  \par\indent
Recall from Lemma 6.2 that we have fixed a universal choice of Floer data for $\overline{\mathcal{R}}^{\bullet+1}$, together with cubical subdivisons $s([K])$ of $[K]$, and chains $s([K])\in C_*(\mathcal{F}(\overline{\mathcal{R}}^{\bullet+1}))$, for $K=\overline{\mathcal{R}}^{d+1}, d\geq 2$.
From now on, we fix a universal choice of Floer data for $K=\overline{\mathcal{R}}^{\bullet+1}$ as well as $s_K, \tilde{s}_K',\tilde{s}_K$ as in Lemma 6.2. \par\indent
As a preliminary step, we observe there exists a universal choice of Floer data for the spaces of disks with forgotten marked points $\overline{\mathcal{R}}^{d+1,f_i}, d\geq 2$
(topologically this is a copy of $\overline{\mathcal{R}}^{d+1}$, but with the $i$-th marked point viewed as `forgotten'), satisfying the usual consistency conditions at the boundary plus the following:
\begin{itemize}
    \item  the Floer data on the unique element of $\overline{\mathcal{R}}^{3,f_i}$ is translation invariant after forgetting the $i$-th point;
    \item for $d>2$, the Floer data on $\overline{\mathcal{R}}^{d+1,f_i}$ is pulled back from the forgetful map $\overline{\mathcal{R}}^{d+1,f_i}\rightarrow \overline{\mathcal{R}}^{d}$.
\end{itemize}
The two conditions above uniquely determines the choice of Floer data.\par\indent
On the other hand, one can identify
\begin{equation}
\overline{\mathcal{R}}^{d+1,f_i}\cong \overline{\mathcal{R}}^{d}\times [0,1]
\end{equation}
as topological spaces (the identification is smooth in the open part but not over the corners).
Thus the product of $s([\overline{\mathcal{R}}^{d}])$ with the interval $[0,1]$ forms a cubical subdivision of $\overline{\mathcal{R}}^{d+1,f_i}$, smooth in codimension $0$, and we can that to be the definition of $s([\overline{\mathcal{R}}^{d+1,f_i}])$.
Next, we define
\begin{equation}
\tilde{s}'([\overline{\mathcal{R}}^{d+1,f_i}]):=0\in C_*(\mathcal{F}(\overline{R}^{d-1})).
\end{equation}
Morally, the induced cubical subdivision from (6.41) is equipped a degenerate Floer data (i.e. constant along the second $[0,1]$-coordinate),
hence becomes $0$ in the normalized cubical chains $C_*(\mathcal{F}(\overline{\mathcal{R}}^{d}))$.
This agrees with the fact that the Floer operation associated to $\overline{\mathcal{R}}^{d+1,f_i}, d>2$ is zero, cf \cite[Proposition. 15]{Gan2}.
Here and in the rest of the proof, we only specify the underlying map $b$ of a cube $\mathfrak{o}\in\mathcal{F}_*(\overline{\mathcal{R}}^{d+1,f_i})$, with the associated $(b_f,g_f)_f$ defined in a similar fashion as in Lemma 6.2.\par\indent
Next, we inductively choose Floer data for ${}_{k}\overline{\check{\mathcal{R}}^1}_{d+1}, {}_{k}\overline{\hat{\mathcal{R}}^1}_{d+1}$ together with cubical subdivisons
$s([K])$'s and lifts $\tilde{s}'([K])$'s satisfying the desired requirements. The choice of Floer data for $(d=0,k=0)$ can be arbitrary. Suppose the choice is made for all
\begin{equation}
{}_{k}\overline{\check{\mathcal{R}}^1}_{d+1}, {}_{k}\overline{\hat{\mathcal{R}}^1}_{d+1}, (d,k)<(d_0,k_0).
\end{equation}
Condition Lemma 6.3. 1) uniquely determines a choice of Floer data for $K$ a codimension one boundary component of type (2.22); we define $s([K])$ to be
\begin{equation}
s([\overline{\mathcal{R}}^{s+1}])\times_i s([{}_k\overline{\check{\mathcal{R}}^1}_{d_0-s+1}]),
\end{equation}
and $\tilde{s}'([K])$ to be the product of the underlying cubes equipped with the product of Floer data/gluing atlases.
Next, we consider the case $K={}_{k_0-1}(\overline{\mathcal{R}}^1)^{S^1}_{d_0+1}$, type (2.23). For the locus ${}_{k_0-1}(\overline{\mathcal{R}}^1)^{S^1_{d_0,1}}_{d_0+1}$,
where $p_{k_0}$ points between $z_{d_0}$ and $z_1$, we define
\begin{equation}
s([{}_{k_0-1}(\overline{\mathcal{R}}^1)^{S^1_{d_0,1}}_{d_0+1}]):=h_*s([{}_{k_0-1}\overline{\hat{\mathcal{R}}^1}_{d_0}]),
\end{equation}
where $h$ is the real blow down map
\begin{equation}
{}_{k_0-1}\overline{(\mathcal{R}^1)}^{S^1_{d_0,1}}_{d_0+1}\leftarrow {}_{k_0-1}\overline{\hat{\mathcal{R}}}^1_{d_0}
\end{equation}
which is a diffeomorphism on the open part. The Floer data on ${}_{k_0-1}(\mathcal{R}^1)^{S^1_{d_0,1}}_{d_0+1}$ is pulled back along $h$ on the open part.
When $p_{k_0}$ approaches either $z_{d_0}$ or $z_1$ (note the compactifications of the left and right hand side of (6.46) differs),
the Floer data is pulled back from ${}_{k_0-1}\check{\overline{\mathcal{R}}}^1_{d_0+1}$. Accordingly, we set
\begin{equation}
\tilde{s}'([{}_{k_0-1}(\overline{\mathcal{R}}^1)^{S^1_{d_0,1}}_{d_0+1}]):=\tilde{s}'([{}_{k_0-1}\overline{\hat{\mathcal{R}}^1}_{d_0}])
\end{equation}
since the (induced map on cubical subdivision of) blowdown $h$ only collapse degenerate cubes, which are $0$ in the symmetric normalized cubical complex.
To satisfy equivariance, we define the Floer data on ${}_{k_0-1}(\overline{\mathcal{R}}^1)^{S^1_{i,i+1}}_{d_0+1}$, the locus where $p_{k_0}$ points between $z_i$ and $z_{i+1}$,
to be pulled back from ${}_{k_0-1}(\overline{\mathcal{R}}^1)^{S^1_{d_0,1}}_{d_0+1}$ via cyclic permutation of boundary marked points, and accordingly set
\begin{equation}
\tilde{s}'([{}_{k_0-1}(\overline{\mathcal{R}}^1)^{S^1_{i,i+1}}_{d_0+1}]):=\tau^i \tilde{s}'([{}_{k_0-1}(\overline{\mathcal{R}}^1)^{S^1_{d_0,1}}_{d_0+1}]).
\end{equation}
This gives a choice of Floer data on ${}_{k_0-1}(\overline{\mathcal{R}}^1)^{S^1}_{d_0+1}$, and we accordingly define
\begin{equation}
\tilde{s}'([{}_{k_0-1}(\overline{\mathcal{R}}^1)^{S^1}_{d_0+1}]):=\sum_i \tilde{s}'([{}_{k_0-1}(\overline{\mathcal{R}}^1)^{S^1_{i,i+1}}_{d_0+1}])=(1+\tau+\cdots+\tau^{d_0-1}) \tilde{s}'([{}_{k_0-1}(\overline{R}^1)^{S^1_{d_0,1}}_{d_0+1}]).
\end{equation}
For $K={}_{k_0}^{i,i+1}\overline{\check{\mathcal{R}}^1}_{d_0+1}$ of type (2.24), we pull back the Floer data via the map
\begin{equation}
\pi_i: {}_{k_0}^{i,i+1}\overline{\check{\mathcal{R}}^1}_{d_0+1}\rightarrow {}_{k_0-1}\overline{\check{\mathcal{R}}^1}_{d_0}
\end{equation}
that forgets $p_{i+1}$. We define
\begin{equation}
s([{}_{k_0}^{i,i+1}\overline{\check{\mathcal{R}}^1}_{d_0}]):=s([{}_{k_0-1}\overline{\check{\mathcal{R}}^1}_{d_0+1}])\times [0,1],
\end{equation}
where the right hand side is viewed as a cubical subdivision of ${}_{k_0}^{i,i+1}\overline{\check{\mathcal{R}}^1}_{d_0+1}\cong {}_{k_0-1}\overline{\check{\mathcal{R}}^1}_{d_0+1}\times S^1$ via the map $[0,1]\rightarrow S^1$
identifying endpoints on the second coordinate. Having in mind that this is a cube equipped with degenerate Floer data, we set
\begin{equation}
\tilde{s}'([{}_{k_0}^{i,i+1}\overline{\check{\mathcal{R}}^1}_{d_0+1}]):=0.
\end{equation}\par\indent
Having inductively defined Floer data and associated cubical chains $s([K])$ for $K$ a codimension $1$ boundary component,
a similar argument involving tubular neighborhoods and gluing charts as in the proof of Lemma 6.2 allows us to extend the Floer data to
all of ${}_{k_0}\overline{\check{\mathcal{R}}^1}_{d_0+1}$, to construct $s([{}_{k_0}\overline{\check{\mathcal{R}}^1}_{d_0+1}])$
and $\tilde{s}'([{}_{k_0}\overline{\check{\mathcal{R}}^1}_{d_0+1}])$. \par\indent
By the inductive nature of the construction, $\partial s([K])=s([\partial K])$ (similarly for $\tilde{s}'$), i.e. $s$ is a chain map. Condition 3a) is clear from construction.
When $K$ is not of type (2.24), we equipped $s(K)$ with the Floer data restricted from the universal choice of Floer data, which ensures condition 3b)
for those $K$. When $K$ is of type (2.24), the universal choice of Floer data restricted to $K$ is degenerate, i.e. pulled
back from a lower dimensional strata, and we have accordingly set $\tilde{s}'(K)=0$, which ensures 3b) and 3c) for those cells. \par\indent
We can continue the induction to choose a Floer data for ${}_{k_0}\overline{\hat{\mathcal{R}}^1}_{d_0+1}$ and construct $s([{}_{k_0}\overline{\hat{\mathcal{R}}^1}_{d_0+1}])$ and $\tilde{s}'([{}_{k_0}\overline{\hat{\mathcal{R}}^1}_{d_0+1}])$.
Note that in its codimension 1 boundary components (2.26)-(2.29), the cells equipped with degenerate Floer data are (products with) one of $\{\overline{\mathcal{R}}^{d+1,f_i},{}_{k-1}\overline{\hat{\mathcal{R}}^1}_{d+1}^{S^1},{}_k^{i,i+1}\overline{\hat{\mathcal{R}}^1}_{d}\}$.
We accordingly set $s$ to be zero on these cells.
The rest is completely analogous to the case for ${}_{k_0}\overline{\check{\mathcal{R}}^1}_{d_0+1}$.  \qed\par\indent
6.3. \textbf{Proof of Theorem 6.1. 2)}. Consider the parameter spaces of disks $\overline{\mathcal{R}}^1_{k_1,\cdots,k_p}$ from section 4.1. 3). Similar to before,
one can form the chain complexes $C_*(\tilde{\mathcal{F}}^{reg}(\overline{\mathcal{R}}^1_{k_1,\cdots,k_p}))$ using operadic Floer theory. Analogous to Definition 4.8, the assignment
\begin{equation}
[k_1,\cdots,k_p]\mapsto C_{-*}(\tilde{\mathcal{F}}^{reg}(\overline{\mathcal{R}}^1_{k_1,\cdots,k_p}))
\end{equation}
defines a dg functor from $\overrightarrow{{}_p\Lambda}^{op}\rtimes \mathcal{A}_{\infty}^{oper,dg}$ to $\mathrm{Mod}_k$. We denote this dg functor by ${}_p\overline{\mathfrak{R}}^1$.
There are inclusions of spaces $\overline{\mathcal{R}}^1_{k_1,k_2,\cdots,k_p}\hookrightarrow \overline{\mathcal{R}}^1_{k_1+\cdots+k_p+p-1}$, which are homotopy equivalences (as both are contractible),
inducing a homotopy equivalence of dg $\overrightarrow{\Lambda}^{op}_p\rtimes \mathcal{A}_{\infty}^{oper,dg}$-modules
\begin{equation}
{}_p\overline{\mathfrak{R}}^1\rightarrow (\overrightarrow{j^{op}_{\mathcal{A}_{\infty}^{oper,dg}}})^*\overline{\mathfrak{R}}^1.
\end{equation}
Pulling back (6.54) along the dg functor $\tilde{s}^{op}_{\overrightarrow{{}_p\Lambda}}: \overrightarrow{\Lambda}^{op}_p\rtimes \mathcal{A}_{\infty}^{dg}\rightarrow \overrightarrow{\Lambda}^{op}_p\rtimes \mathcal{A}_{\infty}^{oper,dg}$ induced by $\tilde{s}$ (cf. Lemma 6.2),
we obtain a homotopy equivalence of dg $\overrightarrow{\Lambda}^{op}_p\rtimes \mathcal{A}_{\infty}^{dg}$-modules
\begin{equation}
(\tilde{s}^{op}_{\overrightarrow{{}_p\Lambda}})^*{}_p\overline{\mathfrak{R}}^1\rightarrow (\tilde{s}^{op}_{\overrightarrow{{}_p\Lambda}})^*(\overrightarrow{j^{op}_{\mathcal{A}_{\infty}^{oper,dg}}})^*\overline{\mathfrak{R}}^1.
\end{equation}
We apply $\Big(\int^{\mathbb{Z}/p}(\iota^{op}_{{}_p\Lambda})_*((P^{op}_{\overrightarrow{{}_p\Lambda}})^{-1})^*(-)\Big)_{h\mathbb{Z}/p}$ to
\begin{equation}
(\tilde{s}^{op}_{\overrightarrow{{}_p\Lambda}})^*{}_p\overline{\mathfrak{R}}^1\xrightarrow{(6.55)} (\tilde{s}^{op}_{\overrightarrow{{}_p\Lambda}})^*(\overrightarrow{j^{op}_{\mathcal{A}_{\infty}^{oper,dg}}})^*\overline{\mathfrak{R}}^1 \xrightarrow{(\tilde{s}^{op}_{\overrightarrow{{}_p\Lambda}})^*(\overrightarrow{j^{op}_{\mathcal{A}_{\infty}^{oper,dg}}})^*OC^{oper}} (\tilde{s}^{op}_{\overrightarrow{{}_p\Lambda}})^*(\overrightarrow{j^{op}_{\mathcal{A}_{\infty}^{oper,dg}}})^*\mathrm{Map}(A^{\sharp}_L,QH),
\end{equation}
and obtain a chain map
\begin{equation}
\Big(\int^{\mathbb{Z}/p}(\iota_{{}_p\Lambda}^{op})_*((P^{op}_{\overrightarrow{{}_p\Lambda}})^{-1})^*(\tilde{s}^{op}_{\overrightarrow{{}_p\Lambda}})^*{}_p\overline{\mathfrak{R}}^1\Big)_{h\mathbb{Z}/p}\rightarrow \mathrm{Map}(\Big(\int^{\mathbb{Z}/p}(\iota_{{}_p\Lambda})_!((P^{op}_{\overrightarrow{{}_p\Lambda}})^{-1})^*(\tilde{s}^{op}_{\overrightarrow{{}_p\Lambda}})^*(\overrightarrow{j^{op}_{\mathcal{A}_{\infty}^{oper,dg}}})^*A^{\sharp}_L\Big)^{h\mathbb{Z}/p},QH).
\end{equation}
The left hand side has cohomology canonically isomorphic to $k[\tilde{t},\tilde{\theta}]$, and denote by $\tilde{OC}^{\mathbb{Z}/p,oper}_{2i}$ (resp. $\tilde{OC}^{\mathbb{Z}/p,oper}_{2i+1}$)
the chain map $\Big(\int^{\mathbb{Z}/p}(\iota_{{}_p\Lambda})_!(P^{-1}_{\overrightarrow{{}_p\Lambda}})^*\tilde{s}^*_{\overrightarrow{{}_p\Lambda}}(\overrightarrow{j_{\mathcal{A}_{\infty}^{oper,dg}}})^*A^{\sharp}_L\Big)^{h\mathbb{Z}/p}\rightarrow QH$ given by the (cohomological level) image of $\tilde{t}^k$ (resp. $\tilde{t}^k\tilde{\theta}$)
under (6.57). \par\indent
Define the chain map $\tilde{OC}^{\mathbb{Z}/p,oper}$ to be
\begin{equation}
\sum_{i\geq0}(\tilde{OC}^{\mathbb{Z}/p,oper}_{2i}+\tilde{OC}^{\mathbb{Z}/p,oper}_{2i+1}\theta)t^i: \Big(\int^{\mathbb{Z}/p}(\iota_{{}_p\Lambda})_!(P^{-1}_{\overrightarrow{{}_p\Lambda}})^*\tilde{s}^*_{\overrightarrow{{}_p\Lambda}}(\overrightarrow{j_{\mathcal{A}_{\infty}^{oper,dg}}})^*A^{\sharp}_L\Big)^{h\mathbb{Z}/p}\rightarrow QH(M)[[t,\theta]].
\end{equation}
There are canonical equivalences
\begin{equation}
\Big(\int^{\mathbb{Z}/p}(\iota_{{}_p\Lambda})_!(P^{-1}_{\overrightarrow{{}_p\Lambda}})^*(\overrightarrow{j_{\mathcal{A}_{\infty}^{dg}}})^*\tilde{s}^*_{\overrightarrow{\Lambda}}A^{\sharp}_L\Big)^{h\mathbb{Z}/p}\simeq\Big(\int^{\mathbb{Z}/p}(\iota_{{}_p\Lambda})_!(P_{\overrightarrow{{}_p\Lambda}}^{-1})^*\tilde{s}^*_{\overrightarrow{{}_p\Lambda}}(\overrightarrow{j_{\mathcal{A}_{\infty}^{oper,dg}}})^*A^{\sharp}_L\Big)^{h\mathbb{Z}/p}\simeq CC^{\mathbb{Z}/p,oper}((\overrightarrow{j_{\mathcal{A}_{\infty}^{oper,dg}}})^*A^{\sharp}_L),
\end{equation}
where the first equivalence follows from $\tilde{s}_{\overrightarrow{\Lambda}}\circ\overrightarrow{j_{\mathcal{A}_{\infty}^{dg}}}=\overrightarrow{j_{\mathcal{A}_{\infty}^{oper,dg}}}\circ\tilde{s}_{\overrightarrow{{}_p\Lambda}}$ and the second equivalence follows from $P_{\overrightarrow{{}_p\Lambda}}=\tilde{P}_{\overrightarrow{{}_p\Lambda}}\circ \tilde{s}_{\overrightarrow{{}_p\Lambda}}$. Therefore, it suffices to show that the following diagram commutes
\begin{equation}
\begin{tikzcd}[row sep=1.2cm, column sep=0.8cm]
\Big(\int^{\mathbb{Z}/p}(\iota_{{}_p\Lambda})_!(P^{-1}_{\overrightarrow{{}_p\Lambda}})^*\tilde{s}^*_{\overrightarrow{{}_p\Lambda}}(\overrightarrow{j_{\mathcal{A}_{\infty}^{oper,dg}}})^*A^{\sharp}_L\Big)^{h\mathbb{Z}/p}\arrow[rr,"{\tilde{OC}^{\mathbb{Z}/p,oper}}"]& &QH(M)[[t,\theta]]\\
CC^{\mathbb{Z}/p}(A_L)\arrow[urr,"OC^{\mathbb{Z}/p}"]\arrow[u,"\simeq"]& &
\end{tikzcd},
\end{equation}
where the vertical quasi-isomorphism is given by Proposition 3.19; recall that (cf. Definition 3.17) $CC^{\mathbb{Z}/p}(A_L)=CC^{\mathbb{Z}/p}(\tilde{s}^*_{\overrightarrow{{}_p\Lambda}}(\overrightarrow{j_{\mathcal{A}_{\infty}^{oper,dg}}})^*A^{\sharp}_L)$.
The following is an analogue of Lemma 6.3, and we omit the proof.
\begin{lemma}
There exists \par\indent
1) A universal choice of Floer data for the $\mathbb{Z}/p$-equivariant open-closed map (cf. section 2.4).\par\indent
2) For a cell
\begin{equation}
K\in \mathbb{Z}/p\;\;\textrm{orbit of}\;\;\{\overline{\mathcal{R}}^1_{k_1,\cdots,k_p}\times \Delta_{2i}, \overline{\mathcal{R}}^1_{k_1,\cdots,k_p}\times \Delta_{2i+1}\}_{k_1,\cdots,k_p,i\geq 0},
\end{equation}
one can assign a cubical subdivision of $s([K])$ of $K$, compatible with product and boundary structures.
In particular these assignments fit together to chain maps $s: C_*^{cell}(\overline{\mathcal{R}}^1_{k_1,\cdots,k_p}\times S^{\infty})\rightarrow C_*^{\Box}(\overline{\mathcal{R}}^1_{k_1,\cdots,k_p}\times S^{\infty})$);\par\indent
3) There are chain maps $\tilde{s}: C_*^{cell}(\overline{\mathcal{R}}^1_{k_1,\cdots,k_p}\times S^{\infty})\rightarrow C_*(\tilde{\mathcal{F}}^{reg}(\overline{\mathcal{R}}^1_{k_1,\cdots,k_p}))$, for $k_1,\cdots,k_p\geq 0$, compatible with product structures, that satisfy:
\begin{enumerate}[label=3\alph*)]
    \item these fit into the commutative diagrams
\begin{equation}
\begin{tikzcd}[row sep=1.2cm, column sep=0.8cm]
C_*^{cell}(\overline{\mathcal{R}}^1_{k_1,\cdots,k_p}\times S^{\infty})\arrow[r,"{\tilde{s}}"]\arrow[d,"{s}"]& C_*(\tilde{\mathcal{F}}^{reg}(\overline{\mathcal{R}}^1_{k_1,\cdots,k_p}))\arrow[d,"{\pi}"]\\
C_*^{\Box}(\overline{\mathcal{R}}^1_{k_1,\cdots,k_p}\times S^{\infty})\arrow[r,"{\mathrm{proj}_1}"]& C_*^{\Box}(\overline{\mathcal{R}}^1_{k_1,\cdots,k_p})
\end{tikzcd}
\end{equation}
   \item $\tilde{s}([K])$ has the same degree as the dimension of $K$, and
\begin{equation}
(\overrightarrow{j^{op}_{\mathcal{A}_{\infty}^{oper,dg}}})^*OC^{oper}(\tilde{s}([K])): A^{\otimes k_1+\cdots+k_p+p-1}_L\rightarrow QH
\end{equation}
agrees with the operation induced by counting the paramtrized moduli problem (where the domain curves vary over $\overline{\mathcal{R}}^1_{k_1,\cdots,k_p}$) with respect to the universal choice of Floer data restricted to $K$.
\end{enumerate}
\qed
\end{lemma}
\emph{Proof of Theorem 6.1. 2) given Lemma 6.4}. By Proposition 3.28,
$\Big(\int^{\mathbb{Z}/p}(\iota^{op}_{{}_p\Lambda})_*((P^{op}_{\overrightarrow{{}_p\Lambda}})^{-1})^*(\tilde{s}^{op}_{\overrightarrow{{}_p\Lambda}})^*(\overrightarrow{j^{op}_{\mathcal{A}_{\infty}^{oper,dg}}})^*{}_p\overline{\mathfrak{R}}^1\Big)_{h\mathbb{Z}/p}$ is computed by the cohomology of
\begin{equation}
CC^{\mathbb{Z}/p,\vee}((\tilde{s}^{op}_{\overrightarrow{{}_p\Lambda}})^*(\overrightarrow{j^{op}_{\mathcal{A}_{\infty}^{oper,dg}}})^*{}_p\overline{\mathfrak{R}}^1)=\bigoplus_{k_1,\cdots,k_p\geq 0}C_{-*}^{cell}(\tilde{\mathcal{F}}^{reg}(\overline{\mathcal{R}}^1_{k_1+\cdots+k_p+p-1}))[-k_1-\cdots-k_p][\tilde{t},\tilde{\theta}],
\end{equation}
with differential
\begin{equation}
\begin{cases}
x\mapsto d_{{}_pCC}x\\
x\tilde{t}^k\mapsto d_{{}_pCC}x\tilde{t}^k+(-1)^{|x|}(\tau-1)^{p-1}x\tilde{t}^{k-1}\tilde{\theta}\\
x\tilde{t}^k\tilde{\theta}\mapsto d_{{}_pCC}x\tilde{t}^k\tilde{\theta}+(-1)^{|x|}(\tau-1)x \tilde{t}^k
\end{cases}.
\end{equation}
The fact that $\tilde{s}$ of Lemma 6.4. 3) is a chain map implies that the generator of
$$H^*(\Big(\int^{\mathbb{Z}/p}(\iota^{op}_{{}_p\Lambda})_*((P^{op}_{\overrightarrow{{}_p\Lambda}})^{-1})^*\tilde{s}_{\overrightarrow{{}_p\Lambda}}^*(\overrightarrow{j_{\mathcal{A}_{\infty}^{oper,dg}}})^*{}_p\overline{\mathfrak{R}}^1\Big)_{h\mathbb{Z}/p})\cong HH_*^{\mathbb{Z}/p,\vee}(\tilde{s}_{\overrightarrow{{}_p\Lambda}}^*(\overrightarrow{j_{\mathcal{A}_{\infty}^{oper,dg}}})^*{}_p\overline{\mathfrak{R}}^1)\cong k[\tilde{t},\tilde{\theta}]$$
corresponding to $\tilde{t}^i\tilde{\theta}$
has a chain representative in $CC^{\mathbb{Z}/p,\vee}(\tilde{s}_{\overrightarrow{{}_p\Lambda}}^*(\overrightarrow{j_{\mathcal{A}_{\infty}^{oper,dg}}})^*{}_p\overline{\mathfrak{R}}^1)$ given by
\begin{equation}
\sum_{j=0}^i\Big(\prod \tilde{s}([\overline{\mathcal{R}}^1_{k_1,\cdots,k_p}]\times \Delta_{2j})\tilde{\theta}+\prod \tilde{s}([\overline{\mathcal{R}}^1_{k_1,\cdots,k_p}]\times \Delta_{2j+1})\Big) \tilde{t}^{i-j},
\end{equation}
and the generator corresponding to $\tilde{t}^i$ has a representative
\begin{equation}
\sum_{j=0}^i\Big(\prod \tilde{s}([\overline{\mathcal{R}}^1_{k_1,\cdots,k_p}]\times \Delta_{2j})+(\tau-1)^{p-2}\prod \tilde{s}([\overline{\mathcal{R}}^1_{k_1,\cdots,k_p}]\times \Delta_{2j-1})\tilde{\theta}\Big) \tilde{t}^{i-j}.
\end{equation}
By Lemma 6.4. 3b), under the vertical quasi-isomorphism of (6.60), $\tilde{OC}^{\mathbb{Z}/p,oper}_{2i}$ (resp. $\tilde{OC}^{\mathbb{Z}/p,oper}_{2i+1}$) agrees with the $t^i$-th (resp. $t^i\theta$-th) coefficient of the $\mathbb{Z}/p$-equivariant open-closed map $OC^{\mathbb{Z}/p}$ defined in section 2.4.\qed\par\indent
6.4. \textbf{Proof of Theorem 6.1. 3)}. This reduces to the commutativity of
\begin{equation}
\begin{tikzcd}[row sep=1.2cm, column sep=0.8cm]
\Big(\int_{\mathbb{T}}(\iota_{\Lambda})_!(P^{-1}_{\overrightarrow{\Lambda}})^*\tilde{s}^*_{\overrightarrow{\Lambda}}A_L^{\sharp}\Big)^{h\mathbb{T}}\langle 1,\theta\rangle\arrow[d]\arrow[rr,"\phi_p"]& &\Big(\int_{\mathbb{Z}/p}j^*(\iota_{{}_p\Lambda})_!(P^{-1}_{\overrightarrow{\Lambda}})^*\tilde{s}^*_{\overrightarrow{\Lambda}}A_L^{\sharp}\Big)^{h\mathbb{Z}/p}\arrow[d]\\
\Big(\int_{\mathbb{T}}(\iota_{\Lambda})_!(\tilde{P}^{-1}_{\overrightarrow{\Lambda}})^*A_L^{\sharp}\Big)^{h\mathbb{T}}\langle1,\theta\rangle\arrow[rr,"\phi_p"] & &\Big(\int_{\mathbb{Z}/p}j^*(\iota_{{}_p\Lambda})_!(\tilde{P}^{-1}_{\overrightarrow{\Lambda}})^*A_L^{\sharp}\Big)^{h\mathbb{Z}/p}
\end{tikzcd},
\end{equation}
where the top horizontal arrow is $\phi_p$ applied to $(\iota_{\Lambda})_!(P^{-1}_{\overrightarrow{\Lambda}})^{op}\tilde{s}^*_{\overrightarrow{\Lambda}}A_L^{\sharp}$, and the bottom horizontal arrow is $\phi_p$ applied to $(\iota_{\Lambda})_!(\tilde{P}^{-1}_{\overrightarrow{\Lambda}})^{op}A_L^{\sharp}$; the vertical arrows are induced by $\tilde{P}_{\overrightarrow{\Lambda}}\circ\tilde{s}_{\overrightarrow{\Lambda}}=P_{\overrightarrow{\Lambda}}$. The commutativity of (6.68) is then a consequence of the naturality of $\phi_p$ proved in Proposition 5.3. 1).\qed

\begin{appendices}
\renewcommand{\theequation}{A.\arabic{equation}}
\setcounter{equation}{0}
\section{From operad to multicategory}
Throughout section 3 to section 6 of this paper, we worked in a simplified setting by only consider one object $L$ of the monotone Fukaya category $\mathrm{Fuk}(M)_{\lambda}$. We showed that $A:=CF^*(L,L)$ is naturally a dg algebra over the dg operad $\mathcal{A}_{\infty}^{oper,dg}$, cf. Lemma 4.6. Moreover, we endowed $A$ with the structure of a classical $\mathcal{A}_{\infty}$-algebra by choosing a `section' $\tilde{s}: \mathcal{A}_{\infty}^{dg}\rightarrow \mathcal{A}_{\infty}^{oper,dg}$, cf. section 6.1. \par\noindent
In Appendix A, we generalize this by allowing multiple objects of the monotone Fukaya category. To do this, we need to recall the mutliple-object version of an operad, a \emph{multicategory}.
\begin{mydef}
A multicategory $\mathcal{C}$ enriched in a symmetric monoidal category $(V,\otimes,I)$ consists of
\begin{enumerate}[label=\arabic*)]
    \item an object set $X$,
\item an object $\mathcal{C}(\overrightarrow{x},y)\in V$, considered as the $n$-ary multimorphism object from $\overrightarrow{x}\in X^n$ to $y\in X$,
\item a distinguished morphism $id_x: I\rightarrow \mathcal{C}(x,x)$,
\item multimorphism composition laws
\begin{equation}
\circ_i: \mathcal{C}(\overrightarrow{x}, y)\circ_i\mathcal{C}(\overrightarrow{y},z)\rightarrow \mathcal{C}(\overrightarrow{x}\circ_i\overrightarrow{y},z),
\end{equation}
where $\overrightarrow{x}\circ_i\overrightarrow{y}$ denotes the replacement of the $i$-th element in $\overrightarrow{y}$ by the sequnce $\overrightarrow{x}$.
\end{enumerate}
The compositions are required to satisfy associativity and identity axioms similar to that of an operad, cf. \cite[Definition B.8]{AGV}.
\end{mydef}
An operad is a multicategory with one object. \par\indent
Fix a set of objects $\mathcal{L}=\{L_i\}_{i\in I}$ in $\mathrm{Fuk}(M)_{\lambda}$.
\begin{mydef}
The dg multicatgeory $\mathcal{A}_{\infty,\mathcal{L}}^{dg}$ has
\begin{itemize}
    \item  objects given by pairs $(L_0,L_1), L_0,L_1\in \mathcal{L}$
    \item multimorphism complexes given by $\mathcal{A}_{\infty,\mathcal{L}}^{dg}\Big(((L_0^0,L_1^0),(L^1_0,L^1_1),\cdots,(L^d_0,L^d_1)),(L_0,L_1)\Big)=0$ unless $L_0=L_0^0, L_1^0=L^0_1, L_2^0=L^1_1,\cdots, L^{d-1}_1=L^d_0, L^d_1=L_1$, in which case
\begin{equation}
\mathcal{A}_{\infty,\mathcal{L}}^{dg}\Big(((L_0,L_1),(L_1,L_2),\cdots,(L_{d-1},L_d)),(L_0,L_d)\Big):=(\mathcal{A}_{\infty}^{dg})_{d}.
\end{equation}
\end{itemize}
Recall this is just a copy of $C_{-*}^{cell}(\overline{\mathcal{R}}^{d+1})$ if $d>1$ and $k$ (generated by the identity morphism) if $d=1$. However, we think of the disks being equipped with Lagrangians $L_0,\cdots,L_d$ labeling the boundary components. The multimorphism compositions are given by operadic compositions of $\mathcal{A}_{\infty}^{dg}$, but we record the Lagrangian labels as we concatenate disks.
\end{mydef}
We now consider the version of $\mathcal{A}_{\infty,\mathcal{L}}^{dg}$ where the Reimann surfaces are equipped with Floer data. Recall from section 2.1 that each pair of objects $(L_0,L_1)$ is associated with a time-dependent Hamiltonian $H_{L_0,L_1}$ and almost complex structure $J_{L_0,L_1}$. Given a sequence $\mathbf{L}=(L_0,\cdots,L_d)$ of objects in $\mathrm{Fuk}(M)_{\lambda}$, we define a $0$-cube of the symmetric cubical set $\mathcal{F}_{\mathbf{L}}(\overline{\mathcal{R}}^{d+1})$ to consists of data: (compare with section 4.2)
\begin{enumerate}
    \item  A stable disk $\Sigma\in \overline{\mathcal{R}}^{d+1}$, a labeling of each interior and boundary marked points as input/output, and a labeling of each boundary component of $\Sigma$ using elements of $\mathbf{L}$.
    \item For each component $\Sigma_v$ of $\Sigma$, and each boundary marked point $p$ of $\Sigma_v$, a choice of strip-like ends at $p$
\begin{equation}
\epsilon_p^{+}:[0,\infty)\times [0,1]\rightarrow \Sigma_v\;\;\mathrm{or}\;\;\epsilon_p^-:(-\infty,0]\times[0,1]\rightarrow \Sigma_v
\end{equation}
depending on whether $p$ is an output or input.
\item For each component $\Sigma_v$, a pair $(K_v,J_v)$ where $K_v\in \Omega^1(\Sigma_v, \mathcal{H}), J_v\in C^{\infty}(\Sigma_v,\mathcal{J})$ satisfying
\begin{equation}
(\epsilon^{\pm}_p)^*K_v=H_{t,L_0,L_1}dt\;\;,\;\;(\epsilon^{\pm}_p)^*J_v=J_{t,L_0,L_1},
\end{equation}
$L_0,L_1$ are the two Lagrangians incidence at $p$.
\end{enumerate}
The definition of a general $n$-cube that follows exactly as in section 4.2. There is a version of this consisting of regular cubes (cf. Definition 4.4) which we denote as $\tilde{F}^{reg}(\overline{\mathcal{R}}^{d+1})$.
\begin{mydef}
The dg multicatgeory $\mathcal{A}_{\infty,\mathcal{L}}^{oper,dg}$ has
\begin{itemize}
    \item  objects given by pairs $(L_0,L_1), L_0,L_1\in \mathcal{L}$
    \item multimorphism complexes given by $\mathcal{A}_{\infty,\mathcal{L}}^{oper, dg}\Big(((L_0^0,L_1^0),(L^1_0,L^1_1),\cdots,(L^d_0,L^d_1)),(L_0,L_1)\Big)=0$ unless $L_0=L_0^0, L_1^0=L^0_1, L_2^0=L^1_1,\cdots, L^{d-1}_1=L^d_0, L^d_1=L_1$, in which case
\begin{equation}
\mathcal{A}_{\infty,\mathcal{L}}^{oper, dg}\Big(((L_0,L_1),(L_1,L_2),\cdots,(L_{d-1},L_d)),(L_0,L_d)\Big):=C_{-*}(\tilde{F}^{reg}_{(L_0,\cdots,L_d)}(\overline{\mathcal{R}}^{d+1})),
\end{equation}
if $d>1$, and $:=k$ (generated by the identity morphism) if $d=1$.
\end{itemize}
The multimorphism compositions are induced by concatenation of disks equipped with Floer data (cf. (4.14)), while recording the Lagrangian labelings on the concatenation.
\end{mydef}
Similar to the proof of Lemma 4.6, one can show that
\begin{lemma}
There is a map of dg multicategories
\begin{equation}
\mathcal{A}_{\infty,\mathcal{L}}^{oper,dg}\rightarrow Mod_k
\end{equation}
which on objects is given by $(L_0,L_1)\mapsto CF^*(L_0,L_1,H_{L_0,L_1},J_{L_0,L_1})$.
\end{lemma}
We call (A.6) the \emph{operadic monotone Fukaya catgeory with set of objects $\mathcal{L}$}. When $\mathcal{L}=\{L\}$, this recovers (4.15). \par\indent
Analogous to Lemma 6.2, there is a homotopy equivalence of dg multicategories $\tilde{s}: \mathcal{A}_{\infty,\mathcal{L}}^{dg}\rightarrow \mathcal{A}_{\infty,\mathcal{L}}^{oper,dg}$. Precomposing (6) with $\tilde{s}$, we obtain a map of dg multicategories $\mathcal{A}_{\infty,\mathcal{L}}^{dg}\rightarrow Mod_k$, which recovers the classical notion of the monotone Fukaya $\mathcal{A}_{\infty}$-category with set of objects $\mathcal{L}$.\par\indent
Finally, note that we can also readily adapt the variants of cyclic categories to the multiple object setting. As an example, we consider the following generalization of $\overrightarrow{\Lambda}\rtimes \mathcal{A}_{\infty}^{oper,dg}$.
\begin{mydef}
Define $\overrightarrow{\Lambda}\rtimes \mathcal{A}_{\infty,\mathcal{L}}^{oper,dg}$ to be the following dg category.
\begin{itemize}
    \item Objects $([n],(L_0,\cdots,L_n))$ are specified by a nonnegative integer $n$ and a sequence of $n+1$ Lagrangians in $\mathrm{Fuk}(M)_{\lambda}$. We think of this sequence as a Lagrangian labeling of the boundary components of the disk corresponding to $[n]$.
    \item  The morphism chain complex from $([n],(L_0,\cdots,L_n))$ to $([m],(L'_0,\cdots,L'_m))$ is defined to be
\begin{equation}
\bigoplus_{f\in \overrightarrow{\Lambda}([n],[m])} \bigotimes_{i=0}^m \mathcal{A}_{\infty,\mathcal{L}}^{oper,dg}(\mathbf{L}_{f^{-1}(i)}, \mathbf{L}'_i),
\end{equation}
where $\mathbf{L}'_i$ denotes the pair of Lagrangians incident to the $i$-th marked point on the disk corresponding to $[m]$, and $\mathbf{L}_{f^{-1}(i)}$ denotes the sequence of pairs of Lagrangians incident to the sequence of marked points in $f^{-1}(i)$.
\end{itemize}
Composition combines the composition of maps in $\overrightarrow{\Lambda}$ with the multicategory structure of $\mathcal{A}^{oper,dg}_{\infty,\mathcal{L}}$ in a way that is analogous to Definition 3.1.
\end{mydef}
Moreover, there is an obvious generalization of the Hochschild functor
\begin{equation}
\overrightarrow{\Lambda}\rtimes \mathcal{A}_{\infty,\mathcal{L}}^{oper,dg}\rightarrow \mathrm{Mod}_k
\end{equation}
which on objects sends
\begin{equation}
([n],(L_0,\cdots,L_n))\mapsto CF^*(L_0,L_1)\otimes CF^*(L_1,L_2)\otimes\cdots\otimes CF^*(L_n,L_0).
\end{equation}
Using this framework, it is easy to see that all the results in this paper generalizes to monotone Fukaya category with multiple objects.

\renewcommand{\theequation}{B.\arabic{equation}}
\setcounter{equation}{0}
\section{Simplicial sets and multi-simplicial sets}
\subsection*{Simplicial sets}
In this subsection, we recall some basic constructions of simplicial sets. Let $\Delta$ be the category whose objects are $[n],n\in\mathbb{N}$, where $[n]$ is viewed as
partially ordered set $[0<1<\cdots<n]$, and whose morphisms are ordering preserving functions.
\begin{mydef}
A \emph{simplicial set} is a functor $X: \Delta^{op}\rightarrow \mathrm{Sets}$.
\end{mydef}
One can define a more general notion of simplicial objects in a category by replacing Sets with the desired category. \par\indent
Let $\Delta^n=\{(t_0,t_1,\cdots,t_n)|0\leq t_i\leq 1, \sum_{i=0}^n t_i=1\}\subset \mathbb{R}^{n+1}$ be the standard $n$-simplex. Then there are linear maps $\delta_i: \Delta_{n-1}\rightarrow \Delta_n, 0\leq i\leq n$
induced by skipping the $i$-th vertex of $\Delta^n$, and $\sigma_i: \Delta_{n+1}\rightarrow \Delta_n, 0\leq i\leq n$ induced by doubling the $i$-th vertex of $\Delta^n$.
\begin{mydef}
Let $X$ be a simplicial set. The \emph{geometric realization} of $X$ is defined to be the topological space
\begin{equation}
|X|:=\coprod X_n\times \Delta^n/\sim,
\end{equation}
where $\sim$ is the equivalence relation generated by
\begin{equation}
(d_i(x),y)\sim(x,\delta_i(y)),
\end{equation}
\begin{equation}
(s_i(x),y)\sim(x,\sigma_i(y)).
\end{equation}
\end{mydef}

In practice, due to the non-existence of strict units, it is often convenient to consider the notion of semi-simplicial sets. To this end, let $\overrightarrow{\Delta}\subset \Delta$
be the subcategory of $\Delta$ with the same objects, but morphisms that are \emph{injective} order preserving functions.
\begin{mydef}
A \emph{semi-simplicial} set is a functor $X:\overrightarrow{\Delta}^{op}\rightarrow \mathrm{Sets}$.
\end{mydef}
Combinatorially, it is a collection of sets $\{X_n\}_{n\geq 0}$ with face maps $d_i:X_n\rightarrow X_{n-1}, 0\leq i\leq n$ that satisfy the relations in (B.1). In other words, a semi-simplcial
set is a simplicial set without degeneracies.
\begin{mydef}
Let $X$ be a semi-simplicial set. The \emph{geometric realization} of $X$ is defined to be the topological space
\begin{equation}
|X|:=\coprod X_n\times \Delta^n/\sim,
\end{equation}
where $\sim$ is the equivalence relation generated by
\begin{equation}
(d_i(x),y)\sim(x,\delta_i(y)).
\end{equation}
\end{mydef}

Most of the (semi)-simplicial objects we consider in this paper are (semi)-simplicial chain complexes. Let $\mathrm{Mod}_k$ denote the category of chain complexes over $k$.
\begin{mydef}
Let $X: \overrightarrow{\Delta}^{op}\rightarrow \mathrm{Mod}_k$ be a semi-simplicial chain complex over $k$. As a graded vector space, the \{cyclic bar complex\} of $X$, denoted $CC(X)$, is
\begin{equation}
CC(X):=\bigoplus_{n\geq 0} X_n[n].
\end{equation}
The \emph{bar differential} on $CC(X)$ is defined as $b'=\sum_{i=0}^{-1}(-1)^i d_i$ and the \emph{cyclic bar differential} is defined as $b=\sum_{i=0}^n(-1)^id_i$. The complexes $(CC(X),b')$ and
$(CC(X),b)$ are called the \emph{bar complex} and \emph{cyclic bar complex} of $X$, respectively.
\end{mydef}
If $X:\Delta^{op}\rightarrow \mathrm{Mod}_k$ is a simplcial chain complex, we define its bar and cyclic bar complex to be that of its underlying semi-simplicial chain complex. However,
in this case there is a subcomplex $D(X)\subset (CC(X),b)$ consisting of degenerate elements (elements in the image of some $s_i$). $D(X)$ is called the \emph{degenerate subcomplex}, and
it is well known that $D(X)$ is acyclic. In particular, the projection
\begin{equation}
(CC(X),b)\rightarrow (CC(X),b)/D(X)=:\overline{CC}(X)
\end{equation}
is a quasi-isomorphism. $\overline{CC}(X)$ is called the \emph{normalized cyclic bar complex} of $X$.

\begin{rmk}
More often in the simplicial literature, $(CC(X),b)$ is called the standard chain complex of $X$, and $\overline{CC}(X)$ the normalized chain complex of $X$.
Our choice of terminology and notation is motivated by the study of Hochschild homology using simplicial methods. In particular, it comes from the following example: let $A$ be a
strictly unital dg algebra over $k$. Then there is a simplicial chain complex $A^{\sharp}$ defined by $[n]\mapsto A^{\otimes n+1}$,
\begin{equation}
d_i(a_0\otimes a_1\otimes\cdots\otimes a_n)=a_0\otimes \cdots \otimes a_ia_{i+1}\otimes \cdots\otimes a_n,
\end{equation}
\begin{equation}
s_i(a_0\otimes a_1\otimes\cdots\otimes a_n)=a_0\otimes \cdots\otimes a_i\otimes 1\otimes\cdots\otimes a_n.
\end{equation}
Then, $CC(A^{\sharp})=CC(A)$ is the cyclic bar complex computing Hochschild homology of $A$. In section 2, we generalize this construction to $A_{\infty}$-algebras.
\end{rmk}

\subsection*{Multi-simplicial sets}
We now consider a slight generalization of simplicial set. Let $N\geq 1$ be an integer.
\begin{mydef}
An \emph{N-fold simplicial set} is a functor $X: (\Delta^{op})^p\rightarrow \mathrm{Sets}$.
\end{mydef}
Combinatorially, $X$ is the data of a set $X_{k_1,\cdots,k_N}$ for each tuple of non-negative integers $(k_1,\cdots,k_N)$, together with face maps
\begin{equation}
d^l_i:X_{k_1,\cdots,k_l,\cdots,k_N}\rightarrow X_{k_1,\cdots,k_l-1,\cdots,k_N}, 0\leq i\leq k_l, 1\leq l\leq N
\end{equation}
and degeneracy maps
\begin{equation}
s^l_i:X_{k_1,\cdots,k_l,\cdots,k_N}\rightarrow X_{k_1,\cdots,k_l+1,\cdots,k_N}, 0\leq i\leq k_l, 1\leq l\leq N
\end{equation}
such that for each $l$, the maps $d^l_i, s^l_i$ satisfy the relations (B.1)-(B.5).
\begin{mydef}
Let $X$ be an $N$-fold simplicial set. The \emph{geometric realization} of $X$ is defined as the topological space
\begin{equation}
|X|:=\coprod X_{k_1,\cdots,k_N}\times \Delta^{k_1}\times\cdots\times \Delta^{k_N}/\sim,
\end{equation}
where $\sim$ is the equivalence relation generated by
\begin{equation}
(d^l_i(x),y_1,\cdots,y_l,\cdots,y_n)\sim(x,y_1,\cdots,\delta_i(y_l),\cdots,y_n),
\end{equation}
\begin{equation}
(s^l_i(x),y_1,\cdots,y_l,\cdots,y_n)\sim(x,y_1,\cdots,\sigma_i(y_l),\cdots,y_n).
\end{equation}
\end{mydef}
For example, the geometric realization of the representable $N$-fold simplicial set $\Delta^N(-,[k_1,\cdot,k_N])$ is homeomorphic to $\Delta^{k_1}\times\cdots\times \Delta^{k_N}$.
\begin{mydef}
Let $X\in \mathrm{Fun}((\Delta^{op})^p,\mathrm{Mod}_k)$ be an $N$-fold simplicial chain complex. As a graded vector space, the \emph{N-fold cyclic bar complex} of $X$, denoted ${}_NCC(X)$, is
\begin{equation}
{}_NCC(X):=\bigoplus_{k_1,\cdots,k_N\geq 0} X_{k_1,\cdots,k_N}[k_1+\cdots+k_N].
\end{equation}
The \emph{$N$-fold cyclic bar differential} on ${}_NCC(X)$ is given by
\begin{equation}
b^N:=\sum_{l=0}^N\sum_{i=0}^{k_l}(-1)^{i+k_1+\cdots+k_{l-1}}d^l_i.
\end{equation}
The complex $({}_NCC(X),b^N)$ is called the \emph{N-fold cyclic bar complex} of $X$.
\end{mydef}
As in the ordinary simplicial case, one can define the $N$-fold bar differential, the degenerate subcomplex $D_N(X)\subset ({}_NCC(X),b^N)$ and the quotient ${}_N\overline{CC}(X)$.
It should also be clear that the $N$-fold (cyclic) bar complex can be defined for any $N$-fold semi-simplicial chain complex $X:(\overrightarrow{\Delta}^{op})^p\rightarrow \mathrm{Mod}_k$.\par\indent
Finally, there is a functor $o:\Delta^p\rightarrow \Delta$, called the \emph{ordinal sum functor}, which on objects is given by
\begin{equation}
[k_1,\cdots,k_N]\mapsto [k_1+\cdots+k_N+N-1].
\end{equation}
One should think of this as concatenating the partially ordered sets $[0<1<\cdots<k_1],[0<1<\cdots<k_2],\cdots,[0<1<\cdots<k_N]$ one after another into a new partially order set
\begin{equation}
[0_1<1_1<\cdots<(k_1)_1<0_2<1_2<\cdots<(k_2)_2<\cdots<0_N<1_N<\cdots<(k_N)_N],
\end{equation}
which is just a copy of $[k_1+\cdots+k_N+N-1]$. Pulling back along this functor induces a functor
\begin{equation}
\mathrm{Dec}:\mathrm{Fun}(\Delta^{op},\mathrm{Sets})\rightarrow \mathrm{Fun}((\Delta^{op})^{op},\mathrm{Sets})
\end{equation}
called \emph{total decalage}.
The following fact is well known.
\begin{lemma}
Let $X$ be a simplicial set, then there is a natural homotopy equivalence
\begin{equation}
|X|\simeq |\mathrm{Dec}(X)|.
\end{equation}
\qed
\end{lemma}
\begin{rmk}
There is a (strict) $\mathbb{Z}/p$-action on $(\Delta^{op})^p$ given by permuting the $p$ factors. One can heuristically form the `cross product' of $(\Delta^{op})^p$ with $\mathbb{Z}/p$ and get
a new category ${}_p\Lambda$ (see section 2); specifically the classifying space of ${}_p\Lambda$ is homotopy equivalent to $B\mathbb{Z}/p$. There is also an `$S^1$-action' on $\Delta^{op}$ (in a less obvious fashion), and heuristically, the cross product of $\Delta^{op}$ with $S^1$ is
Connes' cyclic category $\Lambda$, whose classifying space is well known to be homotopy equivalent to $BS^1$. The ordinal sum functor $(\Delta^{op})^p\rightarrow \Delta^{op}$ can be enhanced
to a functor ${}_p\Lambda\rightarrow \Lambda$, which on classifying spaces is (homotopic to) the natural map $B\mathbb{Z}/p\rightarrow BS^1$. This will be the key in section 2 for studying
cyclic objects and the underlying finite cyclic object.
\end{rmk}

\renewcommand{\theequation}{C.\arabic{equation}}
\setcounter{equation}{0}
\section{Symmetric cubical sets}
The use of symmetric cubical sets in this paper is solely to apply the operadic Floer theory of \cite{AGV}, which uses the cubical model for homotopy theory of spaces. As such, we follow
Appendix B of loc.cit. for a brief review on the subject. Roughly speaking, a symmetric cubical set $X$ is a collection of sets $\{X_n\}_{n\geq 0}$ with maps between them that model
the types of maps between cubes $[0,1]^n\rightarrow [0,1]^m$ consisting of (i) projection onto some coordinates, (ii) permutation of coordinates and (iii) inclusion of faces. More precisely,
\begin{mydef}
A symmetric cubical set $X$ is a sequence of sets $\{X_n\}_{n\geq 0}$ together with a collection of \emph{face maps}
\begin{equation}
d^{\pm}_{n,i}: X_n\rightarrow X_{n-1},\;\;\;n\geq 1,\;\;1\leq i\leq n,
\end{equation}
\emph{degeneracy maps}
\begin{equation}
s_{n-1,i}:X_{n-1}\rightarrow X_n,\;\;\;n\geq 1,\;\;1\leq i\leq n,
\end{equation}
and \emph{transposition maps}
\begin{equation}
p_{n,i}:X_n\rightarrow X_n,\;\;n\geq 2,\;\;1\leq i\leq n-1.
\end{equation}
There are required to satisfy the following relations for $\mu,\nu\in\{+,-\}$:
\begin{equation}
d^{\mu}_{n-1,i}\circ d^{\nu}_{n,j}=d^{\nu}_{n-1,j-1}\circ d^{\mu}_{n,i},\;\;\;\;\;\;i<j,
\end{equation}
\begin{equation}
s_{n,i}\circ s_{n-1,j}=s_{n,j+1}\circ s_{n-1,i},\;\;\;\;\;\;i\leq j,
\end{equation}
\begin{equation}
p_{n,i}^2=\mathrm{id},\;\;(p_{n,i}\circ p_{n,i+1})^3=\mathrm{id},
\end{equation}
\begin{equation}
p_{n,i}\circ p_{n,j}=p_{n,j}\circ p_{n,i},\;\;\;\;\;\;\;i+1<j,
\end{equation}
\begin{equation}
d^{\mu}_{n,i}\circ s_{n-1,j}=
\begin{cases}
s_{n-2,j-1}\circ d^{\mu}_{n,i},\;\;\;\;i<j,\\
s_{n-2,j}\circ d^{\mu}_{n,i-1},\;\;\;\;i>j,\\
\mathrm{id},\;\;\;\;\;\;\;\;\;\qquad\qquad i=j,
\end{cases}
\end{equation}
\begin{equation}
d^{\mu}_{n,j}\circ p_{n,i}=
\begin{cases}
p_{n-1,i-1}\circ d^{\mu}_{n,j},\;\;\;\;j<i,\\
d^{\mu}_{n,i+1},\;\;\;\;\;\;\quad \qquad j=i,\\
d^{\mu}_{n,i},\;\;\;\;\;\;\quad \qquad \;\;\;\, j=i+1,\\
p_{n-1,i}\circ d^{\mu}_{n,j},\;\;\;\,\quad j>i+1,
\end{cases}
\end{equation}
\begin{equation}
p_{n,i}\circ s_{n-1,j}=
\begin{cases}
s_{n-1,j}\circ p_{n-1,i-1},\;\;\;\;j<i,\\
s_{n-1,i+1},\;\;\;\;\;\;\quad\qquad j=i,\\
s_{n-1,i},\;\;\;\;\;\;\;\quad \qquad \;\;\;j=i+1,\\
s_{n-1,j}\circ p_{n,i},\;\;\;\;\qquad j>i+1.
\end{cases}
\end{equation}
\end{mydef}
Note that the transposition maps $p_{n,i}$ generate an action of the symmetric group $S_n$ on $X_n$. If we forget about the transposition maps, the notion of an ordinary cubical set is
recovered. However, the symmetry provided by the $p_{n,i}$'s is crucial in defining a symmetric monoidal product for cubical sets, which is needed for a category over which we study operads.\\\par\indent
\emph{Example}. Let $X$ be a topological space. The \emph{singular (symmetric) cubical set} of $X$, denoted $\Box(X)$, has $n$-cubes $\Box_n(X)$ the set of all
continuous maps $[0,1]^n\rightarrow X$. For $\sigma\in \Box_n(X)$ we define $d^{\pm}_{n,i}(\sigma)$ to be the precomposition of $\sigma$ with the inclusion
$\iota^{\pm}_{n,i}:[0,1]^{n-1}\rightarrow [0,1]^n$ as the $i$-th front (if $+$) or back (if $-$) face of $[0,1]^n$.
We define $s_{n,i}(\sigma)$ to be $\sigma$ precomposed with the projection $\pi_{n,i}:[0,1]^{n+1}\rightarrow [0,1]^{n}$ forgetting the $i$-th coordinate. The transposition $p_{n,i}(\sigma)$ is
defined to be $\sigma$ precomposed with the map $\tau_{n,i}:[0,1]^n\rightarrow [0,1]^n$ that transposes the $i$th and $i+1$-th coordinates.

\begin{mydef}
Let $X$ be a cubical set, the \emph{geometric realization} $|X|$ is defined as the topological space
\begin{equation}
\coprod X_n\times [0,1]^n/\sim,
\end{equation}
where $\sim$ is the equivalence relation generated by
\begin{equation}
(d^{\pm}_{n,i}(x),y)\sim(x,\iota^{\pm}_{n,i}(y)),\;\;\;(s_{n,i}(x),y)\sim(x,\pi_{n,i}(y)),\;\;\;(p_{n,i}(x),y)\sim(x,\tau_{n,i}(y)).
\end{equation}
\end{mydef}

\subsubsection*{The symmetric monoidal product}
Let $X^1$ and $X^2$ be symmetric cubical sets. Recall in particular that $(X^1)_n,(X^2)_n$ are equipped with $S_n$-actions. We define a new symmetric cubical set $X^1\otimes X^2$ as follows.
\begin{equation}
(X^1\otimes X^2)_n:=\coprod_{n_1+n_2=n}S_n\times_{S_{n_1}\times S_{n_2}}(X^1_{n_1}\times X^2_{n_2}/\sim),
\end{equation}
where $\sim$ is the equivalence relation generated by
\begin{equation}
(s_{n_1-1,n_1}(x^1),x^2)\sim(x^1,(s_{n_2-1,1}(x^2)).
\end{equation}
We omit the definition of the face, degeneracy and transposition maps and instead refer the readers to \cite[B.1.1]{AGV}. $\otimes$ defines a symmetric monoidal product on
the category of symmetric cubical sets, and it satisfies a universal property analogous to that of the tensor product of abelian groups:
\begin{lemma}
There is a natural bijection between maps of symmetric cubical set from $X^1\otimes X^2$ to $X$ and collection of maps $X^1_{n_1}\times X^2_{n_2}\rightarrow X_{n_1+n_2}$, for all $n_1,n_2$, that are\par\indent
(1) $S_{n_1}\times S_{n_2}$-equivariant,\par\indent
(2) intertwine the face maps $d^{\pm}_{i,1}$ with $d_i^{\pm}$ and $d^{\pm}_{j,2}$ with $d^{\pm}_{n_1+j}$,\par\indent
(3) intertwine the degeneracy maps $s_{i,1}$ with $s_i$ and $s_{j,2}$ with $s_{n_1+j}$. \qed
\end{lemma}

Finally, we recall the notion of symmetric normalized cubical chains of a symmetric cubical set.
\begin{mydef}
Let $X$ be a symmetric cubical set, and fix a coefficient ring $k$. Its \emph{symmetric normalized cubical chain} $C_*(X;k)$ is define as
\begin{equation}
C_n(X;k):=\frac{k[X_n]}{\sum_{i=1}^n\mathrm{Im}(s_{n-1,i})+\sum_{i=1}^n\mathrm{Im}(1+p_{n,i})}.
\end{equation}
The differential is given by
\begin{equation}
d=\sum_{i=1}^n(-1)^i(d_i^+-d_i^-).
\end{equation}
\end{mydef}
We refer to \cite[Appendix B]{AGV} for the following lemmas.
\begin{lemma}
The symmetric normalized cubical chain functor is symmetry monoidal.
\end{lemma}
\begin{lemma}
The symmetric normalized singular chain functor, i.e. the composition $C_*(\Box(-)):\mathrm{Top}\rightarrow \mathrm{Mod}_k$, is symmetric monoidal. Moreover, for a topological space $X$,
the homology of $C_*(\Box(X))$ is isomorphic to its singular homology.
\end{lemma}

\renewcommand{\theequation}{D.\arabic{equation}}
\setcounter{equation}{0}

\section{Backgrounds in $\infty$-category}
$\infty$-categories (more precisely $(\infty,1)$-categories) are `categories' with a collection of objects, $1$-morphisms, $2$-morphisms and so on that encode higher homotopical data.
In other words, whereas ordinary category has a \emph{set} of morphisms between any two objects, an $\infty$-category has a \emph{space} of morphisms. These higher structures
allow one to capture the geometry involved in a category often unseen from an ordinary category point of view. For instance, in section 2, we see that the mechanism of $\infty$-categories
naturally extracts an action of the circle group (an $\infty$-group) from any cyclic object. \par\indent
The model we use for $\infty$-categories in Joyal's theory of quasi-category, which has been extensively developed in \cite{Lur1},\cite{Lur2}.
\begin{mydef}
A \emph{quasi-category} is a simplicial set $X$ that has the right lifting property with respect to all inner horns $\Lambda^n_i\hookrightarrow \Delta^n, 0<i<n, n\geq 0$.
\end{mydef}
By an abuse of terminology, we will refer to a quasi-category as simply an $\infty$-category. Any category $\mathcal{C}$ has an associated $\infty$-category $N\mathcal{C}$
called its \emph{nerve}, and moreover functors from $\mathcal{C}$ to $\mathcal{D}$ are in bijection with functors (maps of simplicial sets) from $N\mathcal{C}$ to $N\mathcal{D}$. Therefore,
the theory of ordinary categories is subsumed by that of $\infty$-categories.\par\indent
D.1. \textbf{The dg nerve and the derived category}.
In symplectic geometry, we often work with dg categories, or more generally $A_{\infty}$-categories (e.g. the Fukaya category). Below we recall a construction called the \emph{dg nerve} (cf \cite[Construction 1.3.1.6]{Lur2}),
that produces an $\infty$-category out of a dg category.\\\\
\textbf{Construction}. Let $\mathcal{C}$ be a ($\mathbb{Z}$ or $\mathbb{Z}/2$-graded) dg category over a field $k$. We associate to $\mathcal{C}$ an $\infty$-category $N^{dg}(\mathcal{C})$, called
the \emph{dg nerve} of $\mathcal{C}$, defined as follows. For $n\geq 0$, define $N^{dg}(\mathcal{C})_n$ to be the set of order pairs $(\{X_i\}_{0\leq i\leq n},\{f_I\})$, where\par\indent
(a) For $0\leq i\leq n$, $X_i$ is an object of $\mathcal{C}$.\par\indent
(b) For every subset $I=\{i_-<i_1<\cdots<i_m<i_+\}\subset [n]$, with $m\geq 0$, an element $f_I\in \mathrm{Map}^{-m}(X_{i_-},X_{i_+})$, satisfying the equation
\begin{equation}
df_I=\sum_{1\leq j\leq m}(-1)^j(f_{I-\{i_j\}}-f_{i_j<\cdots<i_m<i_+}\circ f_{i_-<i_1<\cdots<i_j}).
\end{equation}
If $\alpha:[m]\rightarrow [n]$ is an order preserving function, then induced map $N^{dg}(\mathcal{C})_n\rightarrow N^{dg}(\mathcal{C})_m$ is given by
\begin{equation}
(\{X_i\}_{0\leq i\leq n},\{f_I\})\mapsto (\{X_{\alpha(j)}\}_{0\leq j\leq m},\{g_J\}),
\end{equation}
\begin{equation}
\textrm{where}\;g_J=
\begin{cases}
f_{\alpha(J)} &\textrm{if $\alpha|_J$ is injective}\\
\mathrm{id}_{X_i}&\textrm{if $J={j,j'}$ with $\alpha(j)=\alpha(j')=i$}\\
0 &\textrm{otherwise}
\end{cases}
\end{equation}
\begin{prop}
The simplicial set $N^{dg}(\mathcal{C})$ is an $\infty$-category.
\end{prop}
\noindent\emph{Proof}. Cf \cite[Proposition 1.3.1.10]{Lur2}.\qed\\\par\indent
Let $\mathcal{C}$ be a category, we denote by $\mathcal{C}_k$ its associated free $k$-linear category, viewed as a dg category over $k$ concentrated in degree $0$. The following result is well known, and we sketch a proof for completeness.
\begin{prop}
Let $\mathcal{C}$ be a category, and $\mathcal{D}$ a ($\mathbb{Z}$ or $\mathbb{Z}/2$-graded) dg category over $k$. Then there is a bijection
\begin{equation}
(-)_{\Delta}:\mathrm{Fun}^{A_{\infty},u}(\mathcal{C}_k,\mathcal{D})\xrightarrow{\cong}\mathrm{sSet}(N\mathcal{C},N^{dg}\mathcal{D}),
\end{equation}
where $\mathrm{Fun}^{A_{\infty},u}$ denotes the set of strictly unital $A_{\infty}$-functors.
\end{prop}
\noindent\emph{Proof}. Let $F: \mathcal{C}_k\rightarrow \mathcal{D}$ be a strictly unital $A_{\infty}$-functor. Recall that this is the data of \par\indent
$\bullet$ a map on objects $F: \mathrm{Ob}(\mathcal{C})\rightarrow \mathrm{Ob}(\mathcal{D})$ and\par\indent
$\bullet$ for each composable sequence of morphisms $x_0\xrightarrow{f_1} x_1\xrightarrow{f_2}\cdots\xrightarrow{f_d}x_d, d\geq 1$ in $\mathcal{C}$, an assignment $F^d(f_d,\cdots,f_1)\in \mathrm{Hom}_{\mathcal{D}}^{-d}(F(x_0),F(x_d))$
that satisfies the $A_{\infty}$-equations and unital condition.\par\indent
Given such an $F$, we now associate a map of simplicial sets $F_{\Delta}:N\mathcal{C}\rightarrow N^{dg}\mathcal{D}$. We start by describing $F_{\Delta}$ in low degrees. On $0$-simplicies,
we define $(F_{\Delta})_0:N\mathcal{C}_0\rightarrow N^{dg}\mathcal{D}_0$ to be the map $F$ on objects. On $1$-simplicies, we define
\begin{equation}
(F_{\Delta})_1(x_0\xrightarrow{f_1}x_1)_{01}:=F^1(f_1).
\end{equation}
For a $2$-simplex $x_0\xrightarrow{f_1}x_1\xrightarrow{f_2}x_2$, compatibility with face maps dictates that we define
\begin{equation}
\begin{cases}
(F_{\Delta})_2(x_0\xrightarrow{f_1}x_1\xrightarrow{f_2}x_2)_{01}:=F^1(f_1),\\
(F_{\Delta})_2(x_0\xrightarrow{f_1}x_1\xrightarrow{f_2}x_2)_{12}:=F^1(f_2),\\
(F_{\Delta})_2(x_0\xrightarrow{f_1}x_1\xrightarrow{f_2}x_2)_{02}:=F^1(f_2\circ f_1).
\end{cases}
\end{equation}
Finally, we define
\begin{equation}
(F_{\Delta})_2(x_0\xrightarrow{f_1}x_1\xrightarrow{f_2}x_2)_{012}:=F^2(f_2,f_1).
\end{equation}
Condition (D.1) in this case reads
\begin{equation}
dF^2(f_2,f_1)=-F^1(f_2\circ f_1)+F^1(f_2)\circ F^1(f^1),
\end{equation}
which is satisfied since $F$ is an $A_{\infty}$-functor. \par\indent
We can continue this and define $F_{\Delta}$ inductively. By requiring $F_{\Delta}$ to be compatible with face maps, the components of $(F_{\Delta})_d(f_d,\cdots,f_1)$ are completely determined from lower degree data except for
$(F_{\Delta})_d(f_d,\cdots,f_1)_{01\cdots d}$, which we define to be
\begin{equation}
(F_{\Delta})_d(f_d,\cdots,f_1)_{01\cdots d}:=F^d(f_d,\cdots,f_1).
\end{equation}
As before, condition (D.1) is guaranteed by $F$ being an $A_{\infty}$-functor, and by construction, $F_{\Delta}$ is compatible with face maps. Moreover, it is easy to see that
when $F$ is strictly unital, $F_{\Delta}$ is compatible with degeneracy maps. This gives the desired assignment $(-)_{\Delta}: \mathrm{Fun}^{A_{\infty},u}(\mathcal{C}_k,\mathcal{D})\rightarrow\mathrm{sSet}(N\mathcal{C},N^{dg}\mathcal{D})$.
One can easily check that $G\mapsto G_{A_{\infty}}$, where
\begin{equation}
G_{A_{\infty}}^d(f_d,\cdots,f_1):=G_d(x_0\xrightarrow{f_1}\cdots\xrightarrow{f_d}x_d)_{01\cdots d},
\end{equation}
defines an inverse to $F\mapsto F_{\Delta}$. \qed\par\indent
Let $\mathrm{Mod}_k$ denote the dg category of chain complexes over $k$, and $\mathrm{Mod}_{k[\epsilon]}$ the dg category of mixed complexes.
\begin{mydef}
In each case, let $W$ denote the collection of quasi-isomorphisms. The $\infty$-category $N^{dg}\mathrm{Mod}_k[W^{-1}]$ is called the \emph{derived category of chain complexes over $k$}, denoted $Mod_k$. The $\infty$-category $N^{dg}\mathrm{Mod}_{k[\epsilon]}[W^{-1}]$ is called the \emph{derived category of mixed complexes over $k$}, denoted $Mod_{k[\epsilon]}$.
\end{mydef}
\begin{rmk}
By \cite[Proposition 1.3.5.15]{Lur2}, $Mod_k$ is equivalent to the $\infty$-category $N\mathrm{Ch}_k[W^{-1}]$, where $W$ denotes the collection of quasi-isomorphisms in the ordinary category of chain complexes $\mathrm{Ch}_k$. Analogously, $Mod_{k[\epsilon]}$ is equivalent to $N\mathrm{Ch}_{k[\epsilon]}[W^{-1}]$.
\end{rmk}

D.2. \textbf{Kan extensions}. In this subsection, we review the notion of Kan extensions in the $\infty$-categorical setting, a relative form of $\infty$-categorical (co)limits, and some of its properties. Heuristically, given a diagram
\begin{equation}
\begin{tikzcd}[row sep=1.2cm, column sep=0.8cm]
\mathcal{C}\arrow[r,"X"]\arrow[d,"F"]&\mathcal{E}\\
\mathcal{D}\arrow[ur,dashed]&
\end{tikzcd}
\end{equation}
the left Kan extension of $X$ along $F$, denoted $F_!(X)$, is a functor from $\mathcal{D}$ to $\mathcal{E}$ that should be thought of as a base change `$X\otimes_{\mathcal{C}}\mathcal{D}$'.
Alternatively, it serves as a left adjoint to the pullback $F^*:\mathrm{Fun}(\mathcal{D},\mathcal{E})\rightarrow\mathrm{Fun}(\mathcal{C},\mathcal{E})$. The right adjoint to the pullback
is called the \emph{right Kan extension of $X$ along $F$}, denoted $F_*(X)$. In what follows, we discuss the notion of a left Kan extension; there is a dual notion of right Kan extension which we omit, as it
can be obtained from its left counterpart by taking opposite categories. \par\indent
As a preliminary, we recall the definition of a left Kan extension in the ordinary $1$-categorical setting, and two formulas to compute it.
\begin{mydef}
Let $X:\mathcal{C}\rightarrow \mathcal{E}$ (viewed as a diagram) and $F:\mathcal{C}\rightarrow \mathcal{D}$ be functors between $1$-categories. Then the \emph{left
Kan extension} of $X$ along $F$ is the data of a functor $F_!(X): \mathcal{D}\rightarrow \mathcal{E}$ together with a natural transformation $\epsilon:X\Rightarrow F_!(X)\circ F$ that is initial
among all such pairs. In other words, there is a natural bijection $\mathrm{Nat}(F_!(X),S)\cong \mathrm{Nat}(X,S\circ F)$, where $S\in\mathrm{Fun}(\mathcal{D},\mathcal{E})$.
\end{mydef}
By abuse of terminology, we often call $F_!(X)$ the left Kan extension, and make the natural transformation $\epsilon$ implicit.\par\indent
As an example, left Kan extending a functor $X:\mathcal{C}\rightarrow \mathcal{E}$ along $\mathcal{C}\rightarrow *$ exactly
recovers $\mathrm{colim}_{\mathcal{C}}X$. Therefore, left Kan extension can be viewed as a relative form of colimit. The follow lemma makes this precise
by giving a pointwise colimit formula for computing left Kan extensions.\par\indent
Before we state the lemma, we introduce a notation. Let $F:\mathcal{C}\rightarrow \mathcal{D}$ be a functor and $d\in\mathcal{D}$. Then the comma category $\mathcal{C}_{/d}$ (with $F$ implicit) has objects pairs
$(c\in\mathcal{C}, f:F(c)\rightarrow d)$, and morphisms are those $h:c\rightarrow c'$ such that $f'\circ F(h)=f$. There is a forgetful functor $U:\mathcal{C}_{/d}\rightarrow \mathcal{C}$
sending $(c,f)$ to $c$.

\begin{lemma}
Let $X:\mathcal{C}\rightarrow \mathcal{E}$ and $F:\mathcal{C}\rightarrow \mathcal{D}$ be functors between $1$-categories. Suppose for each $d\in\mathcal{D}$, the diagram
$(\mathcal{C}_{/d})\xrightarrow{U} \mathcal{C}\xrightarrow{X} \mathcal{E}$ has a colimit. Then the left Kan extension $F_!(X)$ exists, and satisfies
\begin{equation}
F_!(X)(d)\cong \mathrm{colim}_{\mathcal{C}_{/d}} X\circ U.
\end{equation}
\end{lemma}

If one thinks of colimit heuristically as integrating over a category, then left Kan extensions are the analogue of pushing forward differential forms. Lemma D.2 then says that pushforward can be
computed by `integrating over fibers'. \par\indent
Now we give another formula for left Kan extension using coends, which makes precise the heuristic that left Kan extension acts as a base change `$X\otimes_{\mathcal{C}}\mathcal{D}$'.
\begin{lemma}
Let $X:\mathcal{C}\rightarrow \mathcal{E}$ and $F:\mathcal{C}\rightarrow \mathcal{D}$ be functors between $1$-categories. Suppose for each $c,c'\in\mathcal{C}$ and $d\in\mathcal{D}$,
the copower $\mathcal{D}(F(c),d)\otimes X(c')$ exists, and furthermore that the coend $\int^{c\in \mathcal{C}} \mathcal{D}(F(c),d)\otimes X(c)$ exists for each $d\in\mathcal{D}$. Then the left Kan extension
$F_!(X)$ exists and for each $d\in\mathcal{D}$,
\begin{equation}
F_!(X)(d)\cong \int^{c\in \mathcal{C}} \mathcal{D}(F(c),d)\otimes X(c).
\end{equation}
\end{lemma}

Now we work consider the setting of $\infty$-categories. The definition of left Kan extension in this setting is
easily framed in terms of a pointwise colimit formula (cf Lemma D.2). We first consider the case of left Kan extending along an inclusion.

\begin{mydef}
Let $X:\mathcal{C}\rightarrow \mathcal{E}$ be a map of $\infty$-category, and let $\mathcal{C}^0\subset \mathcal{C}$ be a full subcategory. Let $X^0:=X|_{\mathcal{C}^0}$. We say that
$X$ is \emph{left Kan extended along} $\mathcal{C}^0\subset\mathcal{C}$ if for each $c\in\mathcal{C}$, $X(C)$ is the colimit of the diagram
\begin{equation}
\mathcal{C}^0_{/C}\rightarrow \mathcal{C}^0\xrightarrow{X^0} \mathcal{E}.
\end{equation}
\end{mydef}

By abuse of notation, we denote $X$ as $\iota_!X^0$, where $\iota:\mathcal{C}^0\subset\mathcal{C}$ is the inclusion. We refer the reader to \cite[Def 4.3.3.2]{Lur1} for
a definition of Left Kan extension along a general functor.



The universal property of left Kan extension is characterized by the following lemma (cf \cite[Prop 4.3.3.7]{Lur1}).
\begin{lemma}
Let $\delta:K\rightarrow K'$ be a map of simplicial sets, and $\mathcal{E}$ an $\infty$-category. Then $\delta_!: \mathrm{Fun}(K,\mathcal{E})\rightarrow \mathrm{Fun}(K',\mathcal{E})$
is a left adjoint to the pullback $\delta^*: \mathrm{Fun}(K',\mathcal{E})\rightarrow \mathrm{Fun}(K,\mathcal{E})$.\qed
\end{lemma}
D.3. \textbf{(Relative) cofinality}.
The question of cofinality naturally arises when one asks the following question. When is it possible to compute the colimit of a diagram $F: K\rightarrow \mathcal{E}$ using a
sub-diagram $K'\subset K$?

\begin{mydef}
A map of simplicial sets $v: K'\rightarrow K$ is \emph{cofinal} if for every $\infty$-category $\mathcal{C}$ and every colimit diagram $\overline{p}:K^{\triangleright}\rightarrow \mathcal{C}$,
the induced map $\overline{p}':{K'}^{\triangleright}\rightarrow \mathcal{C}$ is a colimit diagram.
\end{mydef}

The key theorem regarding cofinality is the following criterion, also known as Quillen's theorem A (cf \cite[Theorem 4.1.3.1]{Lur1}) for $\infty$-categories.

\begin{thm}
Let $f:\mathcal{C}\rightarrow \mathcal{D}$ be a map of simplicial sets, where $\mathcal{D}$ is an $\infty$-category. Then $f$ is cofinal if and only if for each $d\in\mathcal{D}$, the
simplicial set $\mathcal{C}_{d/}$ is weakly contractible. \qed
\end{thm}
Given a commutative square of $\infty$-categories
\begin{equation}
\begin{tikzcd}[row sep=1.2cm, column sep=0.8cm]
\mathcal{A}\arrow[r,"f"]\arrow[d,"u"]&\mathcal{B}\arrow[d,"v"]\\
\mathcal{C}\arrow[r,"g"]&\mathcal{D}
\end{tikzcd},
\end{equation}
there is a natural transformation $u_!f^*\rightarrow g^*v_!$, called the \emph{Beck-Chevalley transform}, induced by the adjunction units and counits
\begin{equation}
u_!f^*\rightarrow u_!f^*v^*v_!\xrightarrow{\simeq} u_!u^*g^*v_!\rightarrow g^*v_!.
\end{equation}
There is also a dual Beck-Chevalley transform $u_*f^*\leftarrow g^*v_*$.
Just as left Kan extension is a relative notion of colimit, there is a relative notion of cofinality.
\begin{mydef}
Using notation as in (D.15). The commutative square of $\infty$-categories is a \emph{homotopy exact square} if, assuming left Kan extensions exist, for each $X: \mathcal{B}\rightarrow \mathcal{E}$, the Beck-Chevalley transform $u_!f^*X\rightarrow g^*v_!X$ is an equivalence in $\mathcal{E}$.
\end{mydef}
If one sets $\mathcal{C}=\mathcal{D}=*$, this recovers the notion of cofinality for $f$.
The following theorem is a straightforward generalization of Quillen's theorem A.
\begin{thm}
With notation as above, the commutative square is homotopy exact if and only if for each $b\in \mathcal{B}, c\in\mathcal{C}$ and morphism $\varphi: v(b)\rightarrow g(c)$, the double
comma category $({}_{b/}A_{/c})_{\varphi}$ is weakly contractible.
\end{thm}

\end{appendices}

\end{document}